\documentclass[11pt]{article}

\usepackage{amsmath,amssymb,amsthm}
\usepackage{graphicx}
\usepackage{hyperref}
\usepackage{authblk}

\usepackage{t1enc}

\usepackage{cancel}

\usepackage{soul}

\usepackage{subcaption}

\usepackage{tikz}

\usetikzlibrary{arrows.meta,positioning,fit,calc}

\makeatletter
\DeclareRobustCommand{\gobblesomeargs}[2]{#2}
\renewcommand{\p@equation}{\gobblesomeargs}
\makeatother

\usepackage{lipsum}
\usepackage{caption}
\usepackage{subcaption}
\usepackage[all,dvips]{xy}
\usepackage{color}
\usepackage{tikz-cd}
 
\usetikzlibrary{arrows.meta,positioning,calc}
\usepackage{xcolor}

\usepackage{yhmath}

\usepackage{wasysym}

\usepackage{harmony} 

\usepackage{pifont}

\usepackage{float}

\usepackage{dsfont}\let\mathbb\mathds

\usepackage{fdsymbol}

\usepackage{hyperref}

\usepackage{bbm}

\def\virgp{\raise 2pt\hbox{,}}

\renewcommand{\geq}{\geqslant}
\renewcommand{\leq}{\leqslant}
\def\N{{\mathbb N}}

\def\R{{\mathbb R}}

\def\virgp{\raise 2pt\hbox{,}}
\def\cdotpv{\raise 2pt\hbox{;}}

\def\1{\mathbbm{1}}

\newtheorem{theorem}{Theorem}[section]
\newtheorem{corollary}[theorem]{Corollary}
\newtheorem{proposition}[theorem]{Proposition}

\theoremstyle{remark}
\newtheorem{remark}{Remark}[section]

\theoremstyle{definition}
\newtheorem{definition}{Definition}[section]

\theoremstyle{definition}

\theoremstyle{definition}

\DeclareSymbolFont{yhlargesymbols}{OMX}{yhex}{m}{n}
\DeclareMathAccent{\wideparen}{\mathord}{yhlargesymbols}{"F3}

\usepackage{authblk}

\title{Emergent Loewner Dynamics in Slime Mold Growth}
\date{}

\author[,1]{\textbf{Claire David}\thanks{Corresponding author: \texttt{Claire.David@sorbonne-universite.fr}}}

\author[1]{\textbf{Aur\`ele Boussard}}

\author[2]{\textbf{Nizare Riane} }

\author[3]{\textbf{\mbox{Michel L.~Lapidus}}}

\author[4]{\textbf{Audrey Dussutour}}

\affil[1]{Sorbonne Universit\'e, Universit\'e Paris Cit\'e, CNRS, INRIA, Laboratoire Jacques-Louis Lions, LJLL, F-75005 Paris, France}

\affil[2]{Universit\'e Mohamed V - Rabat, Morocco}
\affil[3]{University of California, Riverside, Department of Mathematics, Skye Hall, 900 University Ave., Riverside, CA 92521-0135,
	USA}

\affil[4]{Universit\'e Paul Sabatier - CNRS, Centre de Recherches sur la Cognition Animale, CRCA,
	118 route de Narbonne
	31062 TOULOUSE,  France}

\title{\textbf{Emergent Loewner Dynamics in Slime Mold Growth 
}}

\begin{document}

	\maketitle

	\thanks{This research is supported by the MITI CNRS 80PRIME, Conditions Extrêmes, Formes Optimales and by the Agence Nationale de la Recherche (ANR) FRACTALS (ANR-24-CE45-3362).}
	
	\vskip 1cm

	\begin{abstract}

		Growth fronts of slime molds are characterized through a direct geometric analysis based on Loewner evolutions, using experimentally acquired time-resolved images. The associated Loewner driving functions reconstructed from expanding pseudopod boundaries display statistical properties consistent with Gaussian-like behavior.\\
		
		A geometric estimate of the diffusivity parameter~$\kappa$ is inferred from fractal scaling, while Brownian diagnostics are assessed on the reconstructed driving signal.\\
		
		These findings show that the boundaries of a growing living organism display statistical and geometric properties consistent with emergent Loewner dynamics over experimentally accessible scales. This study establishes a quantitative framework for analyzing biological growth interfaces and suggests new connections between morphogenesis, stochastic geometry, and network reorganization under varying environmental conditions.
		We provide, to our knowledge, the first explicit reconstruction of a Loewner driving function from a living growth interface, revealing an emergent Brownian-like conformal growth regime at expanding fronts.

	\end{abstract}

	 \vskip 1cm


 	\newpage	 
	 
	\tableofcontents
	
 	\newpage

	\section{Introduction}

	\hskip 0.5cm SLE -- Schramm--Loewner Evolution -- is usually associated with statistical physics, in connection with critical phenomena~\cite{OdedSchrammSteffenRhodBAsicPropertiesOfSLE2005}. In the classical setting, the driving function is Brownian motion (up to a diffusivity parameter~$\kappa$). The traces of the SLE process~\cite{OdedSchrammScalingLimitsOfLoopErasedRandomWalksAndUniformSpanningTrees2000}, consist of curves -- the so-called~\mbox{SLE$_\kappa$}. As is recalled in~\cite{BertrandDuplantierFractalCriticalPhenomenaInTwoDimensionsAndConformalInvariance1989},~\cite{BertrandDuplantierIlliaBinderHarmonicMeasureAndWindingOfRandomConformalPathsACoulombGasPerspective2008}, those (conformally invariant) fractal curves are associated with critical phenomena which arise in statistical systems (for instance, critical Ising or Potts percolation clusters). They arise at \emph{critical phase transitions}.

	In real life, morphological growth may combine stochastic and deterministic components. However, an explicit identification of Loewner-type dynamics in living systems remains limited. For instance, in the case of collective cellular movement, in~\cite{BenjaminAndersenEtAlEvidenceOfRobustUniversalConformalInvarianceInLivingBiologicalMatter2024}, Benjamin~H.~Andersen and his collaborators report experimental evidence that the flows generated by some collective cells not only exhibit conformal invariance, but are also described by the Schramm--Loewner Equation (SLE).  Yet,  their study did not involve an explicit reconstruction of the Loewner driving function.
	
	Previously, an attempt was made to identify a~\mbox{SLE} in turbulence by Denis~Bernard, Guido~Boffetta, Antonio Celani and Gregory Falkovich; see~\cite{BernardBoffettaCelaniFalkovichConformalInvarianceInTwodimensionalTurbulence2006}. However, they did not verify the defining stochastic properties of a Brownian Loewner driving function. 
	It was later shown that such criteria were insufficient to recover an underlying SLE. See the paper by Tom Kennedy~\cite{TomKennedyComputingTheLoewnerDrivingProcessOfRandomCurvesInTheHalfPlane2008}. This illustrates that geometric indicators alone may mimic~\mbox{SLE$_\kappa$} geometry while failing to recover the underlying Brownian dynamics of the driving function.

	Our unicellular model organism, the slime mold \emph{Physarum polycephalum}, forms intricate, highly ramified networks exhibiting pronounced multiscale spatial organization, often described as fractal-like; see~\cite{PhilipRosinaMartinGrubeNovelImageAnalyticApproachRevealsNewInsightsInFineTuningOfSlimeMouldNetworkAdaptation2024},~\cite{PhilipRosinaMartinGrubeAMathematicalModelToPredictNetworkGrowthInPhysarumPolycephalumAsAFunctionOfExtracellularMatrixViscosityMeasuredByANovelViscometer2025}. Its motion and behavior rely on a transport network composed of interconnected veins. The network is driven by synchronized, rhythmic actomyosin contractions which generate pressure gradients and induce a bidirectional cytoplasmic flow (\emph{shuttle-streaming}); see~\cite{AudreyDussutourChloeArsonFlowNetworkAdaptationAndBehaviorInSlimeMolds2024}, along with  the review~\cite{MathieuLeVergeSerandourAndKarenAlimPhysarumPolycephalumSmartnetworkAdaptation2023}. Those contractions regulate both nutrient transport and large-scale cellular motility. Veins continuously appear and disappear in response to local pressure and shear, while pseudopods extend and retract (in coordination with the streaming flow). The resulting evolving growth interface therefore reflects an adaptive, internally driven, oscillatory, and spatially distributed dynamical process. This naturally raises the question of whether the growth interface is governed by an underlying stochastic law.
	  
	  Percolation properties of such networks have been reported~\cite{AdrianFesselChristinaOttmeyerHGDobereinerStructuringPrecedesExtensionInPercolatingPhysarumPolycephalumNetworks2015}. While percolation models are known to exhibit deep connections with conformally invariant scaling limits, the stochastic structure of the evolving growth interface itself has not been explicitly reconstructed in this biological context.

	To address this question, we used time-resolved experimental imaging to extract the evolving growth fronts.   
 From these data, we reconstruct the Loewner driving function, numerically solve the Loewner equation, and estimate a geometric diffusivity. This yields quantitative evidence for an emergent conformal-growth description of slime-mold fronts, and reveals a systematic multilevel modulation of the effective driving statistics by internal structural constraints.

	\section{Experimental Data and Growth Front Extraction}
 	\label{ExperimentalDataAndGrowthFrontExtraction}

	Our identification was carried out on the following representative experiment: a slime mold (\emph{Physarum polycephalum}) of the AUS strain, with an initial diameter of~13~mm, explored a 90-mm-diameter petri dish containing 1\% agar gel for 24~hours. Both rearing (on moistened oat flakes, Quaker Oats Company) and experiments were conducted at 25°C in temperature-controlled incubators. Images were recorded every two minutes using digital cameras (EOS 70D, Canon). 
	The resulting data were analyzed by means of the Python-based software \href{https://github.com/Aurele-B/Cellects}{Cellects}; see~\cite{CellectsASoftwareToQuantifyCellExpansionAndMotion}.

	An important comment to be made is that the purpose of this work is inferential: we test whether experimental fronts admit a Loewner description and whether the inferred driver is consistent with an effective Brownian driving function. 
	

	\begin{figure}
		\centering
		
		\begin{subfigure}{0.32\linewidth}
			\includegraphics[width=\linewidth]{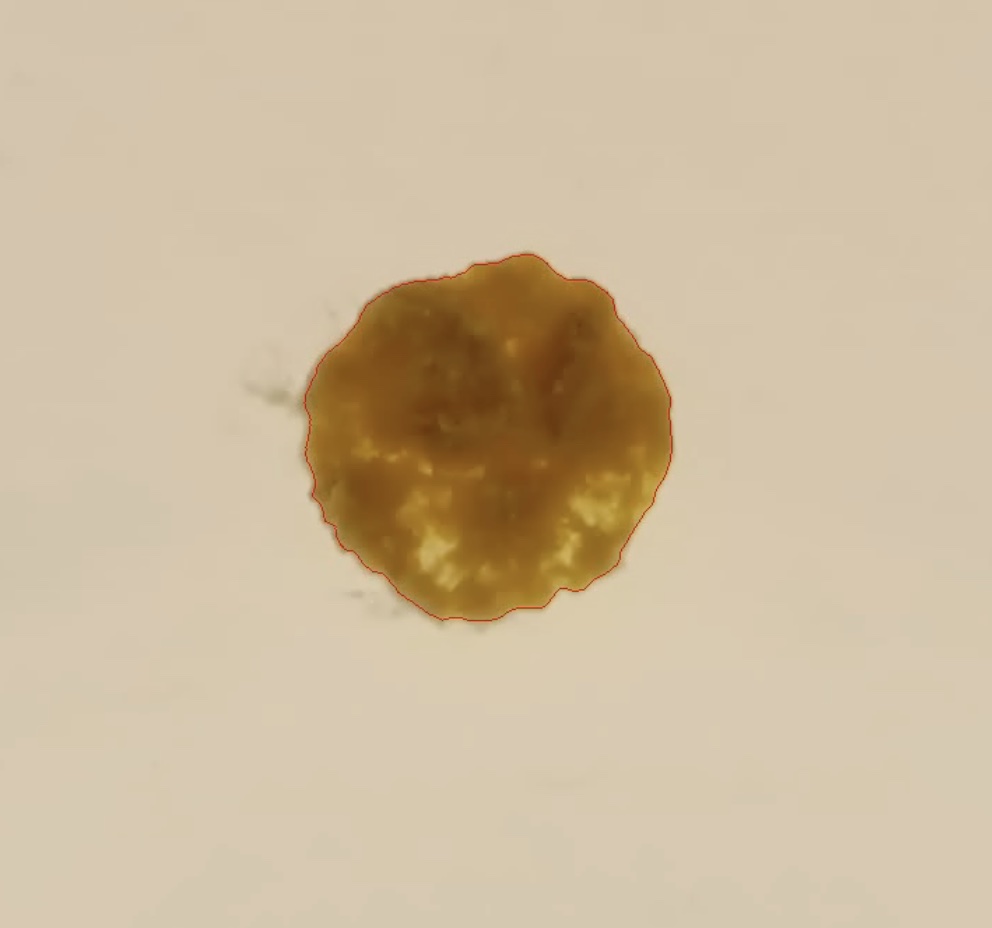}
			\caption{The slime mold at time~\mbox{$t = 0{:}00$}.}
		\end{subfigure}
		\hfill
		\begin{subfigure}{0.32\linewidth}
			\includegraphics[width=\linewidth]{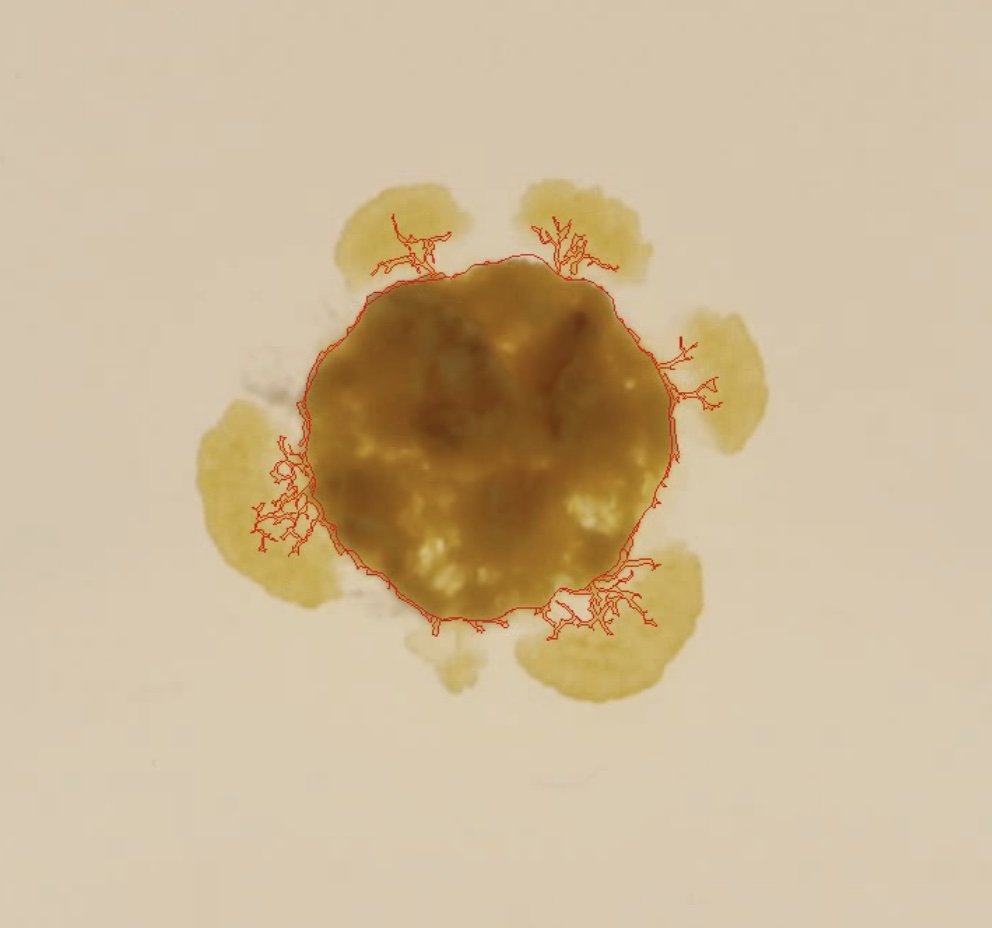}
			\caption{The slime mold at time~\mbox{$t \approx 3{:}20$}.}
		\end{subfigure}
		\hfill
		\begin{subfigure}{0.32\linewidth}
			\includegraphics[width=\linewidth]{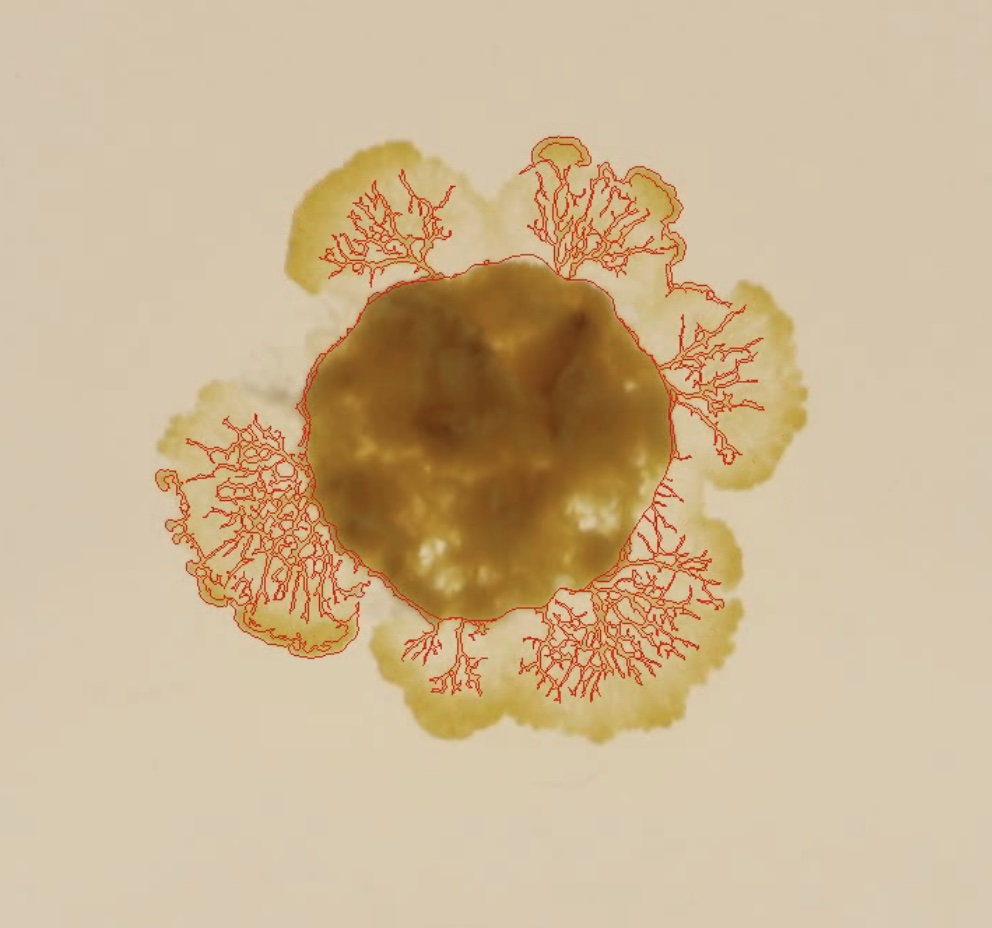}
			\caption{The slime mold at time~\mbox{$t \approx 4{:}40$}.}
		\end{subfigure}
		
		\vspace{0.2cm}
		
		\begin{subfigure}{0.32\linewidth}
			\includegraphics[width=\linewidth]{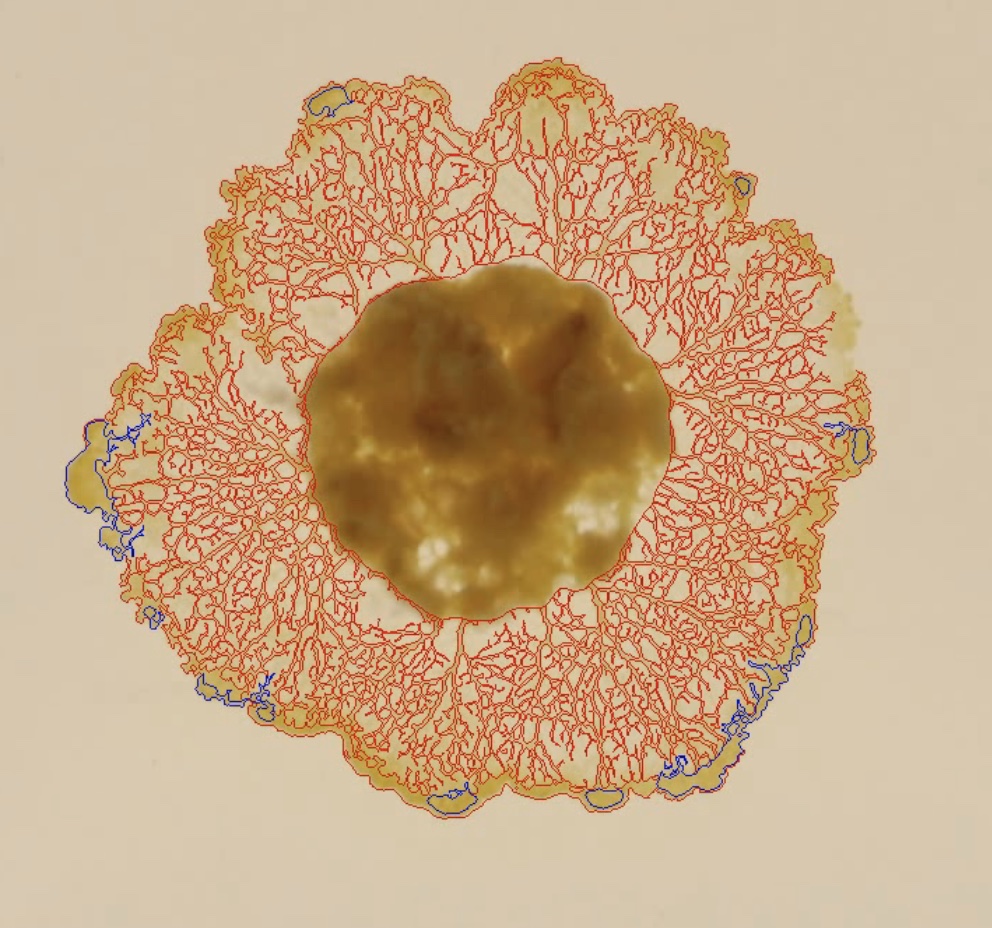}
			\caption{The slime mold at time~\mbox{$t \approx 6{:}40$}.}
		\end{subfigure}
		\hfill
		\begin{subfigure}{0.32\linewidth}
			\includegraphics[width=\linewidth]{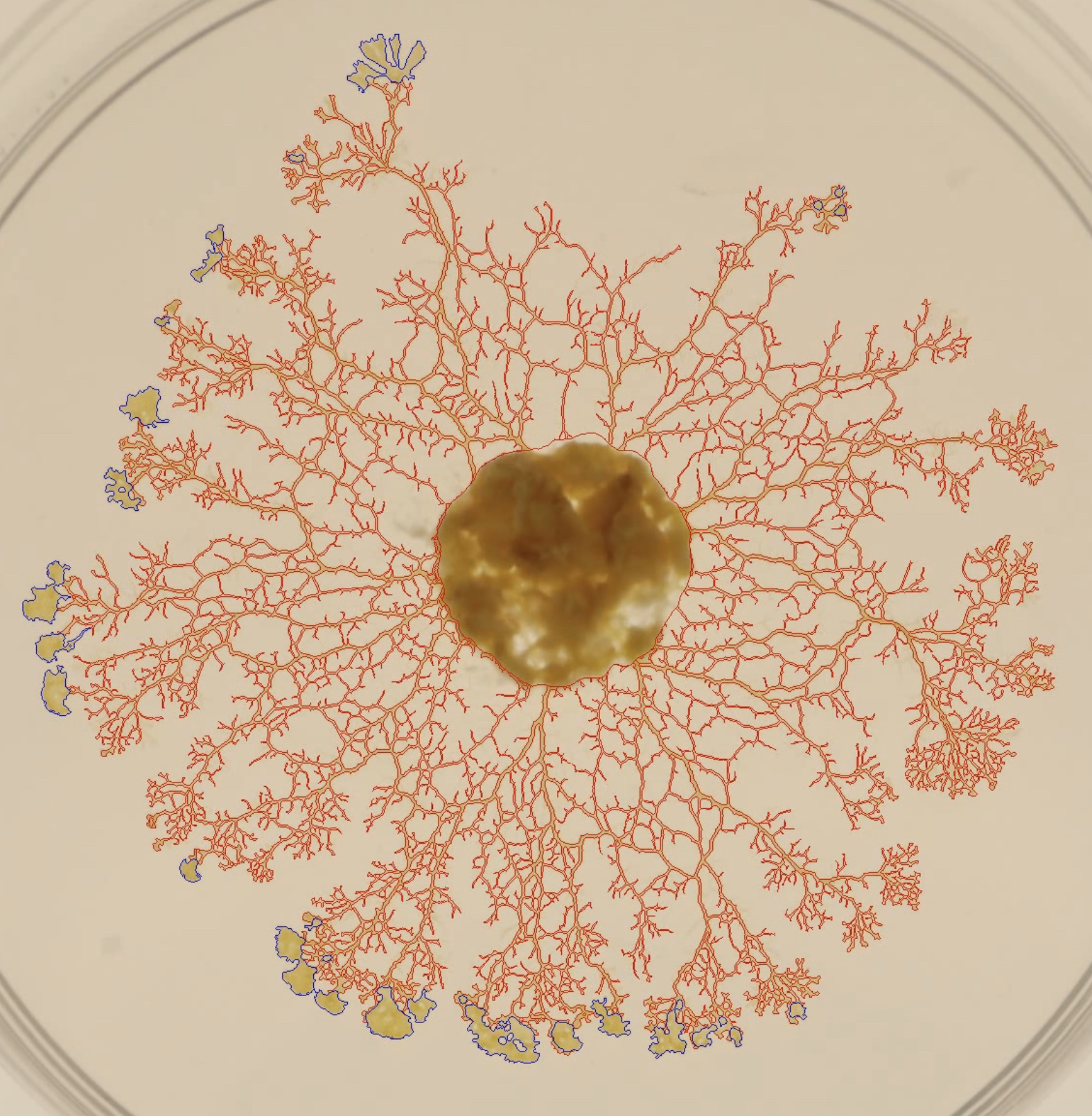}
			\caption{The slime mold at time~\mbox{$t \approx 13{:}20$}.}
		\end{subfigure}
		\hfill
		\begin{subfigure}{0.32\linewidth}
			\includegraphics[width=\linewidth]{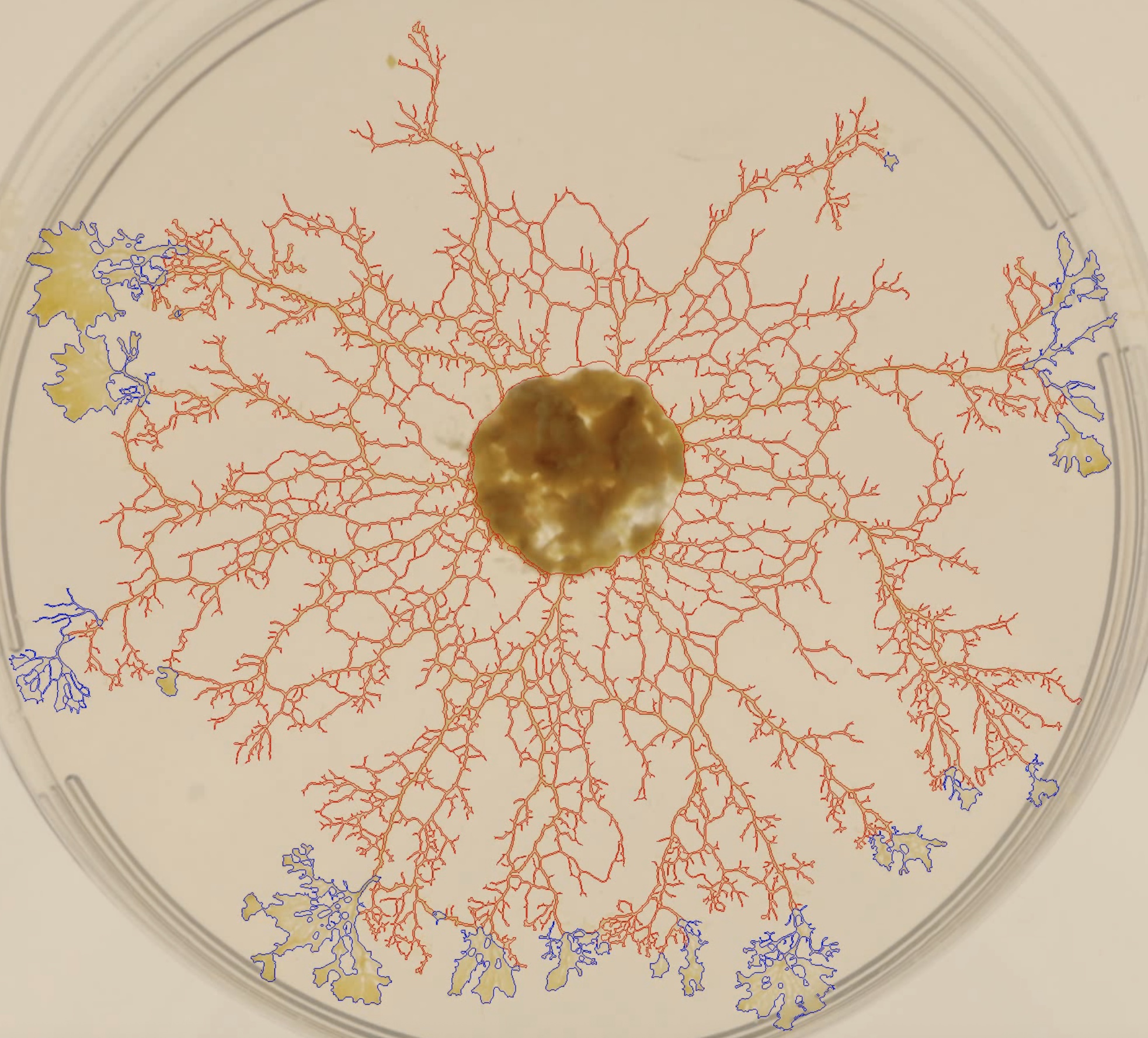}
			\caption{The slime mold at time~\mbox{$t \approx 17{:}20$}.}
		\end{subfigure}

		\caption{\textbf{Time-compressed growth sequence of \emph{Physarum polycephalum}.}
			Snapshots extracted from a 36-s video representing a 24-h experiment (1 frame every 2 min). Times indicate the corresponding experimental time.}

		\label{Snapshots}
	\end{figure}

	\section{Reconstruction of the Loewner Driving Function from Experimental Fronts}
	\label{SectionReconstructionOfTheLoewnerDrivingFunctionFromExperimentalFronts}

	We divide the time-lapse video into several outer and inner windows, defined on a regular mesh, in order to follow in detail the local growth dynamics.
	The objective is to determine whether the (discrete) time-parametrized curve~\mbox{$t \mapsto C_t$} corresponding to the boundary of the spatial region under study can be described by a Loewner-type growth process, and to reconstruct (in an inverse sense) the associated Loewner driving function~\mbox{$U_t$}. To this end, we consider the chordal Loewner equation
	in the upper half-plane,
	\[\label{LoewnerEquation}
	\partial_t g_t(z) = \frac{2}{g_t(z) - U_t},
	\qquad g_0(z) = z,
	\]
	where~\mbox{$U_t$} denotes the (a priori unknown) driving function.
	The central question is therefore whether the experimentally reconstructed observable~\mbox{$U_t$} 
	exhibits the defining stochastic properties of a Brownian motion with variance parameter~$\kappa$. In the ideal Brownian setting this would correspond to~\mbox{SLE$_\kappa$}; note that, in the present biological context, we consider~$\kappa$ as a geometric diffusivity and empirically test Brownian diagnostics.
	
	The outer windows correspond to growing pseudopods, where an active boundary can be directly extracted from the video frames. The exterior boundary~\mbox{$\partial {\cal S}_t = C_t$} of the occupied region provides a natural activity mask, from which we extract the associated (discrete) time-parametrized curve~\mbox{$t \mapsto C_t$} describing the boundary of the expanding fronts.

	For the inner windows, there is no propagating boundary in the usual sense. We bypass this configuration by introducing \emph{an effective surrogate driving observable} in the following way:  first, we determine the main direction of activity by relying on a principal component analysis (PCA) applied to the set of active pixels within the considered spatial window. More precisely, for each window, we compute the covariance matrix of the spatial coordinates of all pixels belonging to the occupied region. The aggregation is performed over the entire duration of the experiment. The first principal component, corresponding to the largest eigenvalue of this covariance matrix, defines a unit vector that captures the dominant direction of growth or oscillatory transport within the window. Hence, projecting the boundary displacements onto this principal axis provides a robust scalar observable~\mbox{$U_t$}, which retains the dominant geometric component of the local dynamics while filtering out transverse fluctuations.
	
	In both configurations, the resulting scalar signal is interpreted as an effective Loewner driving observable.
	
	In practice, we denote by~\mbox{$\Omega_W=  \left ( x_{W}, y_W\right )$} the center of mass of the considered inner window. We then compute, for each time frame, the minimal active radius
	
	$$R_t= \displaystyle \min_{N_t \, \in\, \mathrm{Window} } \left \| N_t- \Omega_W \right  \| =  \left \| M_t- \Omega_W \right  \| 
	 \, .$$

	\noindent where~\mbox{$ M_t  $} denotes the point where the minimum is reached. The projection~of the vector~\mbox{$\displaystyle   \overrightarrow{ \Omega_WM_t} $} onto the main direction (i.e., after rotating the coordinates in order to have this direction become the real axis) yields a scalar radial displacement~\mbox{$ U_t$}.

	\section{Multilevel Loewner Reconstruction}

	\label{SectionMultilevelAnalysis}

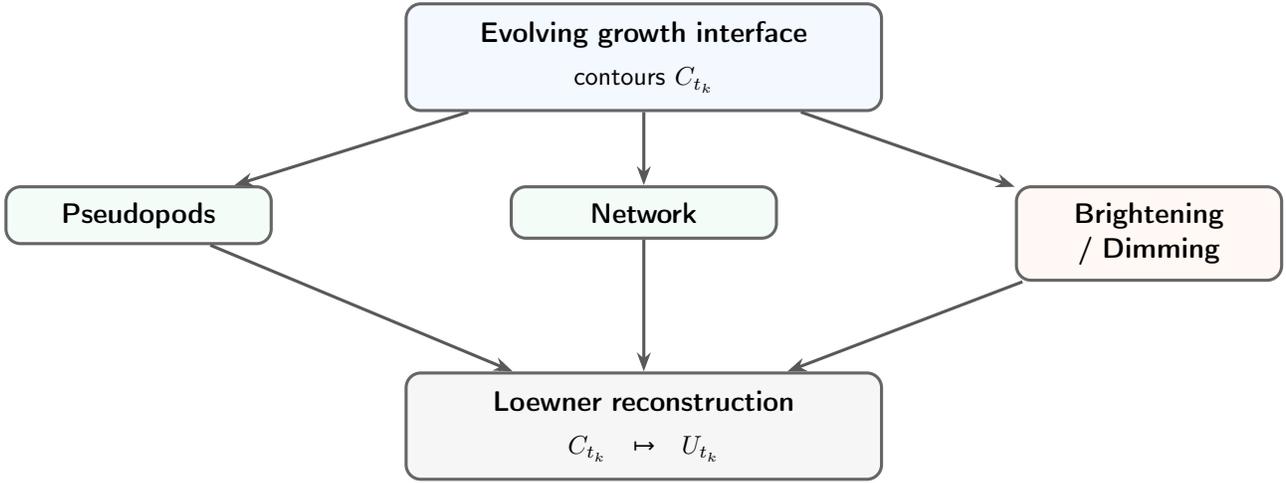
\begin{figure}[t]
 	\centering
 	\resizebox{\columnwidth}{!}{%
 		\begin{tikzpicture}[
 			font=\sffamily,
 			bigbox/.style={
 				rounded corners=2mm,
 				draw=black!60,
 				very thick,
 				fill=#1,
 				inner sep=7pt,
 				align=center,
 				text width=6.0cm
 			},
 			smallbox/.style={
 				rounded corners=2mm,
 				draw=black!60,
 				very thick,
 				fill=#1,
 				inner sep=6pt,
 				align=center,
 				text width=3.2cm
 			},
 			arrow/.style={-{Stealth[length=2.6mm,width=1.9mm]}, very thick, draw=black!65}
 			]
 			
 			\definecolor{cTop}{RGB}{244,248,255}
 			\definecolor{cMid}{RGB}{244,252,248}
 			\definecolor{cDyn}{RGB}{255,248,244}
 			\definecolor{cBottom}{RGB}{246,246,246}
 			
 			\node[bigbox=cTop] (front)
 			{\textbf{Evolving growth interface}\\[1mm]
 				{\small contours $C_{t_k}$}};
 			
 			\node[smallbox=cMid, below left=10mm and 18mm of front] (pseudo)
 			{\textbf{Pseudopods}};
 			
 			\node[smallbox=cMid, below=10mm of front] (network)
 			{\textbf{Network}};
 			
 			\node[smallbox=cDyn, below right=10mm and 18mm of front] (bd)
 			{\textbf{Brightening / Dimming}};
 			
 			\node[bigbox=cBottom, below=18mm of network] (loew)
 			{\textbf{Loewner reconstruction}\\[1mm]
 				{\small $C_{t_k}\mapsto U_{t_k}$}};
 			
 			\draw[arrow] (front) -- (pseudo);
 			\draw[arrow] (front) -- (network);
 			\draw[arrow] (front) -- (bd);
 			
 			\draw[arrow] (pseudo) -- (loew);
 			\draw[arrow] (network) -- (loew);
 			\draw[arrow] (bd) -- (loew);
 			
 		\end{tikzpicture}%
 	}
 	\caption{\textbf{Multilevel Loewner reconstruction.}
 		Time-resolved imaging yields growth contours $C_{t_k}$. The analysis is performed on three geometries: pseudopods, the extracted network, and intensity-based regions (brightening/dimming). Each configuration is mapped to a Loewner driving signal $U_t$ for statistical diagnostics.}
 	\label{DiagrammeMultiLevelAnalysis}
 \end{figure}

	In order to assess the validity of our results and to enable comparison across different geometric configurations, our Loewner analysis is performed on the following geometries, which capture distinct structural levels of organization, ranging from active expansion fronts to internal transport architecture (see Figure~\ref{DiagrammeMultiLevelAnalysis}):

	\begin{enumerate}
		\item[$i$.] The largest connected component of the pseudopods. From the biological point of view, the pseudopods correspond to the main expanding regions of our model organism, which makes them natural candidates for driving-function identification.

		\item[$ii$.] The largest connected component of the extracted network.

		\item[$iii$.] The largest connected component of the brightening regions, defined as areas exhibiting synchronous increases in pixel intensity between successive frames.

		\item[$iv$.] The largest connected component of the dimming regions, defined as areas exhibiting synchronous decreases in pixel intensity between successive frames.

	\end{enumerate}
	
	While the extracted network represents the global transport architecture, brightening and dimming regions encode internally generated dynamical activity patterns rather than geometric displacements. Together, these four configurations provide complementary structural and dynamical projections of the same biological system. 
	This multilevel hierarchy enables us to test whether Loewner-type statistics are restricted to actively expanding interfaces or instead emerge more generally from the internally regulated growth dynamics of the organism.

 
	\section{Statistical Diagnostics of the Reconstructed Driving Signal}

	We now describe the statistical diagnostics used to assess the compatibility of the reconstructed driving signal with Brownian-type Loewner dynamics.
	
\subsection{Data Extraction}

  At a given time~$t$, each window is binarized (occupied region~\mbox{${\cal S}_t$} converted to~$1$, empty region~\mbox{${\cal V}_t$}  to~$0$).
		For a given time discretization~\mbox{$  \displaystyle \underset{  k=1}{\overset{N-1  }{\bigcup}}\,   \left [t_k,t_{k+1} \right]$ } of the initial time interval~\mbox{$[0,T]$}, with~\mbox{$1 \leq N \leq 720$}, we sample the video by selecting nonredundant frames (e.g., by subsampling to reduce near-duplicate segmentations).
	We rely on the segmentation provided by \emph{Cellects}, which provides the spatial coordinates of the detected biological structures. These geometrically defined sets are then used to reconstruct the evolving interface and compute the associated Loewner driving function.
	From the segmented masks provided by \emph{Cellects}, we retain the largest connected component and extract its external boundary, which defines the contour~\mbox{$ C_{t_k}$} used for the Loewner reconstruction.

 	\subsection{Interface Reconstruction}

  For~\mbox{$  1 \leq k\leq N-1$}, we compute
		
		$$\Delta S_k= {\cal S}_{t_{k+1} }- {\cal S}_{t_{k }}\, ,$$
		
	\noindent	along with its barycenter~\mbox{$G_k $}. The difference is understood in the sense of pixel-wise set difference.
			Note that each~\mbox{$ \Delta S_k$}, with~\mbox{$  1 \leq k \leq N-1$}, corresponds to the growth of the expanding front during the time interval~\mbox{$     \left [t_k,t_{k+1} \right]$ }, while~\mbox{$G_k$} represents the attachment point of the growing front.
		The Python script then registers, for each~\mbox{$ k \, \in\, \left \lbrace 1, \ldots, N -1\right \rbrace$}:

		\begin{enumerate}
			\item[$\rightsquigarrow$]  The time stamps, which enable us to convert a frame index into a physical time.
			
			\item[$\rightsquigarrow$]  The Cartesian coordinates~\mbox{$\left (X_k,Y_k\right)$} of the barycenter~\mbox{$G_k $}.
			
			\item[$\rightsquigarrow$] If necessary, size metrics of the growing region.\\
		\end{enumerate}

		\subsection{Reconstruction of the Loewner Driving Function}
		
  Recall that in the Loewner equation~\mbox{(\ref{LoewnerEquation})}, the driving function is a real-valued boundary process. The scalar time series~\mbox{$\left \lbrace {U}_k \right\rbrace _{1\leq k \leq N} $} obtained by projection therefore constitutes a direct numerical reconstruction of the Loewner driving signal associated with the local growth dynamics.

\subsection{Geometric Diffusivity and Brownian Diagnostics} In order to obtain a biological counterpart of the usual~\mbox{SLE$_\kappa$} theory, we use $\kappa$ as a descriptive parameter inferred from geometric fractal scaling. 
Importantly, in the present work,~$\kappa$ is not identified from variance growth in time (video time need not coincide with capacity time), but is inferred independently from geometric fractal scaling (see below).
		We then assess whether the reconstructed driving signal exhibits features compatible with a Loewner-type growth process, by testing the following four criteria, both globally (i.e., for the whole considered geometry), and locally (i.e., per window):

		\begin{enumerate}

			\item[$\rightsquigarrow$] We test the Gaussianity of the driving function by plotting the associated  \emph{Quantile--Quantile} plot (in short, \emph{Q--Q plot}). Points located close to the~\mbox{$45$}-degree reference line indicate that the associated distribution is the usual normal distribution, as expected for a Brownian driving function.

			\item[$\rightsquigarrow$] We examine whether the increments of the driving signal are stationary in time. While strict independence of increments is a defining feature of a Wiener process (see SI Appendix, Section~1), we stress the fact  that this condition is not imposed \emph{a priori} here, but is tested empirically as a diagnostic indicator.

			\item[$\rightsquigarrow$] We analyze the growth of the variance of increments as a function of time. In the Brownian case, this growth is linear (see SI Appendix, Section~1); deviations from linearity are interpreted as signatures of non-Brownian stochastic dynamics or of temporal correlations in the driving process. In video time, linear growth is only indicative and does not provide a direct estimate of~$\kappa$.
			
			\item[$\rightsquigarrow$] 	We compute the power spectral density of the reconstructed driving signal. For a Brownian driving function (not its increments) the power spectrum scales as~\mbox{$\omega^{-2}$} (see SI Appendix, Section~2); more generally, departures from this scaling exponent provide a quantitative measure of non-Brownian behavior.\\
		\end{enumerate}

	We proceed as follows:

		\begin{enumerate}
			\item[$\rightsquigarrow$] We compute the spectral densities~\mbox{$\left \lbrace  P_k \right\rbrace _{1\leq k\leq N-1} $} by means of the Welch method; see SI Appendix, Section~2.

			\item[$\rightsquigarrow$]  We plot the (discrete) parametrized curve~\mbox{$\left ( \log \omega_k , \log P_k  \right)_{1\leq k \leq N-1} $}, and a robust linear regression, in order to obtain the corresponding, adjusted (linear) slope.

			\item[$\rightsquigarrow$] 
			The resulting spectral law is then written in the following form,
			
			$$ \ \ln \left ( P ( \omega) \right)= -\beta \, \ln \omega + C \,,$$

			\noindent where~\mbox{$C \, \in\, \R$} denotes a constant.
			
		\end{enumerate}

		If the estimated spectral exponent~$\beta$ is close to~\mbox{$2$} (corresponding to a value close to~\mbox{$-2$} of the adjusted slope associated with the log--log representation of the power spectrum (see SI Appendix, Section~2), the driving function provides evidence of compatibility with Brownian-type SLE dynamics. Systematic deviations from~\mbox{$\beta= 2$} instead indicate non-Brownian stochastic behavior, reflecting correlated fluctuations, intermittent dynamics, or locally regulated growth mechanisms. \\

		Finally, we note that the video acquisition rate only rescales the frequency axis and the amplitude of the power spectrum, and therefore has no influence on the estimated spectral exponent~$\beta$.

 \subsection{Tail and Correlation Diagnostics}

Further analysis is then required in order to determine whether the observed SLE-like boundary can genuinely be identified with a classical~\mbox{SLE$_\kappa$} trace, or whether it instead belongs to a broader class of non-Brownian Loewner-type growth processes.\\

		For this purpose, we proceed as follows:

		\begin{enumerate}
			\item[$\rightsquigarrow$] We compute the autocorrelation function of the effective increments 
			sequence

		\[ \left \lbrace\Delta \widetilde{U}_{k} \right\rbrace _{ 1\leq k \leq N-1} =\left \lbrace\widetilde{U}_{k+1}- \widetilde{U}_{k }\right\rbrace _{ 1\leq k \leq N-1} \, . \]

			For a Brownian driving function, this autocorrelation is expected to be close to~0 for nonzero lags, reflecting the absence of temporal correlation.

			\item[$\rightsquigarrow$] 	We examine the distribution of the effective increments sequence, and analyze the behavior of their tails using the Hill estimator; see SI Appendix, Section~2.\\

			For a Brownian driving function, the increment distribution is Gaussian and therefore exhibits rapidly decaying tails, corresponding to an effective tail exponent compatible with the domain of attraction of the normal distribution (see SI Appendix, Section~1). 
			
			In contrast, non-Brownian driving processes, including L\'evy-like or intermittently driven dynamics, are characterized by heavy-tailed increment distributions, corresponding to effective tail exponents smaller than~$2$.

			\item[$\rightsquigarrow$] If the above Gaussian assessments of the driving function are satisfied, the next step is to discuss the diffusivity parameter~$\kappa$. Since the experimental recordings provide a physical (video) time, which does not necessarily coincide with the half-plane capacity parameterization of chordal Loewner theory (see SI Appendix, Section~1), the classical linear growth relation~\mbox{$\mathrm{Var}(U_t)=\kappa\,  t$}, which strictly holds in half-plane capacity time, cannot be directly used to identify the SLE parameter~$\kappa$ without an explicit conformal reparametrization.
			Instead, the diffusivity parameter~$\kappa$ is inferred from the relation~\mbox{$D = 1 + \frac{\kappa}{8}$} (valid for~\mbox{$\kappa \le 8$}; see~\cite{VincentBeffaraDimensionOfSLECurves2008}) where the fractal dimension~$D$ is estimated by box-counting.  We obtain~\mbox{$ \kappa= 8\, (D-1)$}.
		This estimate is purely geometric and does not rely on temporal parameterization.

		\end{enumerate}

	\subsection{Geometric Consistency and Univalence}

	An important comment must be made concerning the Loewner univalence of the conformal maps associated with the reconstructed Loewner-type traces, namely the fact that the corresponding curves are simple curves, in connection with the injectivity of the solution of the Loewner equation. In the classical Brownian~\mbox{SLE$_\kappa$} setting, univalence of the Loewner flow is guaranteed for~$\kappa<8$, and the trace is a simple or non-self-crossing curve.

	However, as emphasized by Tom Kennedy~\cite{TomKennedyNumericalComputationsForTheSchrammLoewnerEvolution2009}, the inverse problem consisting in testing whether a given curve is an~\mbox{SLE} trace cannot realistically rely on a direct numerical reconstruction of the conformal maps themselves. Although, in principle, such a reconstruction could be attempted via the characterization of ~\emph{Loewner chains} as \emph{continuously increasing hulls} (see~\cite{GregoryFLawlerConformallyInvariantProcessesInThePlane2005} ~Section~\mbox{$4$}, on page~\mbox{$79$}), it is extremely sensitive to the so-called \emph{lattice effects}, coming from the discretization, and, therefore, rarely performed in practice. These effects include finite spatial resolution, pixel-based discretization, and segmentation artefacts, and they are known to strongly bias numerical conformal maps; see~\cite{TomKennedySimulatingSelfAvoidingWalksInBoundedDomains2013,TomKennedyGregLawlerLatticeEffectsInTheScalingLimitOfTheTwoDimensionalSelfAvoidingWalk2013}.

	Therefore, numerical investigations of Loewner-type dynamics typically focus on the reconstruction and statistical analysis of the driving function rather than on a direct verification of conformal univalence. This strategy is precisely the one adopted in this paper.
	
	In our present setting, the experimentally reconstructed curves are simple and monotone at the spatial and temporal resolutions considered. Thereby, this geometric property is consistent with the absence of self-intersections, a characteristic which supports the applicability of a Schramm--Loewner description of the growth process.

	Altogether, the combination of geometric simplicity of the traces, their fractal scaling properties, and the detailed statistical analysis of the reconstructed driving functions provides an appropriate and robust level of validation for the use of locally regulated Loewner dynamics in our biological growth context.

	\section{Global Results}
	
		\label{SectionGlobalResults}

\begin{figure*}
		\centering
		
		\begin{subfigure}{0.48\linewidth}
			\centering
			\includegraphics[width=\linewidth]{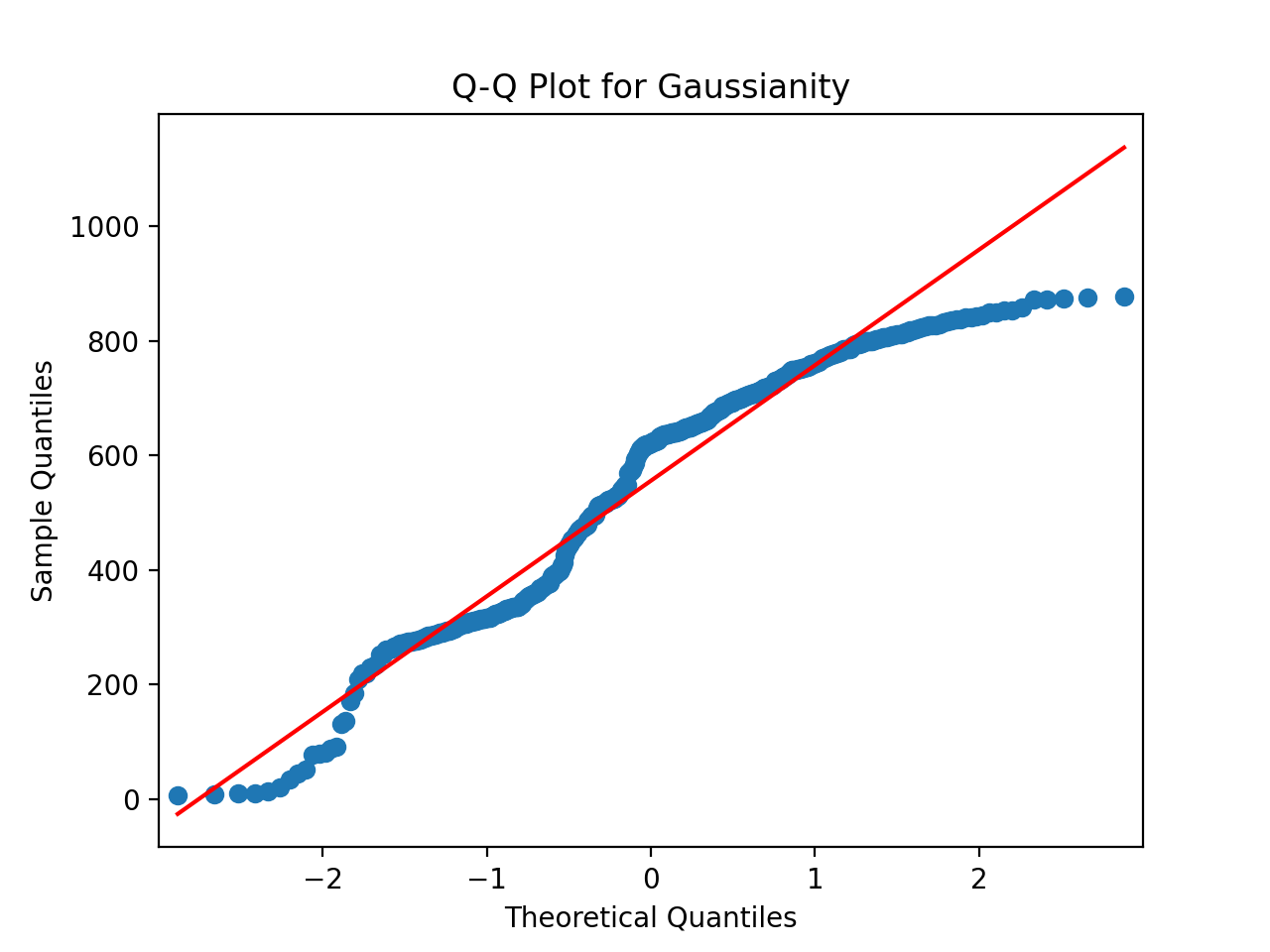}
			\caption{Global  Q-Q plot.}
		\end{subfigure}
		\hfill
		\begin{subfigure}{0.48\linewidth}
			\centering
			\includegraphics[width=\linewidth]{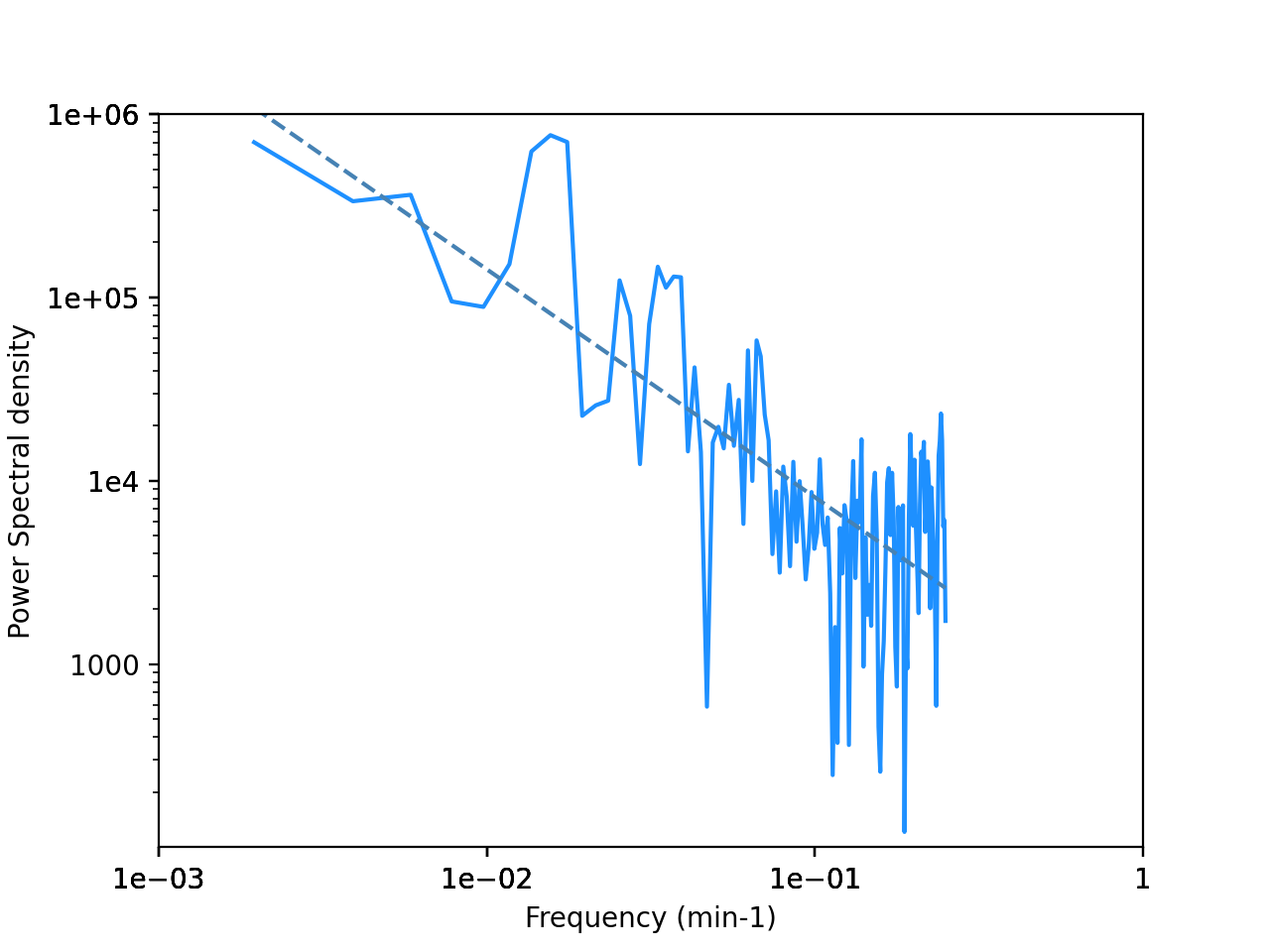}
			\caption{Global power spectral density with regression.}
		\end{subfigure}
		
		\vspace{0.25cm}
		
		\begin{subfigure}{0.48\linewidth}
			\centering
			\includegraphics[width=\linewidth]{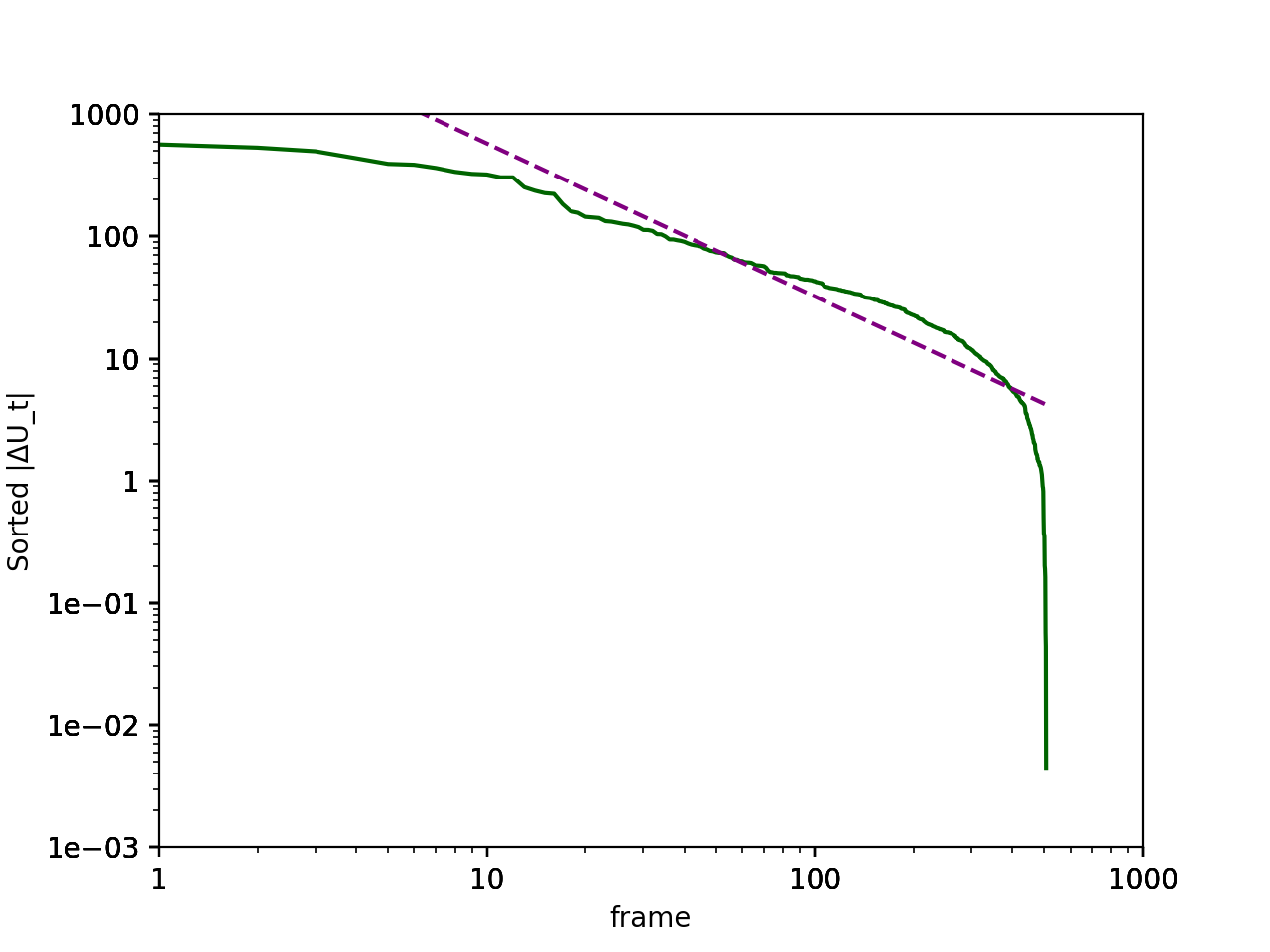}
			\caption{Log-log plot of ordered absolute increments~\mbox{$\Delta U_k$}. The absence of a stable scaling plateau and the rapid tail decay are incompatible with heavy-tailed behavior and support Gaussian statistics.}
		\end{subfigure}
		\hfill
		\begin{subfigure}{0.48\linewidth}
			\centering
			\includegraphics[width=\linewidth]{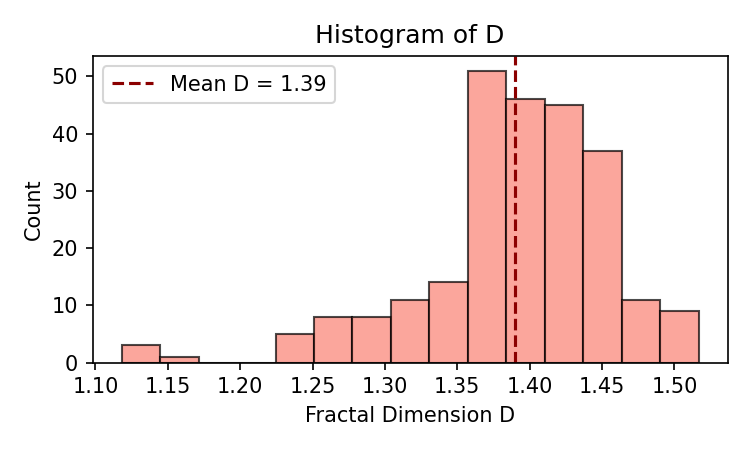}
			\caption{Distribution of the global fractal dimension.}
		\end{subfigure}
		
		\vspace{0.25cm}
		
		\begin{subfigure}{0.48\linewidth}
			\centering
			\includegraphics[width=\linewidth]{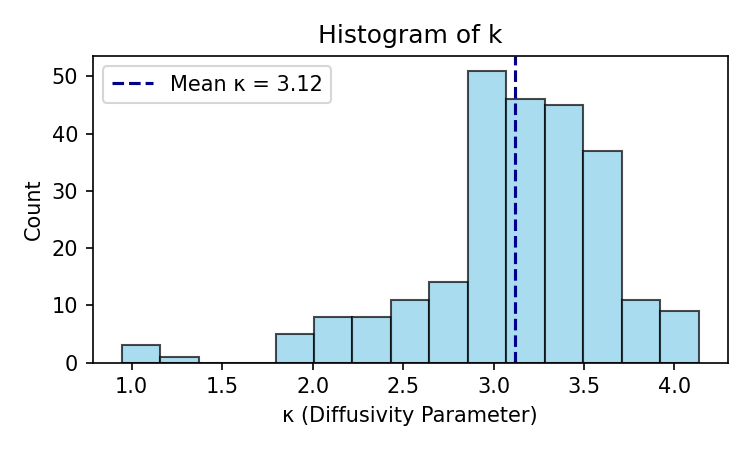}
			\caption{Distribution of the diffusivity parameter~$\kappa$.}
		\end{subfigure}
		
		\caption{\textbf{Global statistical diagnostics of the reconstructed Loewner driving function.}
			 Q-Q plot, power spectral density, Hill-type tail analysis, and distributions of fractal dimension and diffusivity parameter~$\kappa$ computed on the largest connected component of the pseudopod network.}
		\label{QQPlotsPSDEtHistogrammesGlobauxSLEPseudopodsBiggestComponent}
	\end{figure*}

	We now summarize the global statistical behavior of the reconstructed driving signal across the four structural levels introduced above.  Here, \emph{global} is understood in the sense that the reconstruction of the driving function has been made for the largest connected component of the considered structure. Across all four structural levels, the reconstructed driving signals exhibit statistically coherent behavior, including approximately Gaussian marginal distributions and near-linear variance growth, with quantitative variations in the geometrically inferred parameter~$\kappa$ depending on the degree of structural constraint.
	The corresponding diagnostics are shown in Figure~\ref{QQPlotsPSDEtHistogrammesGlobauxSLEPseudopodsBiggestComponent}, including the  Q-Qplot, the power spectral density (PSD) with regression, and the histograms of the global fractal dimension and of the geometrically inferred parameter~$\kappa$. The full set of diagnostics for the remaining structural levels is provided in the SI Appendix.
	
	\begin{enumerate}

		\item[$i$.]  \ul{Pseudopods}:  
		The pseudopod component exhibits the most pronounced Brownian-type statistical signatures.

		\item[$ii$.]  \ul{Network}: 
		The full network displays similar behavior, with moderate quantitative deviations.

		\item[$iii$.]  \ul{Dimming regions}: 
		This level reveals increased dispersion in the estimated diffusivity parameter.
		
		\item[$iv$.]  \ul{Brightening regions}: 
		The expanded network exhibits the strongest structural constraints, reflected in broader variability of the diagnostics.
		
	\end{enumerate}		
	
	Altogether, these results suggest that Brownian-type driving behavior is most prominent at actively expanding fronts and progressively modulated by internal structural constraints.
	
	A detailed comparative discussion of these behaviors is provided below.
	
	\newpage

	
	\clearpage

	 \section{Discussion}
	 \label{SectionDiscussion}

	 \subsection{Main Findings and Interpretation}
	 
	 The present work provides evidence that growth fronts of a living organism can be described, over experimentally accessible spatiotemporal scales, by an emergent Loewner-type dynamics. 
	 By reconstructing an effective Loewner driving signal from time-resolved experimental images, we find that the driving statistics at actively expanding pseudopod fronts exhibit Brownian-type signatures (including approximately Gaussian marginals, near-linear variance growth of increments, and a power spectral density consistent with Brownian scaling). 
	 Within this inferential framework, the estimated diffusivity parameter~$\kappa$ quantifies the intensity of driving fluctuations and provides a compact descriptor of front stochasticity.
	 
	Our multilevel analysis shows that, beyond expanding fronts, the same diagnostics remain broadly coherent across progressively more internal and structurally constrained representations (network, dimming or brightening regions), while displaying systematic quantitative modulations. 
	 Altogether, these results support the fact that Brownian-type driving behavior is most prominent at expansion-driven interfaces, while being progressively shaped by internal transport constraints and network reorganization.
	 
	 To our knowledge, this is the first explicit reconstruction of a Loewner driving function from a living growth interface.

	 \subsection{From Loewner Inference to Biological Growth}
	 
	 The Loewner framework transforms an evolving planar interface into a one-dimensional boundary signal, enabling quantitative comparisons across windows, experiments, and structural levels. 
	 In contrast with geometric descriptors alone (e.g., fractal dimension), the driving function encodes dynamical information and permits direct statistical tests of Brownian-type behavior. 
	 This is particularly relevant in biological contexts, where growth combines stochastic fluctuations, local regulation, and global constraints imposed by transport architecture.

	 \subsection{Relation to Previous Empirical SLE and Conformal Invariance Studies}
	 
	 Experimental signatures of conformal invariance have been reported in biological systems, notably in collective cell flows, where conformal invariance and SLE-consistent statistics were observed without an explicit reconstruction of the driving function~\cite{BenjaminAndersenEtAlEvidenceOfRobustUniversalConformalInvarianceInLivingBiologicalMatter2024}. 
	 In fluid and turbulence settings, SLE-like geometric indicators have been investigated, yet later work emphasized that geometric criteria may be insufficient to recover a genuine SLE mechanism, and that driving-function diagnostics are essential~\cite{TomKennedyComputingTheLoewnerDrivingProcessOfRandomCurvesInTheHalfPlane2008}. 
	 In this perspective, our contribution is to provide an explicit driving-function reconstruction and to evaluate Brownian-type diagnostics directly on experimental growth fronts.
	 A full conformal reparametrization in half-plane capacity time would allow a direct dynamical identification of~$\kappa$, but this lies beyond the scope of the present biological study.

	 \subsection{Methodological Limitations and Robustness}
	 
	 Several limitations are intrinsic to inverse Loewner inference from pixel-based data. 
	 First, discretization and finite resolution may induce lattice effects, segmentation artefacts, and biases in conformal-time parameterization. 
	 As emphasized in the numerical SLE literature, direct reconstruction of conformal maps is highly sensitive to such effects, which is the reason why driving-function based approaches are typically preferred~\cite{TomKennedyNumericalComputationsForTheSchrammLoewnerEvolution2009,TomKennedySimulatingSelfAvoidingWalksInBoundedDomains2013,TomKennedyGregLawlerLatticeEffectsInTheScalingLimitOfTheTwoDimensionalSelfAvoidingWalk2013}. 
	 In the present work, we mitigate these issues by focusing on robust statistical diagnostics of the inferred driving signal, while performing analyses across multiple windows and structural levels, and by using complementary tests (Q--Q plots, increment statistics, variance growth, PSD scaling, tail diagnostics). 
	 Additional methodological details and supplementary diagnostics are provided in the SI Appendix.
	 
	 Second, the interpretation of $\kappa$ is effective: unlike classical Brownian SLE, no assumption is made \emph{a priori} on independence or exact Gaussianity. 
	 
	 This distinction is crucial: geometric similarity to SLE is not sufficient, and only driving-level diagnostics can support such an identification.
	 
	 Accordingly, our results should be understood as evidence for Brownian-like driving signatures rather than as a proof that the biological interface is an exact \mbox{SLE$_\kappa$} trace at all scales.

	 \subsection{Biological Consequences --  Growth Fronts and Transport Architecture}
	 
	 From the biological point of view, pseudopod fronts correspond to active exploration and expansion, while the internal network is used for transport, redistribution, and structural reinforcement. 
	 Our multilevel results suggest a separation between expansion-driven fluctuations at the boundary and constraint-driven correlations associated with internal architecture. 
	 This highlights a possible quantitative link between morphogenesis at the interface and the reorganization of transport networks, potentially modulated by environmental conditions, nutrient distributions, or mechanical constraints.

	 \subsection{Perspectives}
	 
	 Several directions for future work emerge from this framework.
	 First, a systematic exploration of environmental parameters (substrate heterogeneity, nutrient gradients, illumination, repellents/attractants) may reveal controlled transitions between Brownian-like and non-Brownian driving regimes. 
	 Second, the multilevel framework suggests coupling the inferred driving statistics to network descriptors (connectivity, vein thickness distributions, skeleton topology) to quantify how internal reorganization feeds back onto boundary growth. 
	 Finally, the Loewner-based approach introduced here suggests that conformal-growth mechanisms may not be restricted to critical statistical systems, but may also emerge in regulated living matter.

	\section*{Acknowledgments}
	
	The authors thank the interdisciplinary research environment that made this collaboration possible, in particular, ENS Lyon and the INRIA Team MOSAIC.

	\section*{Author contributions}
	
	C.D. initiated the project, identified the Schramm--Loewner type growth process in biological systems, developed the mathematical framework, designed the Loewner analysis pipeline, supervised the numerical implementation, performed the Loewner reconstruction from experimental data, and wrote the manuscript.
	
	A.B. performed the image analysis using the \emph{Cellects} software, developed the Python code, and contributed to data analysis and interpretation.
	
	N.R. reviewed the numerical implementation and contributed to data analysis and discussions of the results.
	
	A.D. designed and performed the laboratory experiments, acquired the experimental data, and contributed to data analysis and interpretation.
	
	M.L.L. contributed to the mathematical analysis and discussions of the results.

	All authors discussed the results and approved the final manuscript.
	
	\paragraph{Funding.}
	This work was supported by the MITI CNRS 80PRIME program, the CNRS interdisciplinary programs \emph{Conditions Extr\^emes} and \emph{Formes Optimales}, and the Agence Nationale de la Recherche (ANR) through the project FRACTALS (ANR-24-CE45-3362).
	The research of M.~L.~L. was supported by the Burton Jones Endowed Chair in Pure Mathematics, as well as by grants from the U. S. National Science Foundation.

	\section*{Competing interests}
	
	The authors declare no competing interests.

	\section*{Data availability}
	
	The experimental data supporting this study are available from the authors upon reasonable request. 
	The numerical scripts are available on the following GitHub link: \url{https://github.com/Aurele-B/PolycephalumSLE.git}.

	\section{Supporting Information}
	
	Supporting Information includes detailed mathematical proofs, additional numerical diagnostics, and supplementary figures.

	\subsection{Experimental Methods}

	\subsubsection{Species}
	The slime mold  Physarum polycephalum is a unicellular organism belonging to the Amoebozoa. Its vegetative stage, the plasmodium, is a giant motile cell consisting of a syncytium of nuclei and an intracellular cytoskeleton that forms a complex cytoplasmic network of veins. Its cytoplasm comprises a viscous phase (ectoplasm) and a fluid phase (endoplasm), characterized by different concentrations of fibrous proteins. The ectoplasm, which contains actin and myosin, forms the contractile walls of the veins. Within these veins flows the endoplasm, which contains organelles such as nuclei and mitochondria. The ectoplasm and endoplasm can interconvert. When starved, a plasmodium can enter a dormant stage called a sclerotium and later revert to the plasmodial state when transferred to a fresh food medium. In this study, we used the Australian strain provided by Southern Biological (Victoria, Australia). One sclerotium was revived in 2021 for this experiment.\\
	
	\noindent \underline{Rearing conditions.}\\
	The slime molds were reared on a 1\% agar medium supplemented with rolled oat flakes. They were fed daily, and the agar medium was replaced each day. The slime molds were two weeks old at the start of the experiment. All experiments were conducted in the dark at a temperature of 25 °C and a relative humidity of 80\% for a duration of 48 h.

	\subsubsection{Experimental setup}
	We observed the behavior of a single slime mold for this paper (diameter: 13 mm) exploring a circular arena consisting of a 90 mm diameter Petri dish filled with plain 1\% agar (homogenous environment). The arena was monitored for 48 h using time-lapse photography. Images were taken at regular intervals using a camera connected to a self-programmed Arduino controller. A LED panel positioned beneath the Petri dish was switched on for 3 s during image acquisition.

	\subsection{Mathematical Context: Brownian and Non-Brownian Stochastics, Loewner Context}
	
	\label{SupportingInformationBrownianStochastics}
	\hskip 0.5cm For the benefit of the reader who may not be familiar with mathematical notions pertaining to Brownian motion, we shall first recall several relevant definitions.

	\begin{definition}[\textbf{L\'evy Process}]

		A~\emph{L\'evy process}~\mbox{$X_t $} is a continuous-time stochastic process such that~\mbox{$X_0 =0$} almost surely, with stationary and independent increments.

		See the book by David~Applebaum~\cite{DavidApplebaumLevyProcessesAndStochasticCalculus2009}, Chapter~I, on page~\mbox{$43$}.
		
	\end{definition}

	\begin{proposition}[\textbf{Brownian Motion as a Wiener Process}]  
		
		\label{PropositionBrownianMotionAsAWienerProcess}
		
		Brownian motion, denoted by~\mbox{$B_t $}, is a \emph{Wiener process}, i.e., a continuous-time stochastic process characterized by the following four properties:
		
		\begin{enumerate}
			\item[$i$.] ~\mbox{$B_0=0 $}.
			\item[$ii$.] It is almost surely continuous.
			
			\item[$iii$.] It has independent increments (which makes it a \emph{L\'evy process}); i.e., for all~\mbox{$(s_1,t_1,s_2,t_2)\, \in\, \R_+^4 $} such that~\mbox{$0 \leq s_1 < t_1 \leq s_2 < t_2  $},~\mbox{$B_{t_1}-B_{s_1}$} and~\mbox{$B_{t_2}-B_{s_2}$} are independent random variables.
			
			\item[$iv$.] For all~\mbox{$(s,t)\, \in\, \R_+^2 $} such that~\mbox{$0 \leq s \leq t$} :
			
			$$B_t-B_s \sim {\cal N} (0, t-s)\, ,$$
			
			where~\mbox{${\cal N} (0, t-s) $} denotes the centered normal law, with mean zero and associated variance~\mbox{$t-s $}.
		\end{enumerate}

		See the book by Ioannis Karatzas and Steven~E.~Shreeve~\cite{IKaratzasSShreeveBrownianMotionAndStochasticCalculus1991}, on page~\mbox{$47$}, at the beginning of Section~2, along with the book by Peter~M\"orters and Yuval Peres~\cite{PeterMortersAndYuvalPeresBrownianMotion2010}.

	\end{proposition}

	\begin{definition}{\textbf{Half-Plane Capacity~\cite{GregoryFLawlerConformallyInvariantProcessesInThePlane2005}}}  
		\label{HalfPlaneCapacity}
		
		Given a compact subset~\mbox{$K \subset \bar{H}$} such that~\mbox{$\bar{H}\setminus K $} is simply connected and~\mbox{$K  \cap \R \neq \emptyset$}, we consider the unique conformal map~\mbox{$g_K \, : \, \bar{H}\setminus K  \to \bar{H}$}, normalized when~\mbox{$ |  z |  \to \infty$} via
		
		$$g_K(z)= z + \displaystyle  \frac{\mathrm{hcap}\, (K)}{z} + {\cal O} \left ( \frac{1}{z^2}\right) \, .$$
		
		The \emph{half-plane capacity}~\mbox{$\text{kcap}\, (K)$} of the subset~$K$ is defined, for all Loewner time~\mbox{$t \geq 0$}, via
		
		$$\text{kcap}\, (K)= 2\, t\, ,$$
		
		\noindent corresponding to the factor of~\mbox{$\displaystyle \frac{1}{z}$} in the above expansion of~\mbox{$g_K(z)$}. See~\cite{GregoryFLawlerConformallyInvariantProcessesInThePlane2005}, Definition~\mbox{$3.\, 35$}, on page~\mbox{$55$}.
		
	\end{definition}

	In the classical chordal SLE$_\kappa$ setting, the Loewner driving function is given by~\mbox{$U_t=\sqrt{\kappa}\, B_t$}, which directly connects Brownian motion to conformally invariant growth.

	\subsection{Methods: Spectral and Tail Diagnostics}
	
	\begin{proposition}[\textbf{Power Spectrum of a Stationary Signal}]  
		
		\label{PropositionPowerSpectrumOfAStationarySignal}
		
		Given a stationary time signal~\mbox{$\xi(t)$} with zero mean~\mbox{$\left \langle \xi(t)\right \rangle = 0$}, the associated \emph{autocorrelation}~\mbox{$R_\xi(\tau)$} is given, for all~\mbox{$\tau \, \in\, \R$}, by
		
		$$R_\xi(\tau)= \left \langle \xi(t)\, \xi (t+\tau)\right \rangle
		=\displaystyle \lim_{T \to \infty} \frac{1}{T}\, \int_0^T \,  \xi(t)\, \xi (t+\tau) \, dt .$$

		The \emph{power spectral density}~\mbox{$S_\xi $} of~$\xi$ is then defined as the following Fourier transform, defined, for all~\mbox{$\omega\, \in\, \R$}, by

		$$S_\xi(\omega) = \displaystyle \int_{-\infty}^\infty R_\xi(\tau) \, e^{-2\, i\,\pi\,\omega\,    \tau } \, d\tau\, .$$

		Equivalently, it can be expressed, for all~\mbox{$\omega\, \in\, \R$}, as
		
		\begin{equation}
			\label{PowerSpectralDensityExpression}
			S_\xi(\omega) = \displaystyle \lim_{T \to \infty} \frac{1}{T}\, \left |  \int_{0}^T  \xi ( \tau)   \, e^{-2\, i\,\pi\,\omega\,    \tau } \, d\tau\right |  ^2\, .\end{equation}

		See the book by Crispin~Gardiner~\cite{CrispinGardinerStochasticMethods1983}, page~\mbox{$59$} of Chapter~\mbox{3} on Markov processes.\\

		For a sampled signal -- i.e., a discrete-time sequence of signals~\mbox{$\left \lbrace \xi_k \right\rbrace _{1\leq k\leq N}=\left \lbrace \xi(t_k)  \right\rbrace _{1\leq k \leq N}$}, with~\mbox{$ N\, \in\, \N^\star$}, and where, for~\mbox{$ 2 \leq k \leq N$},~\mbox{$ t_k=t_{k-1}+\Delta t$}, with~\mbox{$  \Delta t > 0$} -- the associated periodogram estimate of the power spectral density~\mbox{$\left \lbrace S_k(\omega_k) \right\rbrace _{1\leq k \leq N}=\left \lbrace  S_{\xi_k} (t_k)  \right\rbrace _{1\leq k \leq N}$} can be computed by means of the Welch method; see the seminal paper by Peter~D.~Welch~\cite{PeterDWelchTheUseOfFastFourierTransformForTheEstimationOfPowerSpectraAMeethodBasedOnTimeAveragingOverSshortModifiedPeriodograms1967}.\\
		
		By letting~\mbox{$\Delta \omega  = \displaystyle \frac{1}{N\, \Delta t}$}, and, for~\mbox{$ 1\leq k \leq N $},~\mbox{$\omega_k=k\, \Delta \omega = \displaystyle \frac{ k }{N\, \Delta t}$}, the periodogram estimate~\mbox{$\left \lbrace P_k \right\rbrace _{1\leq k \leq N}$}  of the discrete powers associated with the power spectral density (PSD) of the signal~\mbox{$\left \lbrace \xi_k \right\rbrace _{1\leq k \leq N}$} -- namely, the squared modulus of the discrete Fourier transform of the sampled signal~\mbox{$\left \lbrace( \xi_k \right\rbrace _{1\leq k \leq N}$}, normalized by the total duration of the signal (i.e., the signal length) -- is given,  for~\mbox{$ 1 \leq k \leq N$}, by

		$$S_k(\omega_k)= \displaystyle \frac{1}{N\, \Delta \omega}\, \left |  \hat{\xi} (\omega_k) \right|  ^2 \, ,$$

		\noindent where, for~\mbox{$ 1 \leq k \leq N$}, the discrete Fourier transform~\mbox{$\hat{\xi} (\omega_k)$} is given by

		$$\hat{\xi} (\omega_k) = \displaystyle \sum_{j=1}^N \xi_j\,   e^{-2\, i\,\pi\, \frac{j\,k }{N}} =  \displaystyle \sum_{j=1}^N \xi_j\,   e^{-2\, i\,\pi\,j\, \omega_k} \, .$$

		The associated discrete powers~\mbox{$ P_k = P(t_k)  $}, for~\mbox{$ 1 \leq k \leq N$}, are then obtained via
		
		$$P_k= \displaystyle \frac{\left |  \hat{\xi} (\omega_k) \right|  ^2 }{N}\,.$$

		By summing over all frequencies, we then obtain the following Parseval relation, which is the discrete counterpart of the power spectral density given in Relation~\mbox{($\cal R$\ref{PowerSpectralDensityExpression})}, as
		
		$$P = \displaystyle \sum_{k=1}^N\frac{\left |  \hat{\xi} (\omega_k) \right|  ^2 }{N}\,.$$

		\noindent See the book~\cite{JuliusBendatAllanPersolRandomData2010} by Julius Bendat and Allan G.~Piersol.\\

	\end{proposition}

	\begin{corollary}[\textbf{Brownian Power Spectrum}]  
		
		\label{CorollaryBrownianPowerSpectrum}
		
		The power spectral density~\mbox{$S_{B_t}$} of~\mbox{$B_t$} is given, for all~\mbox{$\omega \, \in\, \R$}, by
		
		$$S_{B_t}(\omega)  
		= \displaystyle \frac{1}{4\, \pi^2\, \omega^2}\, .$$	
		
		For the sake of concision, we will use the standard following notation for the associated power law,
		
		$$S_{B_t}(\omega)  	  \propto \omega^{-2}\, .$$
		
	\end{corollary}

	\begin{remark}		
		In the angular frequency convention~\mbox{$\Omega = 2\, \pi \omega$}, the power spectral density~\mbox{$S_{B_t}$} can equivalently be written, in the form~\mbox{$S_{B_t} =    \displaystyle \frac{1}{\Omega^2}$}.
		
	\end{remark}

	In order to test whether the driving increments belong to the domain of attraction of a Gaussian law or of a heavy-tailed stable law, we recall the following definitions.

	\begin{definition}{\textbf{Domain of Attraction of a Probability Distribution~\cite{WilliamFellerAnIntroductionToProbabilityTheoryAndItsApplicationsVolumeII1971}}}  
		
		Let us consider a sequence~\mbox{$\left ( X_n\right)_{n \, \in\,\N^\star }$}, of independent and identically distributed (i.i.d.) random variables, with common distribution function~$F$, along with a probablity law~$G$.\\ The law~$F$ is said to belong to the \emph{domain of attraction} of~$G$ if there exists two sequences\\~\mbox{$\left ( a_n \right)_{n \,\in\, \N^\star} \, \in\, \left ( \R_+^\star\right)^\N$} and~\mbox{$\left ( b_n \right)_{n \,\in\, \N^\star } \, \in\, \R^\N$} (where~\mbox{$ \R_+^\star = \left ]0, + \infty \right[ $}) such that

		$$\displaystyle \frac{\displaystyle\sum_{i=1}^n X_i -b_n}{a_n}\,  \underset{n \to \infty }{\xrightarrow{\displaystyle \mathcal{D}}}  \, Y $$
		
		\noindent with~\mbox{$Y \sim G$}; i.e., the random variables~$Y$ and~$G$ are equal in distribution.\\

		See the book by William~Feller~\cite {WilliamFellerAnIntroductionToProbabilityTheoryAndItsApplicationsVolumeII1971}, Definition~2 of Chapter~\mbox{IV}, on page~\mbox{$172$}.

	\end{definition}

	\begin{theorem}{\textbf{Domain of Attraction of a Stable Law~\cite{WilliamFellerAnIntroductionToProbabilityTheoryAndItsApplicationsVolumeII1971}}}  
		
		\label{TheoremDomainOfAttractionOfTheNormalDistribution}

		A symmetric probability distribution~$F$ not concentrated at a single point is said to belong to \emph{the domain of attraction} of a stable law if and only if its truncated moment function~$\mu$, defined, for all~\mbox{$x \, \in\, \R$}, via~\mbox{$\mu(x)= \displaystyle \int_{-x}^x y^2\, dF(y) $}, varies regularly as~\mbox{$x \to + \infty$}, according to
		
		$$\displaystyle  \mu(x) \sim x^{2-\alpha} \, L (x)  \, , $$

		\noindent where~$L$ denotes a  function varying slowly at~$\infty$, and where~\mbox{$0< \alpha \leq 2$}.
		See the book by William~Feller\\~\cite {WilliamFellerAnIntroductionToProbabilityTheoryAndItsApplicationsVolumeII1971}, Corollary~1 of Chapter~\mbox{XVII}, on page~\mbox{$578$}.\\

		This corresponds to a power-law tail, as~\mbox{$x \to + \infty$}, of the form
		
		$$\mathbb{P} \left ( |  X|  > x \right) \sim C \, x^{-\alpha}\, ,$$
		
		\noindent where~$C$ denotes a strictly positive constant.\\
		
		The case~\mbox{$\alpha=2$} corresponds to the domain of attraction of the normal distribution.\\
		
		The case~\mbox{$0< \alpha \leq 2$} corresponds to~$\alpha$-stable L\'evy processes;  
		see the book by David~Applebaum~\cite{DavidApplebaumLevyProcessesAndStochasticCalculus2009}, Chapter~I, Example~\mbox{$1.\,3.\, 14$}, on page~\mbox{$51$}.

	\end{theorem}

	\begin{proposition}{\textbf{Hill Estimator~\cite{BruceMHillASimpleGeneralApproachToInferenceAboutTheTailOfADistribution1975}}}  
		
		\label{PropositionHillEstimator}
		
		Let us consider a sequence~\mbox{$\left ( X_n\right)_{n \, \in\,\N^\star }$}, of independent and identically distributed (i.i.d.) random variables, issued from a heavy-tailed probability distribution~$\mathbb{P}$ such that, as~\mbox{$x \to \infty$},  
		
		$$\mathbb{P} \left ( |  X|  > x \right) \sim C \, x^{-\alpha}\, ,$$

		\noindent where~\mbox{$\alpha >0$}.
		
		We denote by~\mbox{$\left ( X^{(n)}\right)_{n \, \in\,\N^\star }$} the reordered sequence associated with~\mbox{$\left ( X_n\right)_{n \, \in\,\N^\star }$}, where, for each~\mbox{$n \, \in\,\N^\star  $},
		
		$$X^{(n)} \leq X^{(n-1)} \leq \ldots \leq X^{(2)} \leq X^{(1)} \, .$$

		By choosing a threshold~\mbox{$n_0 \geq 1$}, we obtain the \emph{Hill exponent}~\mbox{$\hat{\alpha}_H $} of the~$\alpha$ tail -- corresponding to the conditional maximum-likelihood estimator of the tail exponent~$\alpha$ -- as
		
		$$\hat{\alpha}_H = (n_0+1)\, \left ( \displaystyle  \sum_{k=1}^{n_0} \ln \left ( \frac{X^{(k)}}{X^{(n_0+1)}} \right) \right)^{-1}\, .$$

		\noindent See the seminal paper by Bruce~M.~Hill~\cite{BruceMHillASimpleGeneralApproachToInferenceAboutTheTailOfADistribution1975}, equation~\mbox{$(3.\,1)$}, on page~\mbox{$1166$}.\\

	\end{proposition}

	In the present work, the Hill estimator is applied to the ordered absolute increments~\mbox{$\left| \Delta U_t \right|$} of the reconstructed Loewner driving function, 
	in order to test for possible heavy-tailed deviations from Gaussian behavior.

	In the classical chordal SLE$_\kappa$ setting, the driving function is~\mbox{$U_t = \sqrt{\kappa} \, B_t$}. Therefore, deviations from the~\mbox{$\omega^{-2}$} spectral law or from Gaussian domain-of-attraction behavior would indicate departures from Brownian-driven conformally invariant growth.
	
 \newpage
	
	\subsection{Supplementary Figures}


	\begin{figure*}[h!]
		\centering
		
		\begin{subfigure}[t]{0.32\textwidth}
			\centering
			\includegraphics[width=\linewidth]{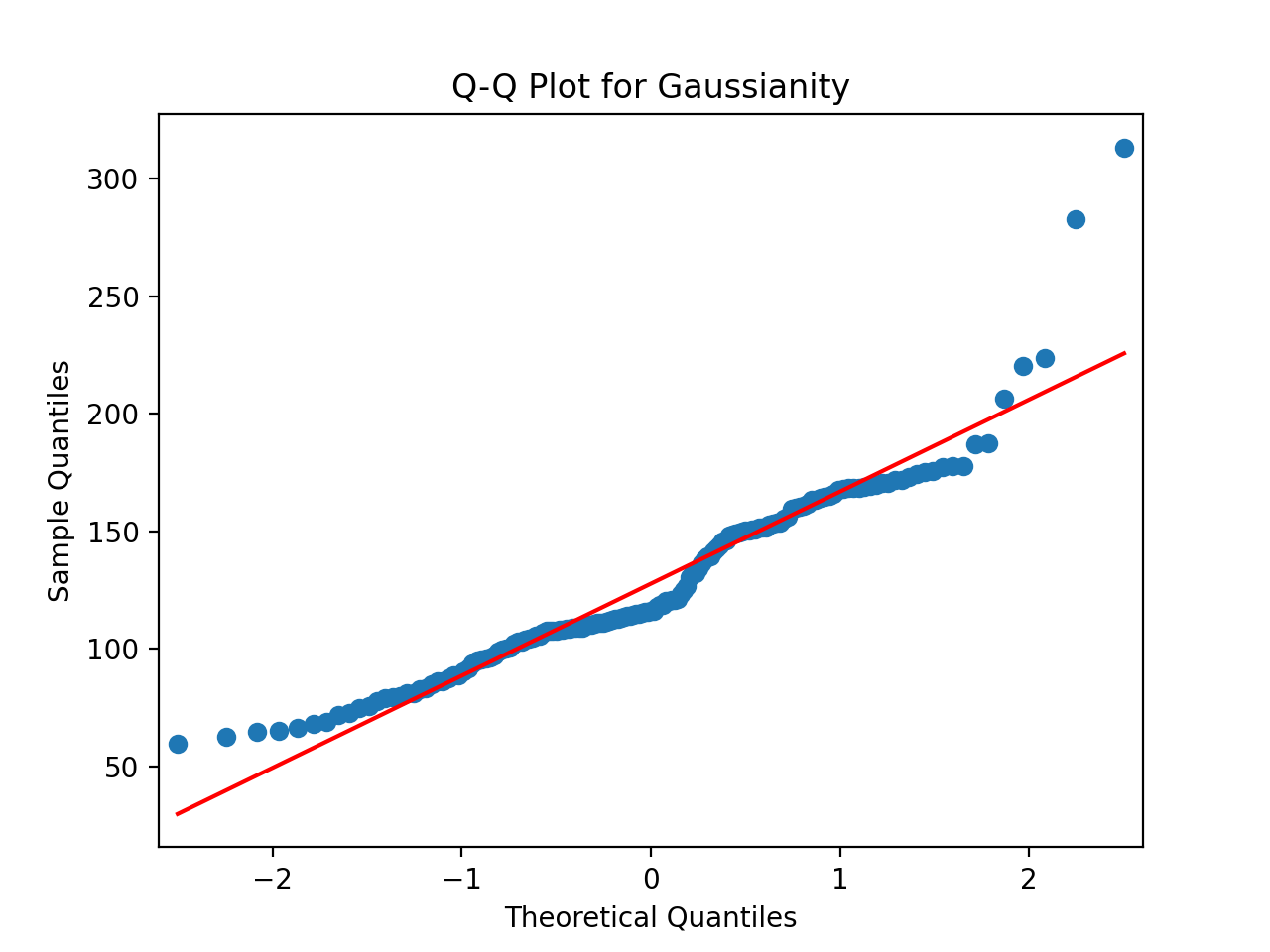}
			\caption{\textbf{Window 2}}
		\end{subfigure}
		\hfill
		\begin{subfigure}[t]{0.32\textwidth}
			\centering
			\includegraphics[width=\linewidth]{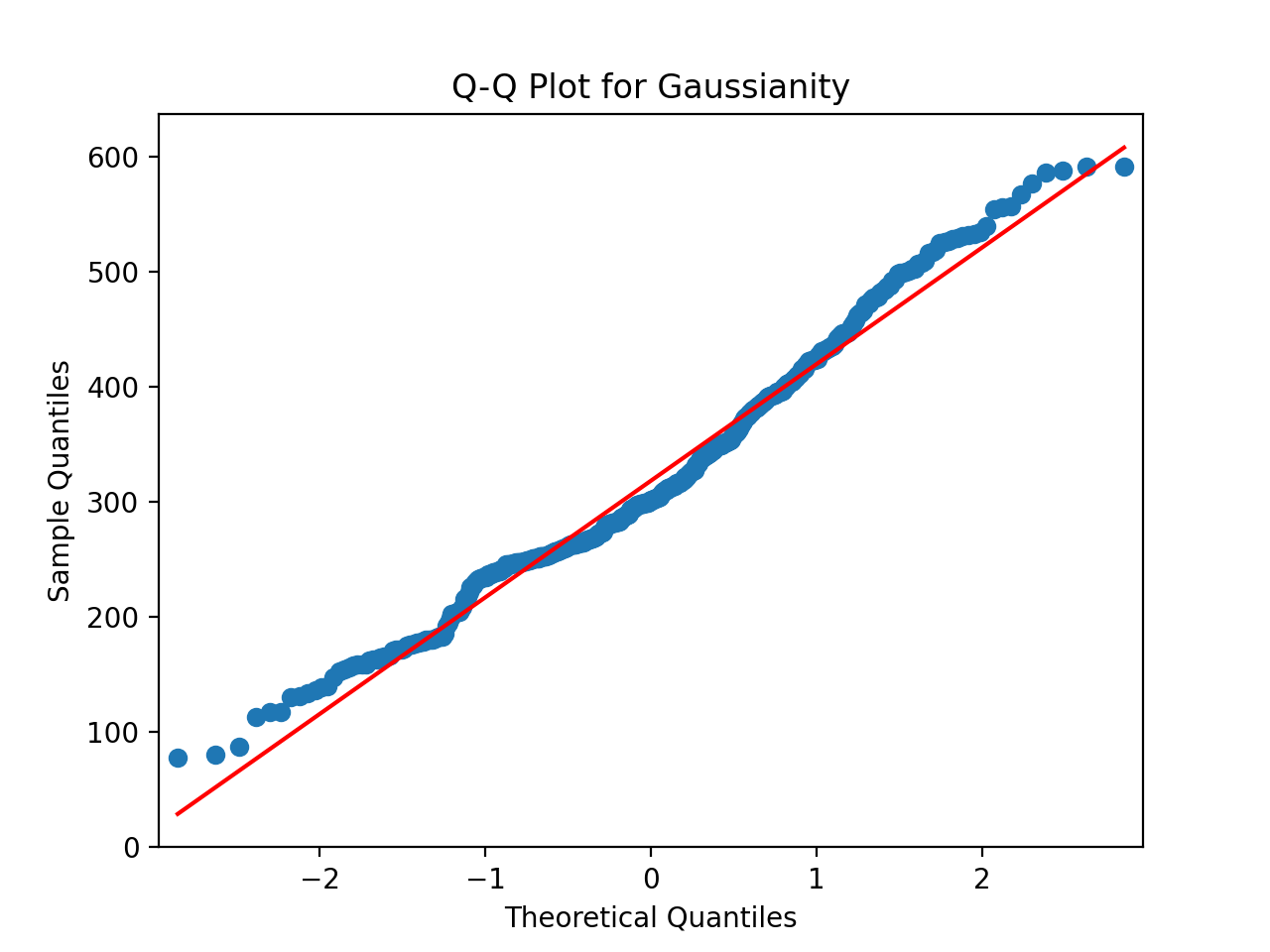}
			\caption{\textbf{Window 3}}
		\end{subfigure}
		\hfill
		\begin{subfigure}[t]{0.32\textwidth}
			\centering
			\includegraphics[width=\linewidth]{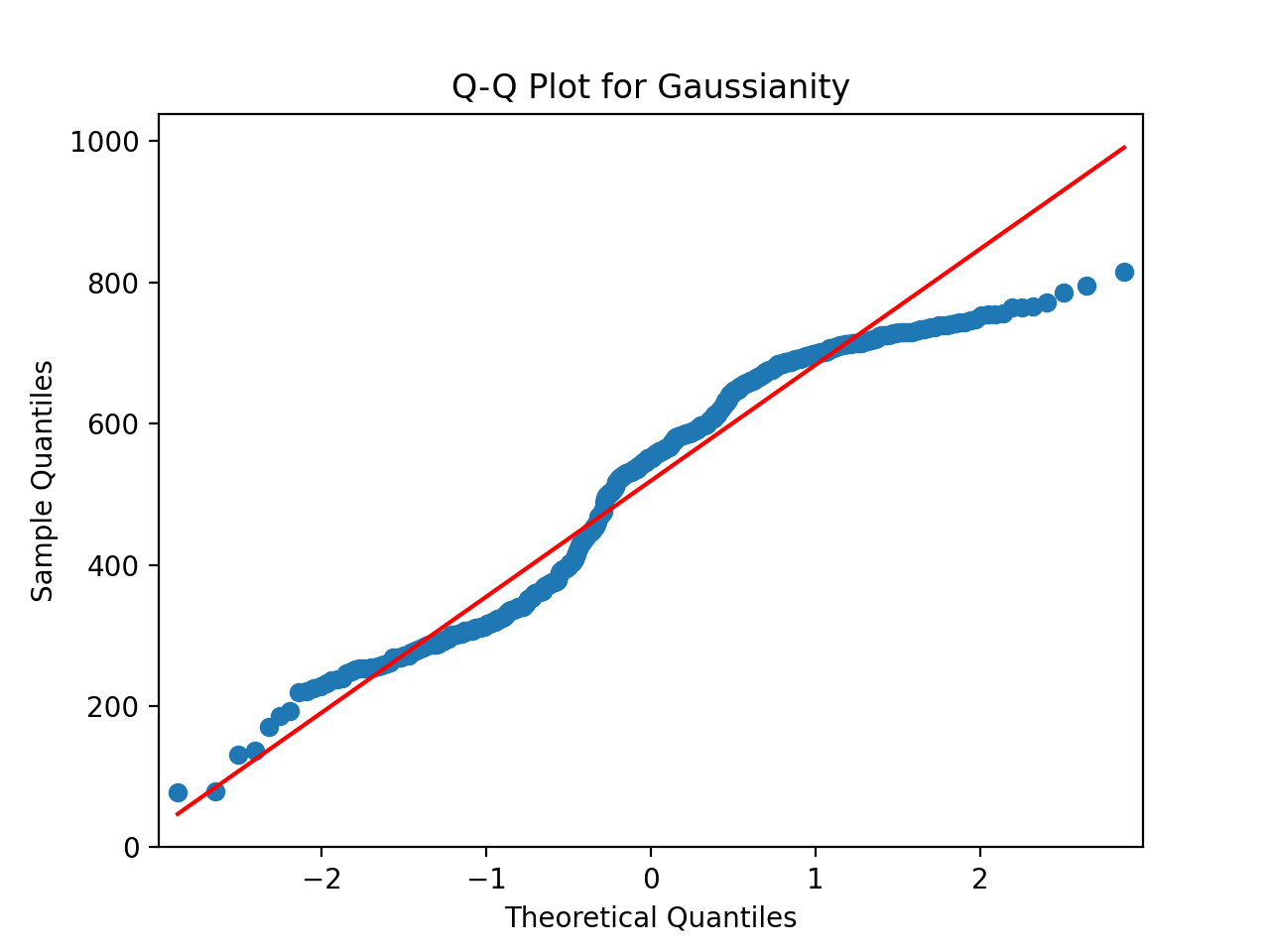}
			\caption{\textbf{Window 4}}
		\end{subfigure}
		
		\vspace{0.35cm}
		
		\begin{subfigure}[t]{0.32\textwidth}
			\centering
			\includegraphics[width=\linewidth]{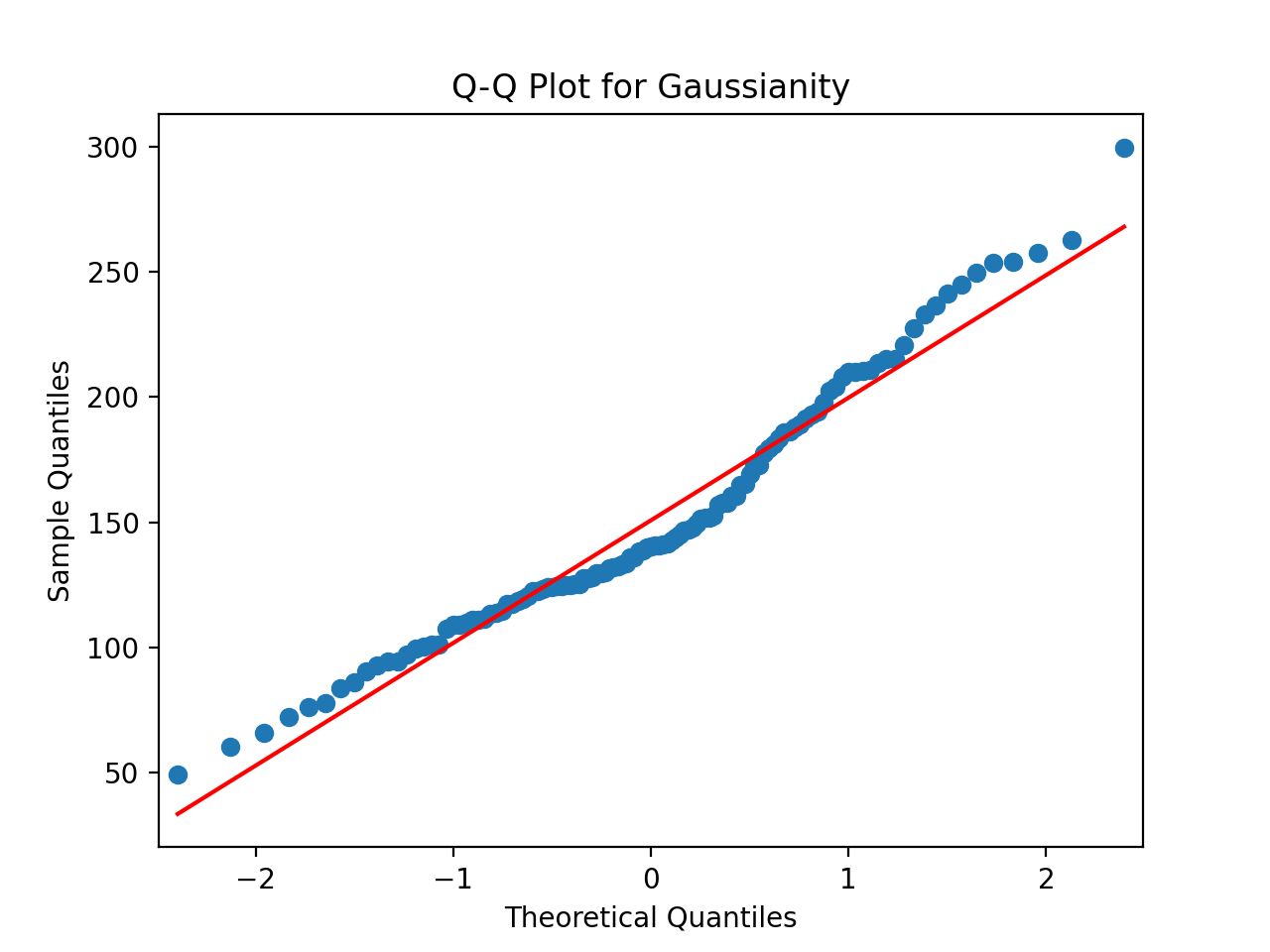}
			\caption{\textbf{Window 7}}
		\end{subfigure}
		\hfill
		\begin{subfigure}[t]{0.32\textwidth}
			\centering
			\includegraphics[width=\linewidth]{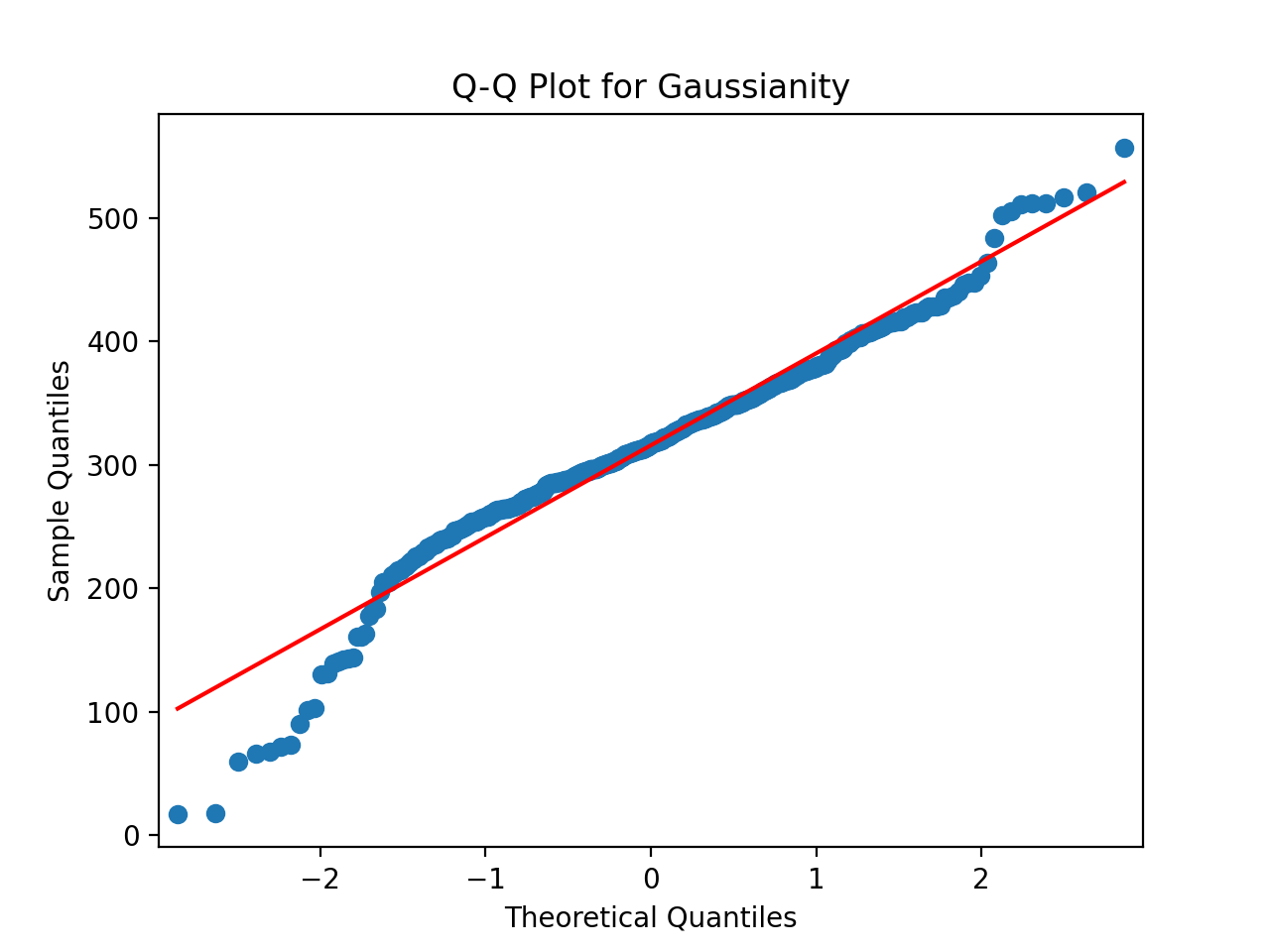}
			\caption{\textbf{Window 8}}
		\end{subfigure}
		\hfill
		\begin{subfigure}[t]{0.32\textwidth}
			\centering
			\includegraphics[width=\linewidth]{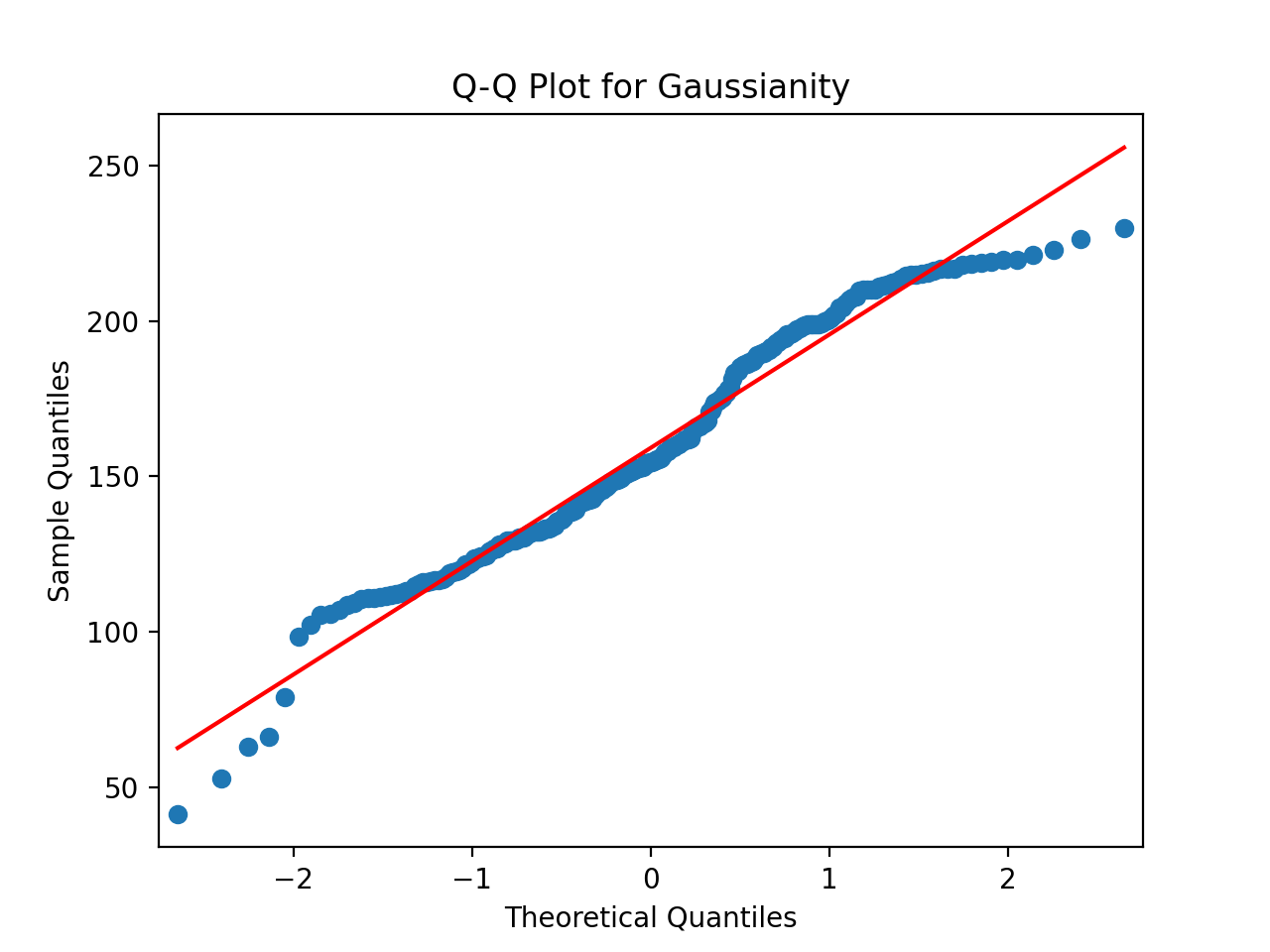}
			\caption{\textbf{Window 12}}
		\end{subfigure}
		
		\caption{\textbf{Local Q-Q plots of the driving function increments $\Delta U_t$ (pseudopods).}
			Representative outer spatial windows (2, 3, 4, 7, 8, 12). The empirical quantiles closely follow the theoretical Gaussian quantiles, with mild deviations at extreme tails.}
		\label{QQPlotsLocauxSLEPseudopodsBiggestComponent}
	\end{figure*}
	

	\begin{figure*}[h!]
		\centering
		
		\begin{subfigure}[t]{0.32\textwidth}
			\centering
			\includegraphics[width=\linewidth]{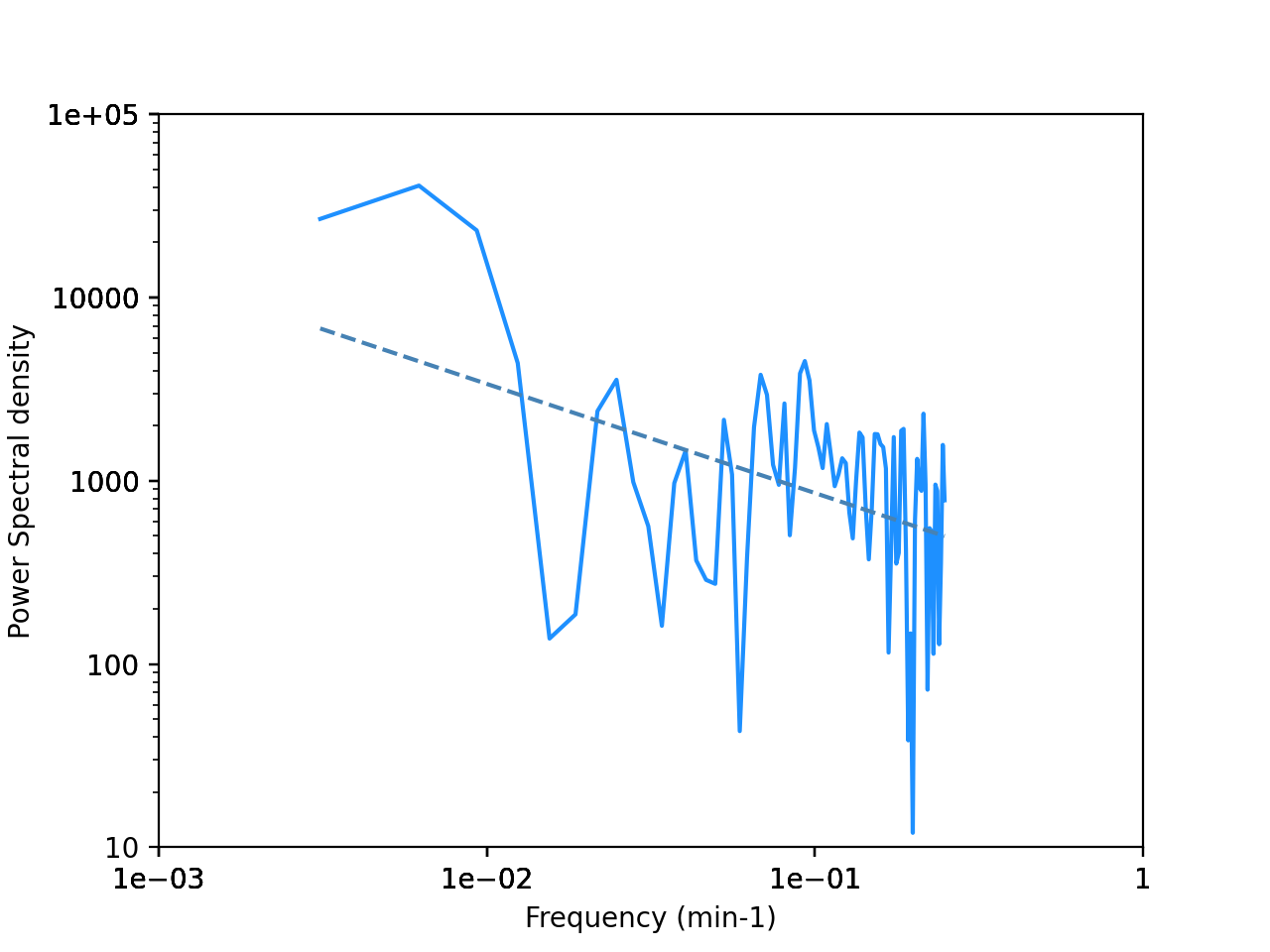}
			\caption{\textbf{Window 2}}
		\end{subfigure}
		\hfill
		\begin{subfigure}[t]{0.32\textwidth}
			\centering
			\includegraphics[width=\linewidth]{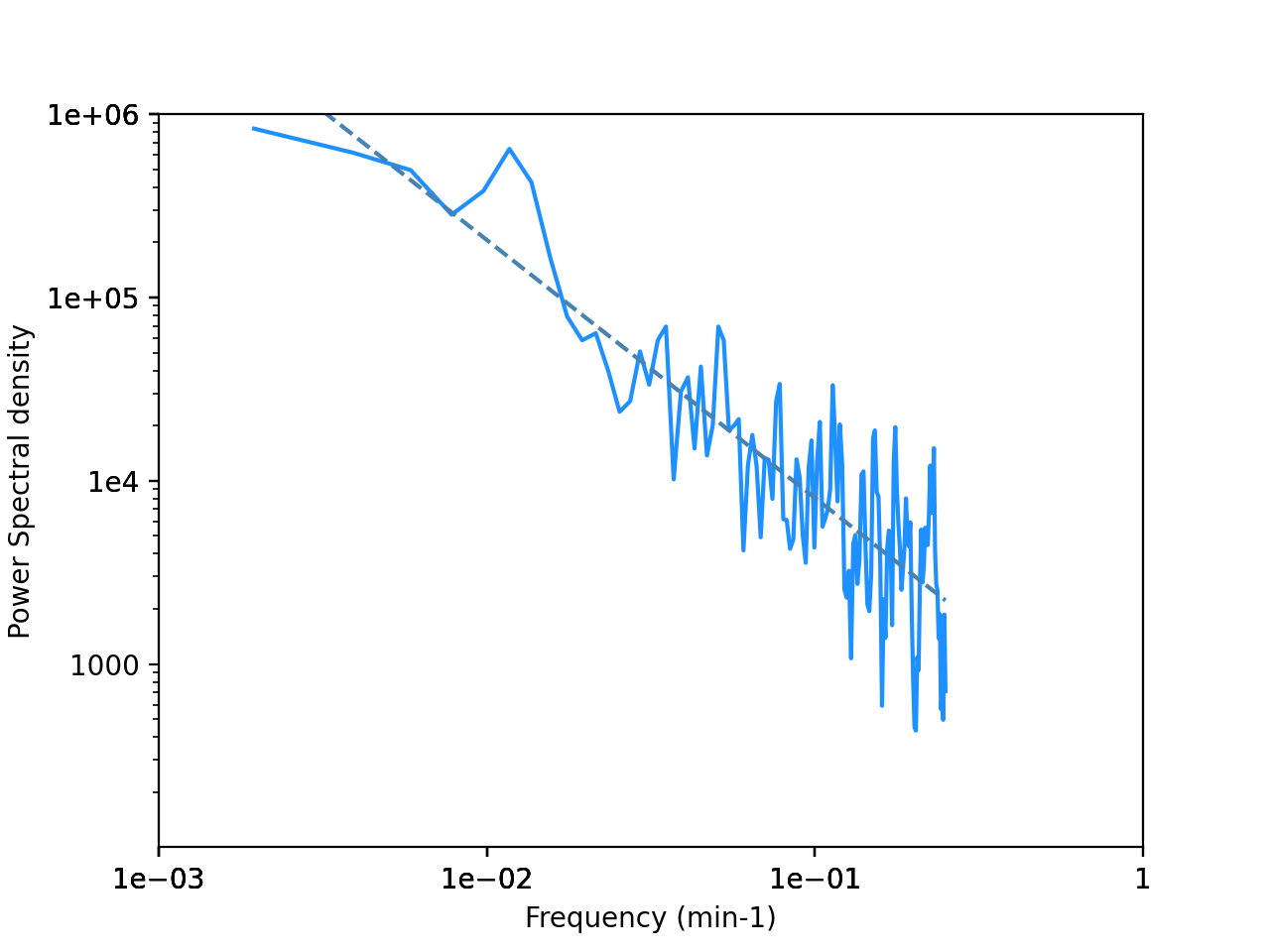}
			\caption{\textbf{Window 3}}
		\end{subfigure}
		\hfill
		\begin{subfigure}[t]{0.32\textwidth}
			\centering
			\includegraphics[width=\linewidth]{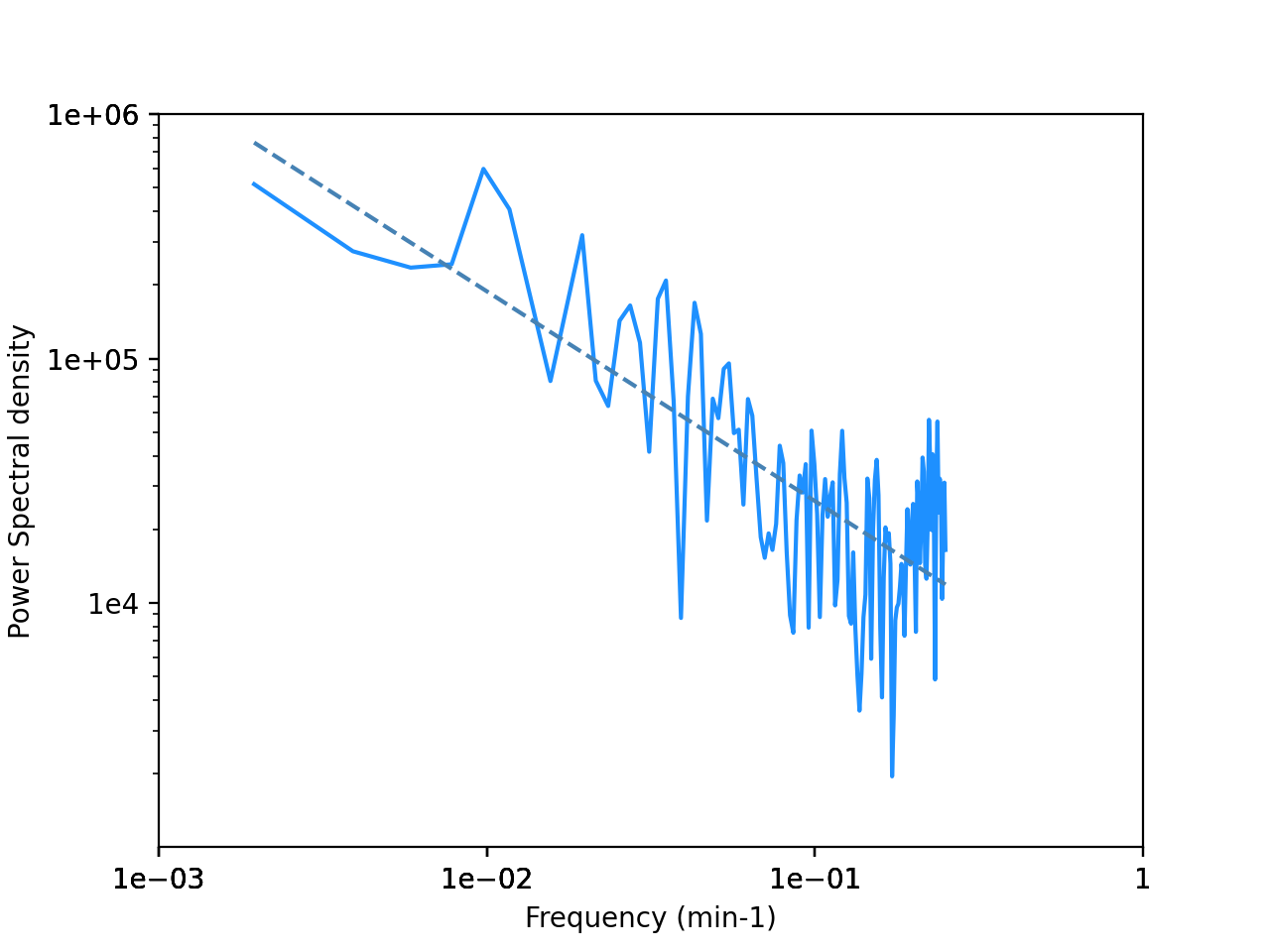}
			\caption{\textbf{Window 4}}
		\end{subfigure}
		
		\vspace{0.35cm}
		
		\begin{subfigure}[t]{0.32\textwidth}
			\centering
			\includegraphics[width=\linewidth]{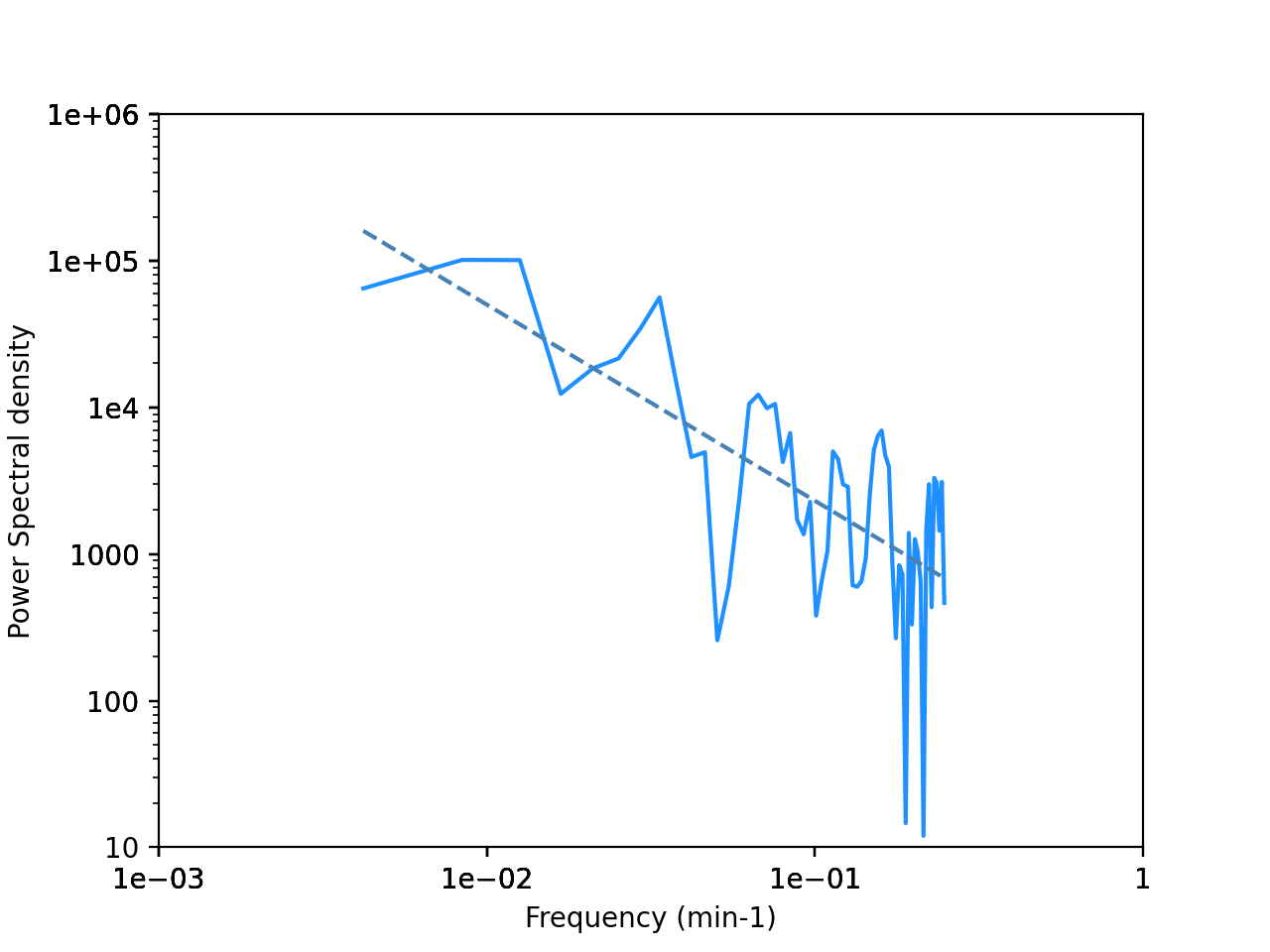}
			\caption{\textbf{Window 7}}
		\end{subfigure}
		\hfill
		\begin{subfigure}[t]{0.32\textwidth}
			\centering
			\includegraphics[width=\linewidth]{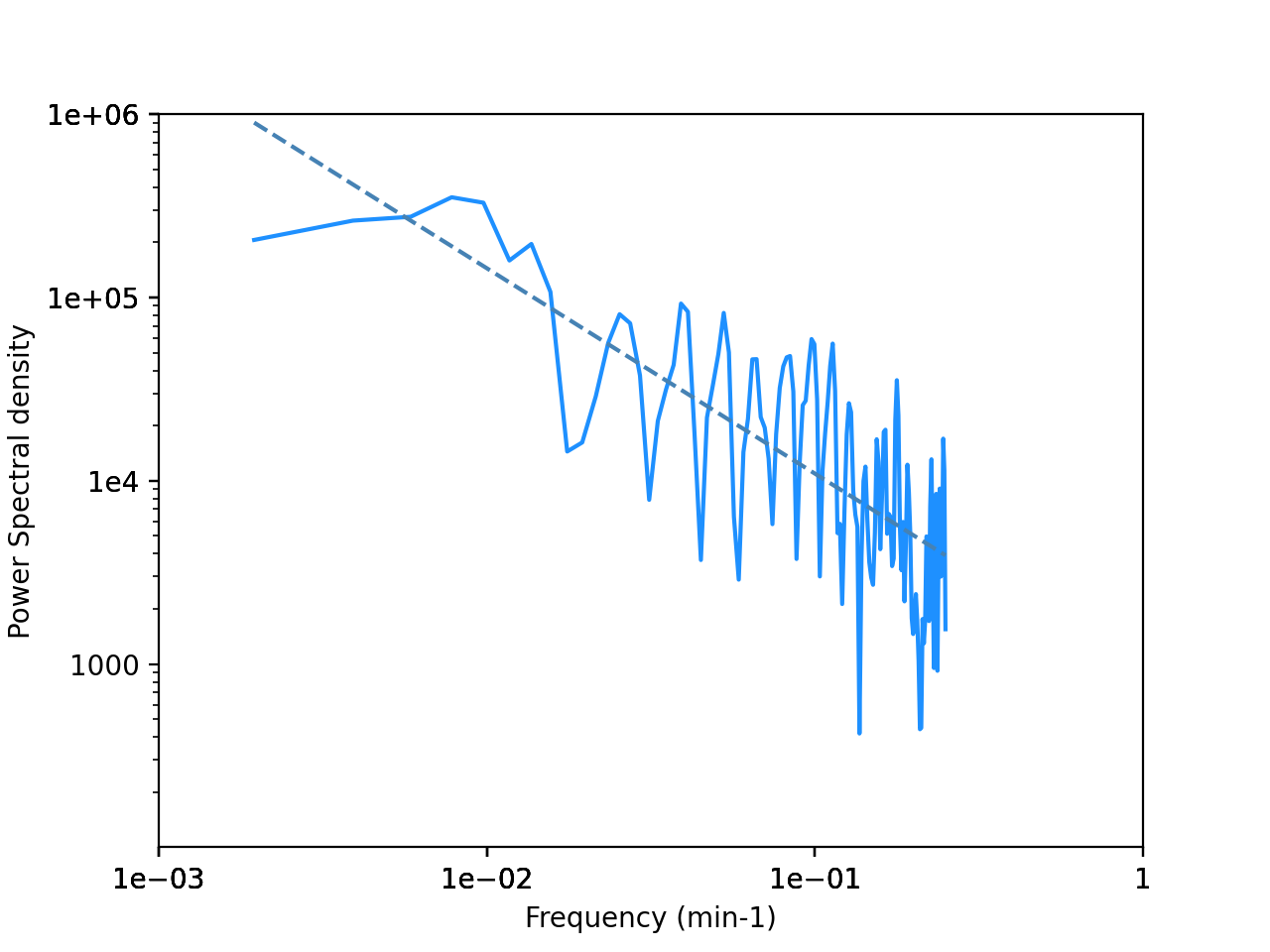}
			\caption{\textbf{Window 8}}
		\end{subfigure}
		\hfill
		\begin{subfigure}[t]{0.32\textwidth}
			\centering
			\includegraphics[width=\linewidth]{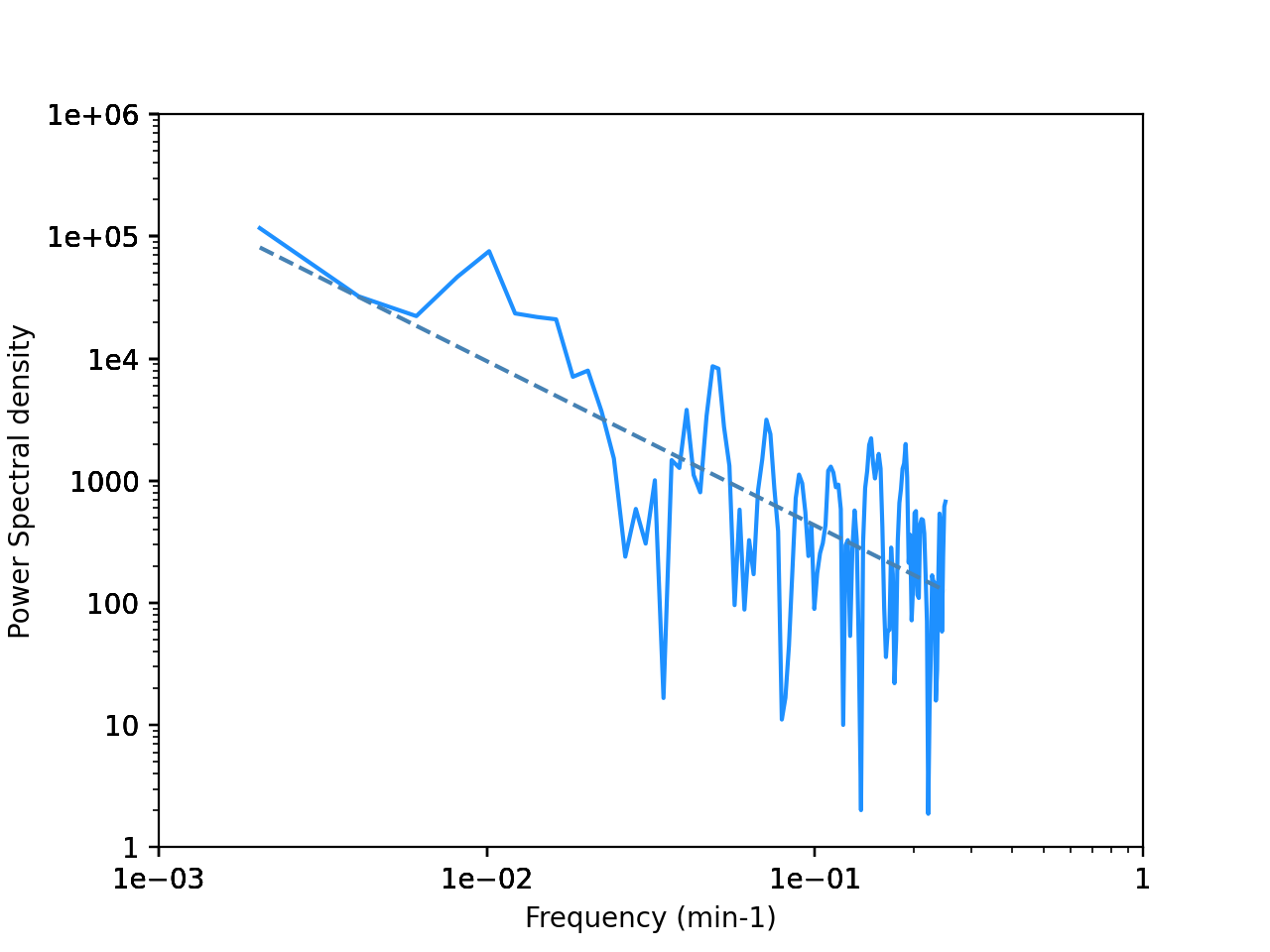}
			\caption{\textbf{Window 12}}
		\end{subfigure}
		
		\caption{\textbf{Local power spectral density (PSD) with regression (pseudopods).}
			Log--log representation of the PSD of the reconstructed driving signal for representative outer windows (2, 3, 4, 7, 8, 12). 
			A linear regression is performed over the selected frequency range, providing an estimate of the local scaling exponent $\beta$ in $S(\omega) \propto \omega^{-\beta}$. 
			The slopes remain broadly compatible with $\beta \approx 2$, while exhibiting inter-window variability attributable to finite-size effects at low frequencies and noise-dominated behavior at high frequencies.}
		\label{PSDWithRegressionLocauxSLEPseudopodsBiggestComponent}
	\end{figure*}
	

	\begin{figure*}[h!]
		\centering
		
		\begin{subfigure}[t]{0.32\textwidth}
			\centering
			\includegraphics[width=\linewidth]{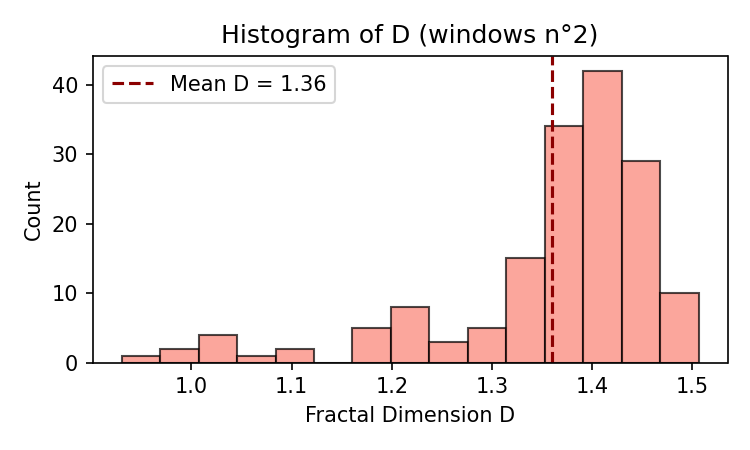}
			\caption{\textbf{Window 2}}
		\end{subfigure}
		\hfill
		\begin{subfigure}[t]{0.32\textwidth}
			\centering
			\includegraphics[width=\linewidth]{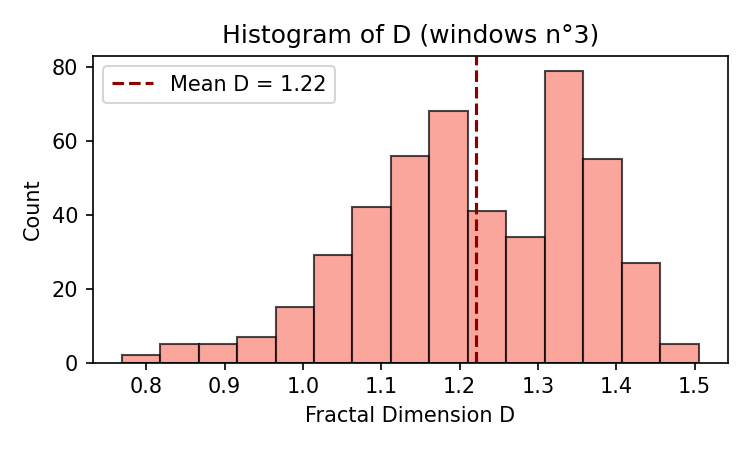}
			\caption{\textbf{Window 3}}
		\end{subfigure}
		\hfill
		\begin{subfigure}[t]{0.32\textwidth}
			\centering
			\includegraphics[width=\linewidth]{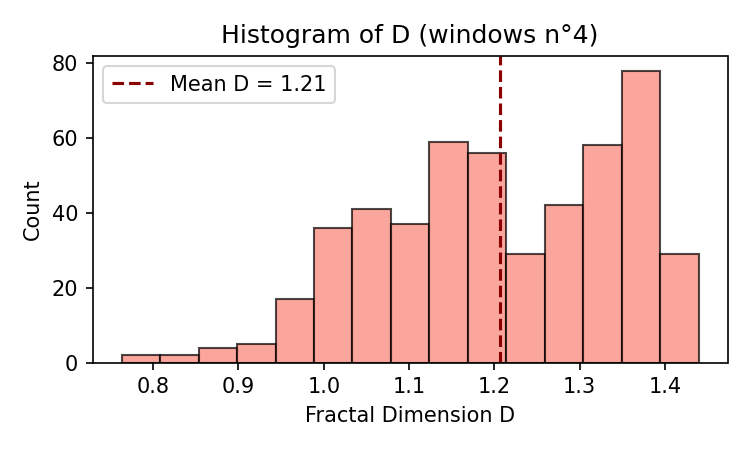}
			\caption{\textbf{Window 4}}
		\end{subfigure}
		
		\vspace{0.35cm}
		
		\begin{subfigure}[t]{0.32\textwidth}
			\centering
			\includegraphics[width=\linewidth]{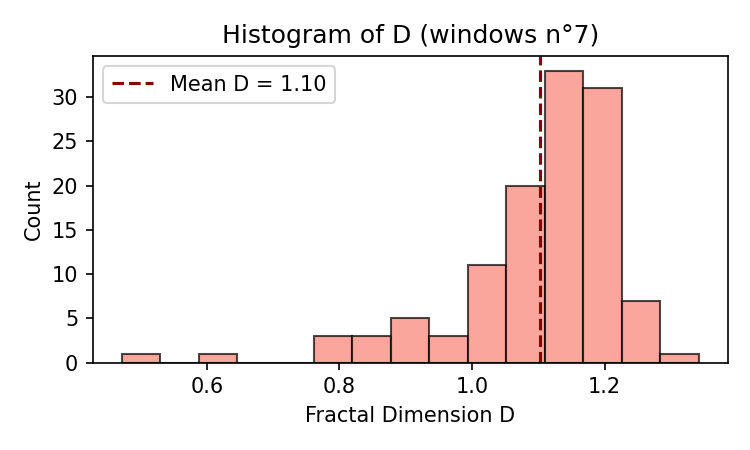}
			\caption{\textbf{Window 7}}
		\end{subfigure}
		\hfill
		\begin{subfigure}[t]{0.32\textwidth}
			\centering
			\includegraphics[width=\linewidth]{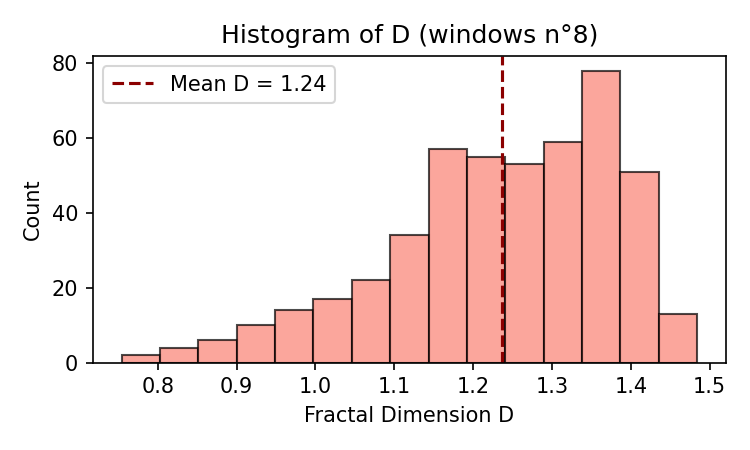}
			\caption{\textbf{Window 8}}
		\end{subfigure}
		\hfill
		\begin{subfigure}[t]{0.32\textwidth}
			\centering
			\includegraphics[width=\linewidth]{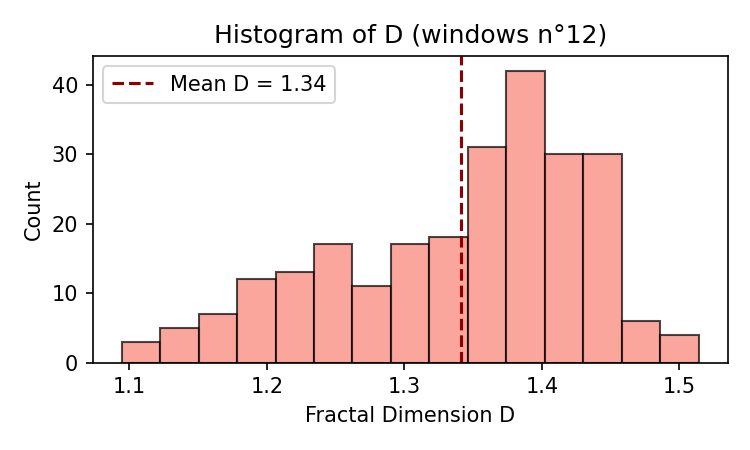}
			\caption{\textbf{Window 12}}
		\end{subfigure}
		
		\caption{\textbf{Local histograms of the fractal dimension (pseudopods).}
			Distributions of the estimated local fractal dimension computed on representative outer windows (2, 3, 4, 7, 8, 12). 
			The histograms exhibit moderate dispersion and window-to-window variability, reflecting geometric heterogeneity across the pseudopod organization.}
		\label{HistogrammesLocauxDimensionFractaleSLEPseudopodsBiggestComponent}
	\end{figure*}
	

	\begin{figure*}[h!]
		\centering
		
		\begin{subfigure}[t]{0.32\textwidth}
			\centering
			\includegraphics[width=\linewidth]{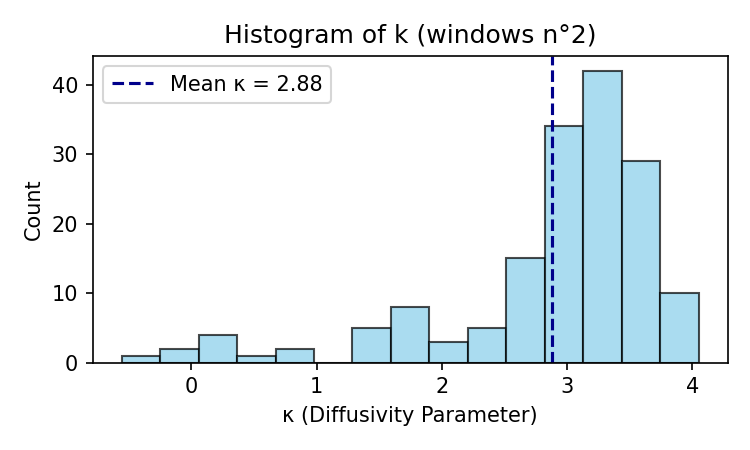}
			\caption{\textbf{Window 2}}
		\end{subfigure}
		\hfill
		\begin{subfigure}[t]{0.32\textwidth}
			\centering
			\includegraphics[width=\linewidth]{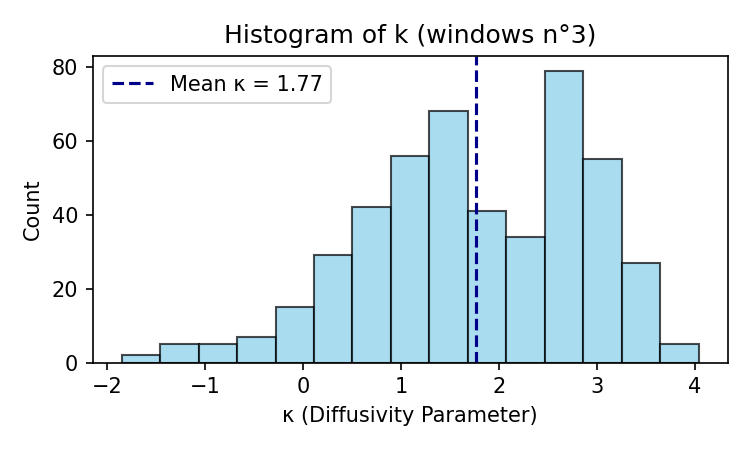}
			\caption{\textbf{Window 3}}
		\end{subfigure}
		\hfill
		\begin{subfigure}[t]{0.32\textwidth}
			\centering
			\includegraphics[width=\linewidth]{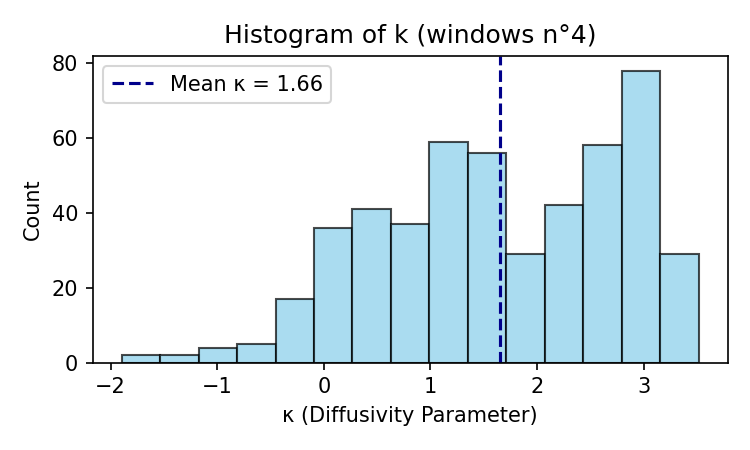}
			\caption{\textbf{Window 4}}
		\end{subfigure}
		
		\vspace{0.35cm}
		
		\begin{subfigure}[t]{0.32\textwidth}
			\centering
			\includegraphics[width=\linewidth]{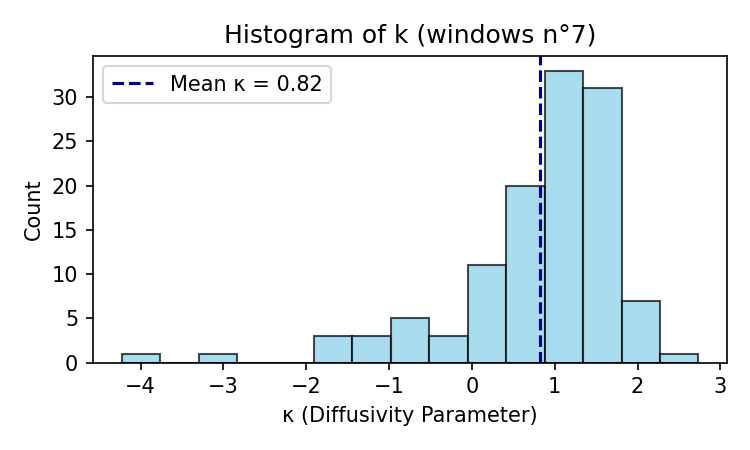}
			\caption{\textbf{Window 7}}
		\end{subfigure}
		\hfill
		\begin{subfigure}[t]{0.32\textwidth}
			\centering
			\includegraphics[width=\linewidth]{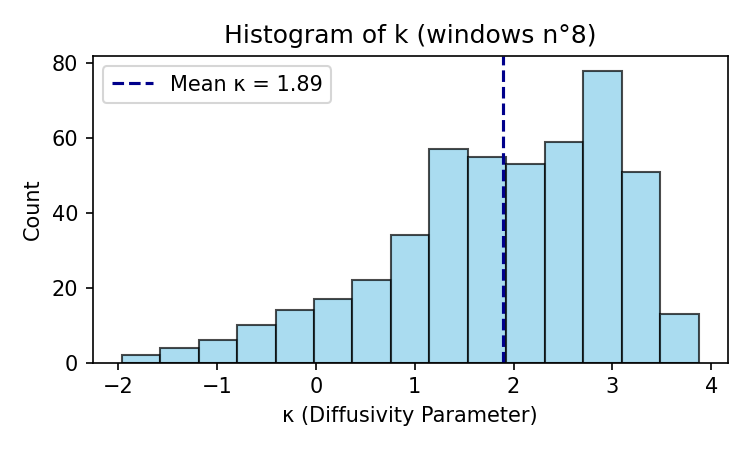}
			\caption{\textbf{Window 8}}
		\end{subfigure}
		\hfill
		\begin{subfigure}[t]{0.32\textwidth}
			\centering
			\includegraphics[width=\linewidth]{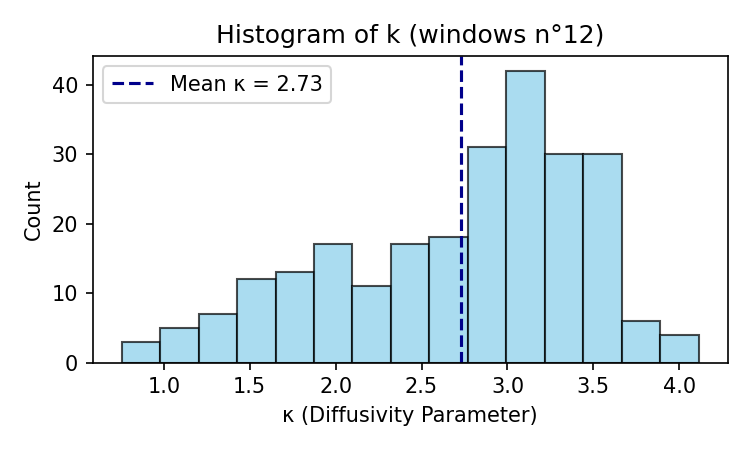}
			\caption{\textbf{Window 12}}
		\end{subfigure}
		
		\caption{\textbf{Local histograms of the diffusivity parameter $\kappa$ (pseudopods).}
			Estimated values of $\kappa$ for representative outer windows (2, 3, 4, 7, 8, 12). 
			The distributions exhibit noticeable dispersion and inter-window variability, reflecting scale-dependent heterogeneity and finite-size effects, while remaining compatible with a Brownian-type scaling regime over the explored ranges.}
		\label{HistogrammesLocauxKappaSLEPseudopodsBiggestComponent}
	\end{figure*}
	
	
	
	\begin{figure*}[h!]
		\centering
		
		\begin{subfigure}[t]{0.48\textwidth}
			\centering
			\includegraphics[width=\linewidth]{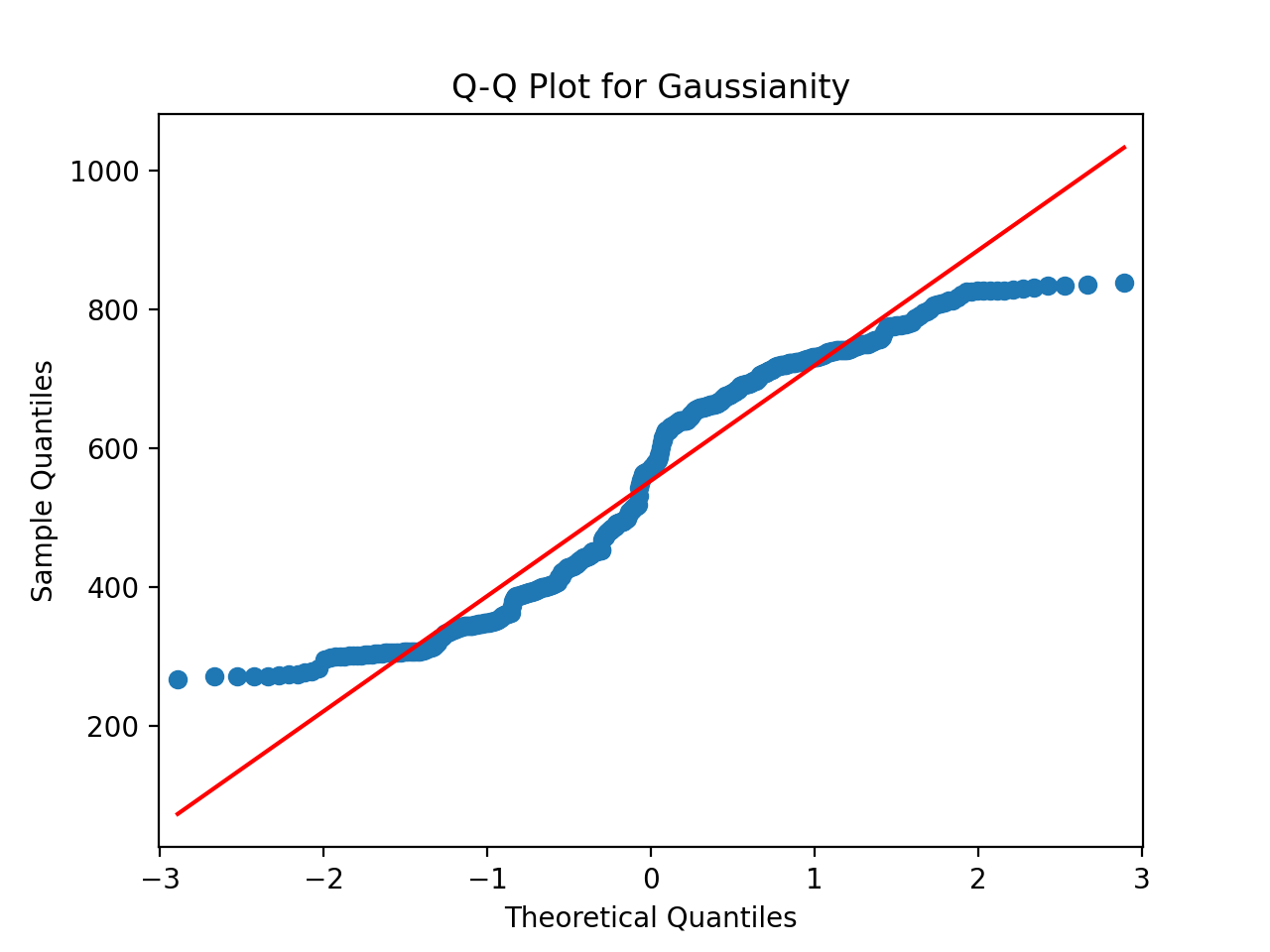}
			\caption{\textbf{Global Q-Q plot}}
		\end{subfigure}
		\hfill
		\begin{subfigure}[t]{0.48\textwidth}
			\centering
			\includegraphics[width=\linewidth]{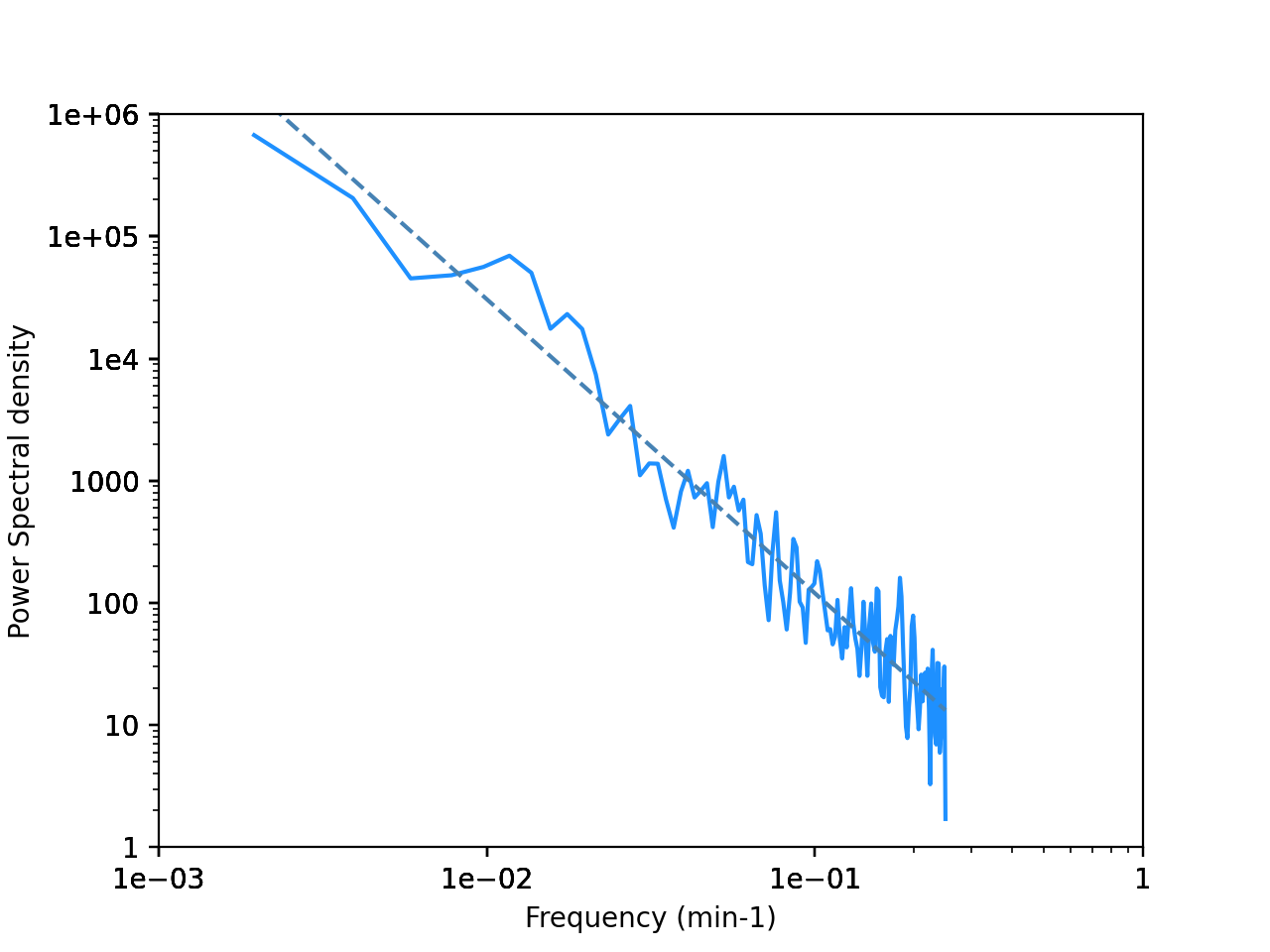}
			\caption{\textbf{Global PSD with regression}}
		\end{subfigure}
		
		\vspace{0.4cm}
		
		\begin{subfigure}[t]{0.60\textwidth}
			\centering
			\includegraphics[width=\linewidth]{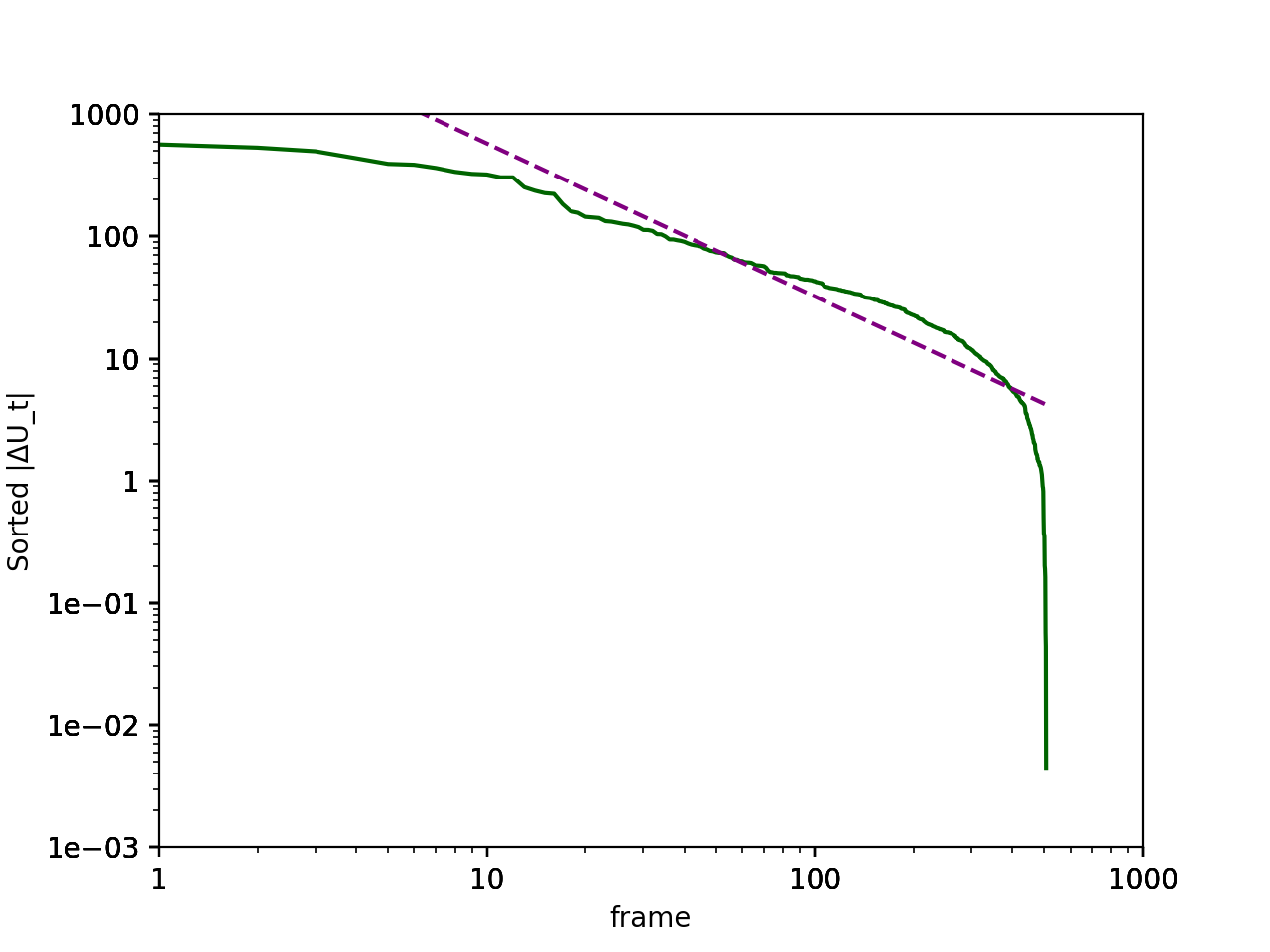}
			\caption{\textbf{Hill-type tail diagnostic.}
				Log–log plot of ordered absolute increments $\Delta U_t$. 
				The absence of a stable scaling plateau and the rapid decay of the tail 
				are incompatible with heavy-tailed behavior and support Gaussian statistics.}
		\end{subfigure}
		
		\vspace{0.4cm}
		
		\begin{subfigure}[t]{0.48\textwidth}
			\centering
			\includegraphics[width=\linewidth]{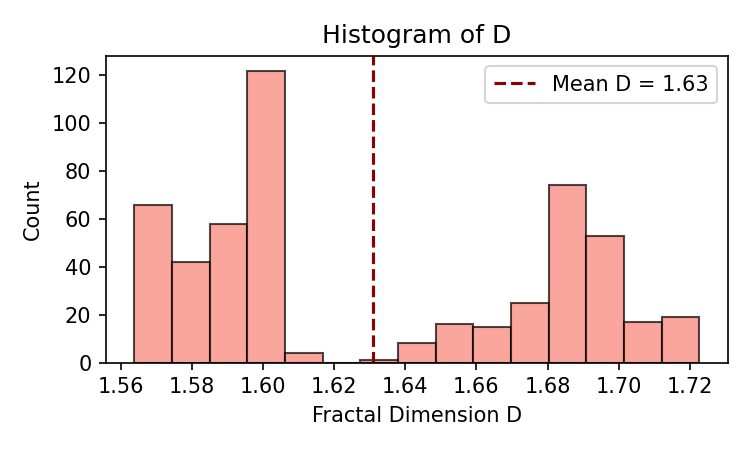}
			\caption{\textbf{Global fractal dimension}}
		\end{subfigure}
		\hfill
		\begin{subfigure}[t]{0.48\textwidth}
			\centering
			\includegraphics[width=\linewidth]{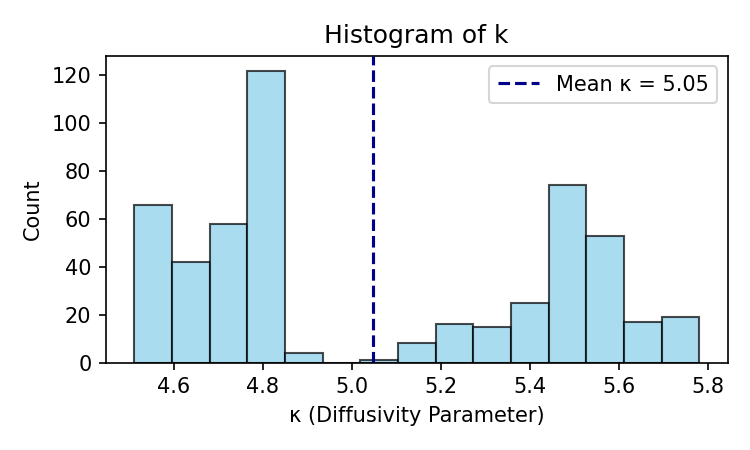}
			\caption{\textbf{Global diffusivity parameter $\kappa$}}
		\end{subfigure}
		
		\caption{\textbf{Global statistical diagnostics for the largest connected component of the extracted network.}
			Q-Q plot, power spectral density with regression, Hill-type tail analysis, and histograms of the global fractal dimension and diffusivity parameter $\kappa$. 
			Here, \emph{global} refers to the reconstruction performed on the largest connected component of the full network.}
		\label{QQPlotsPSDEtHistogrammesGlobauxSLENetworkBiggestComponent}
	\end{figure*}
	

	
	\begin{figure*}[h!]
		\centering
		
		\begin{subfigure}[t]{0.32\textwidth}
			\centering
			\includegraphics[width=\linewidth]{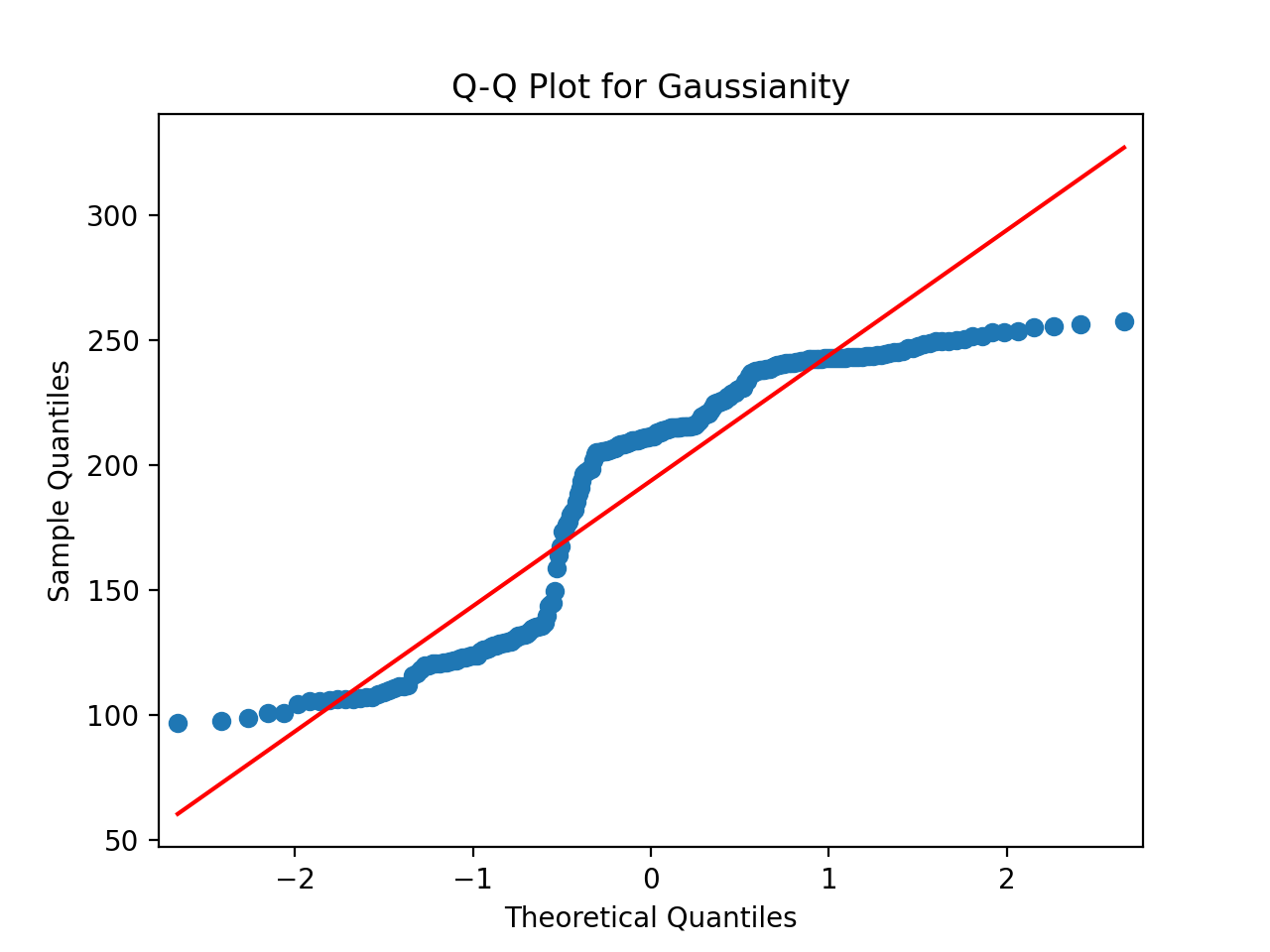}
			\caption{\textbf{Window 2}}
		\end{subfigure}
		\hfill
		\begin{subfigure}[t]{0.32\textwidth}
			\centering
			\includegraphics[width=\linewidth]{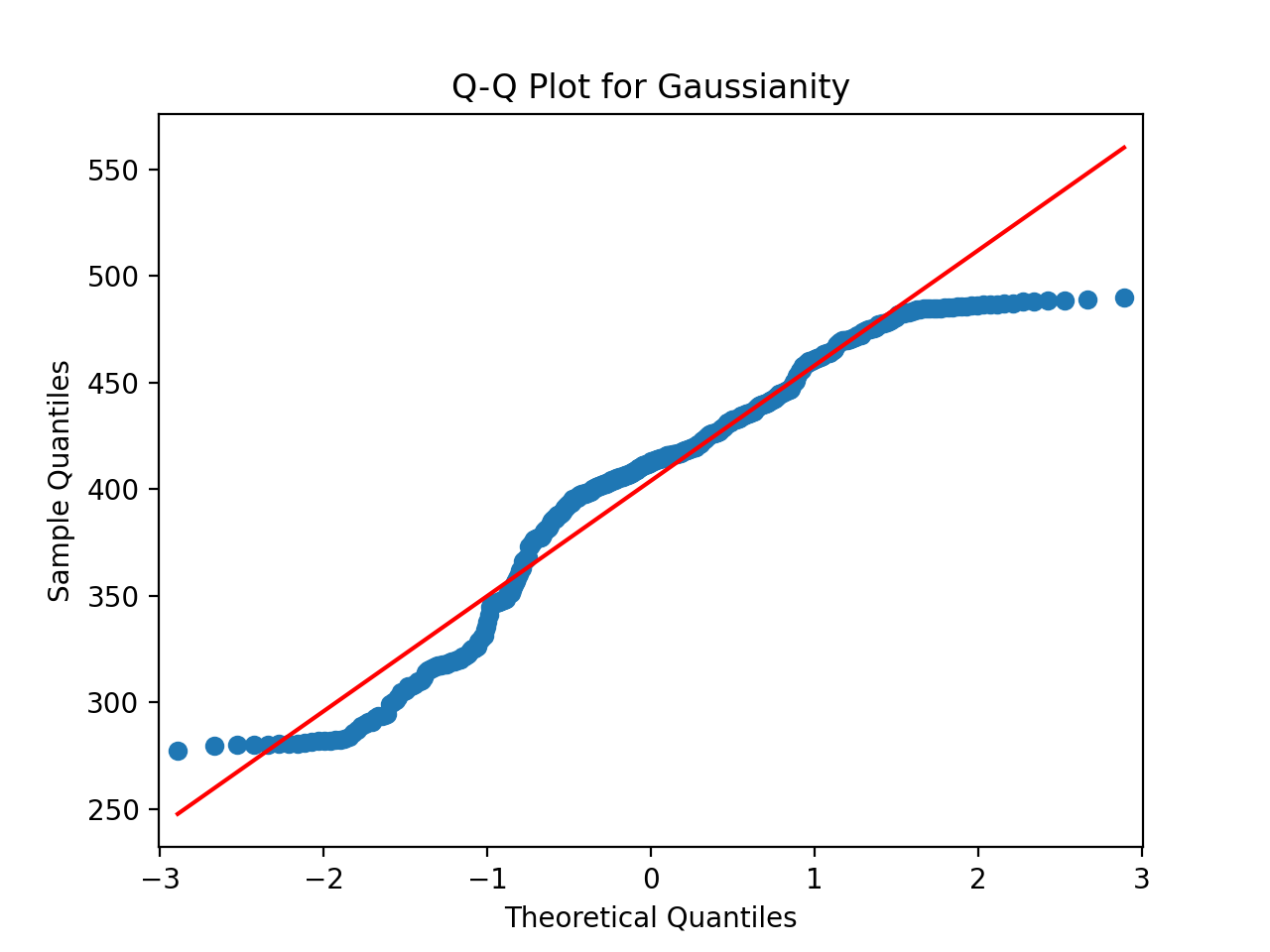}
			\caption{\textbf{Window 3}}
		\end{subfigure}
		\hfill
		\begin{subfigure}[t]{0.32\textwidth}
			\centering
			\includegraphics[width=\linewidth]{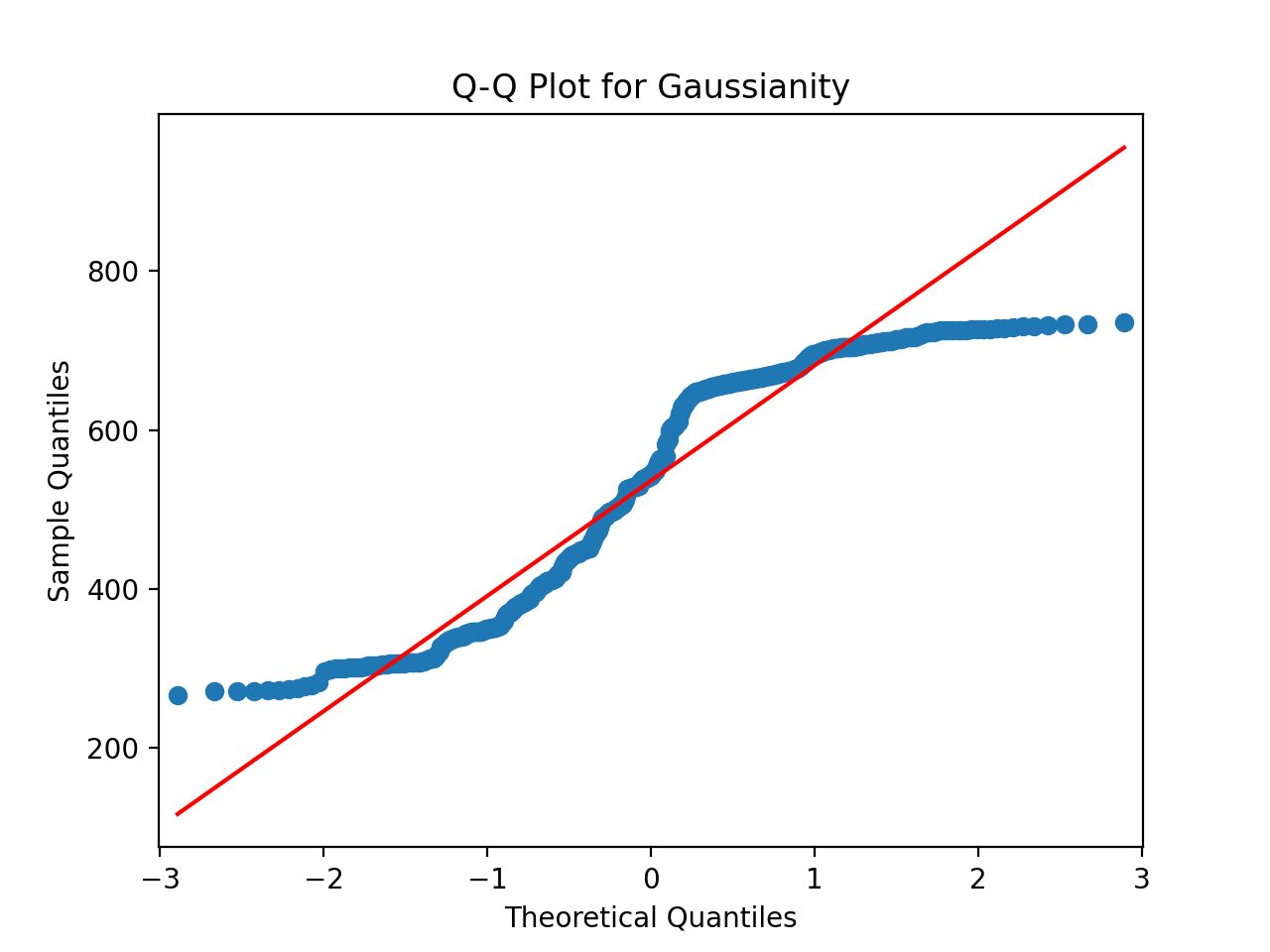}
			\caption{\textbf{Window 4}}
		\end{subfigure}
		
		\vspace{0.35cm}
		
		\begin{subfigure}[t]{0.32\textwidth}
			\centering
			\includegraphics[width=\linewidth]{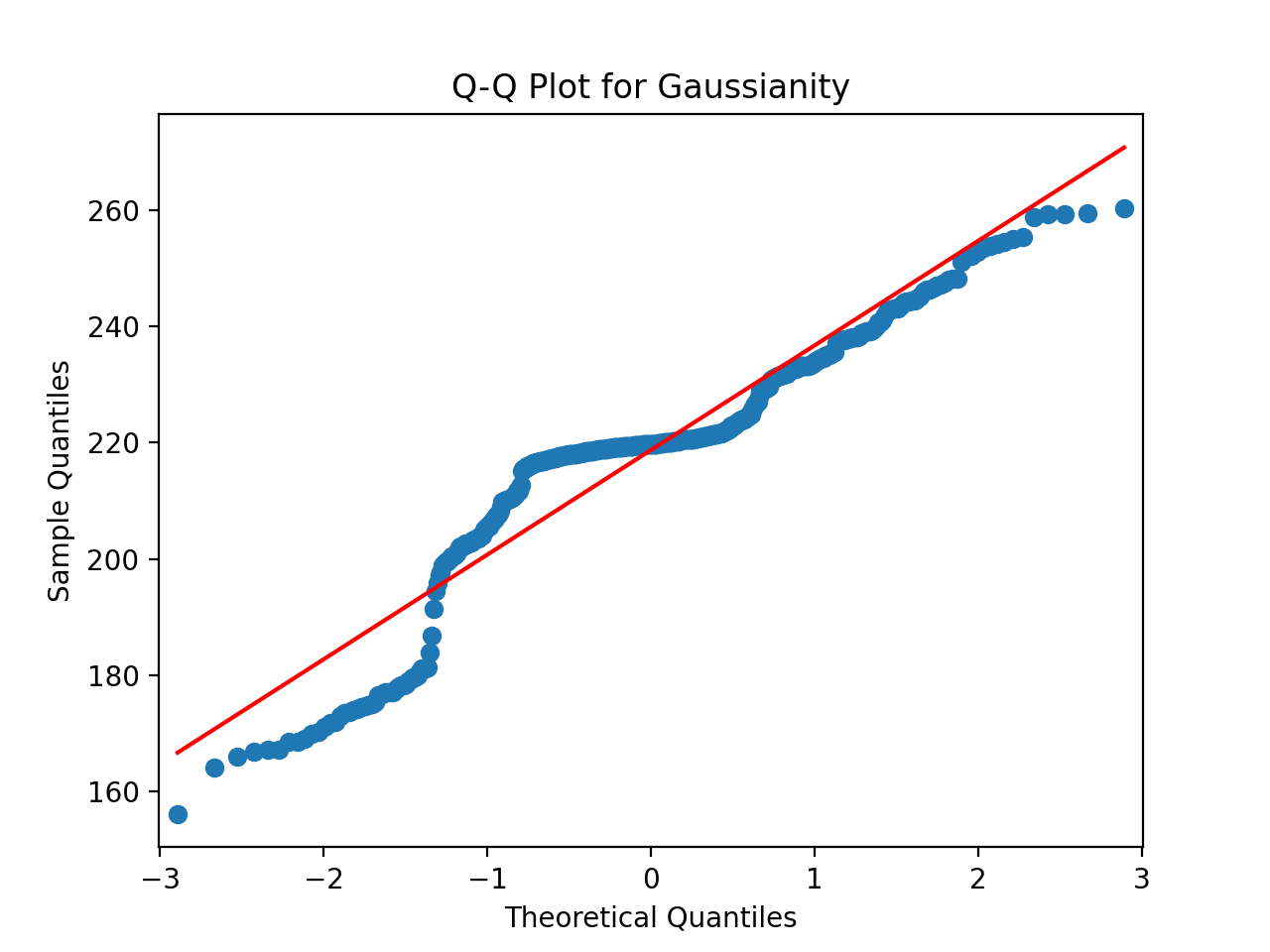}
			\caption{\textbf{Window 7}}
		\end{subfigure}
		\hfill
		\begin{subfigure}[t]{0.32\textwidth}
			\centering
			\includegraphics[width=\linewidth]{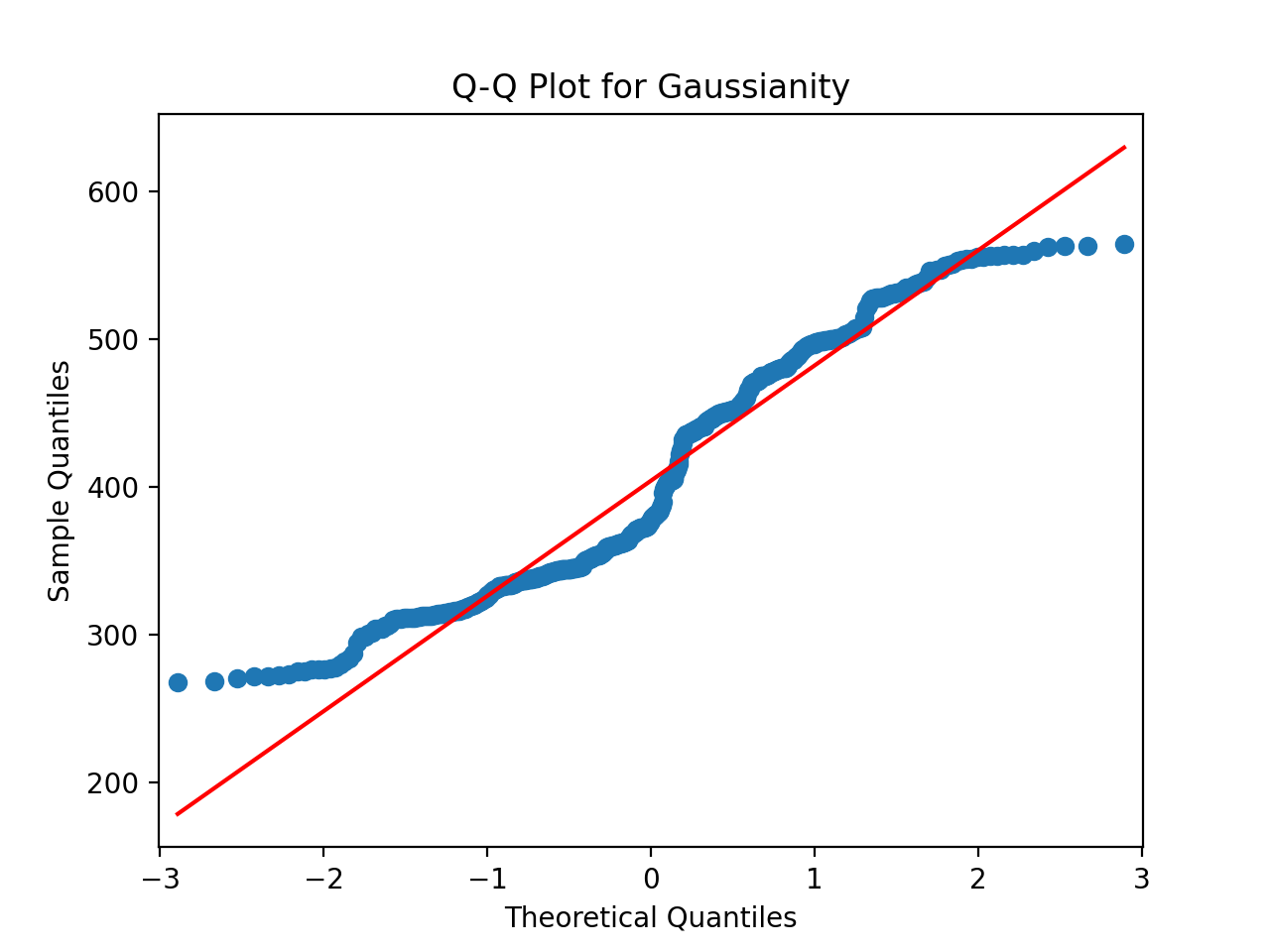}
			\caption{\textbf{Window 8}}
		\end{subfigure}
		\hfill
		\begin{subfigure}[t]{0.32\textwidth}
			\centering
			\includegraphics[width=\linewidth]{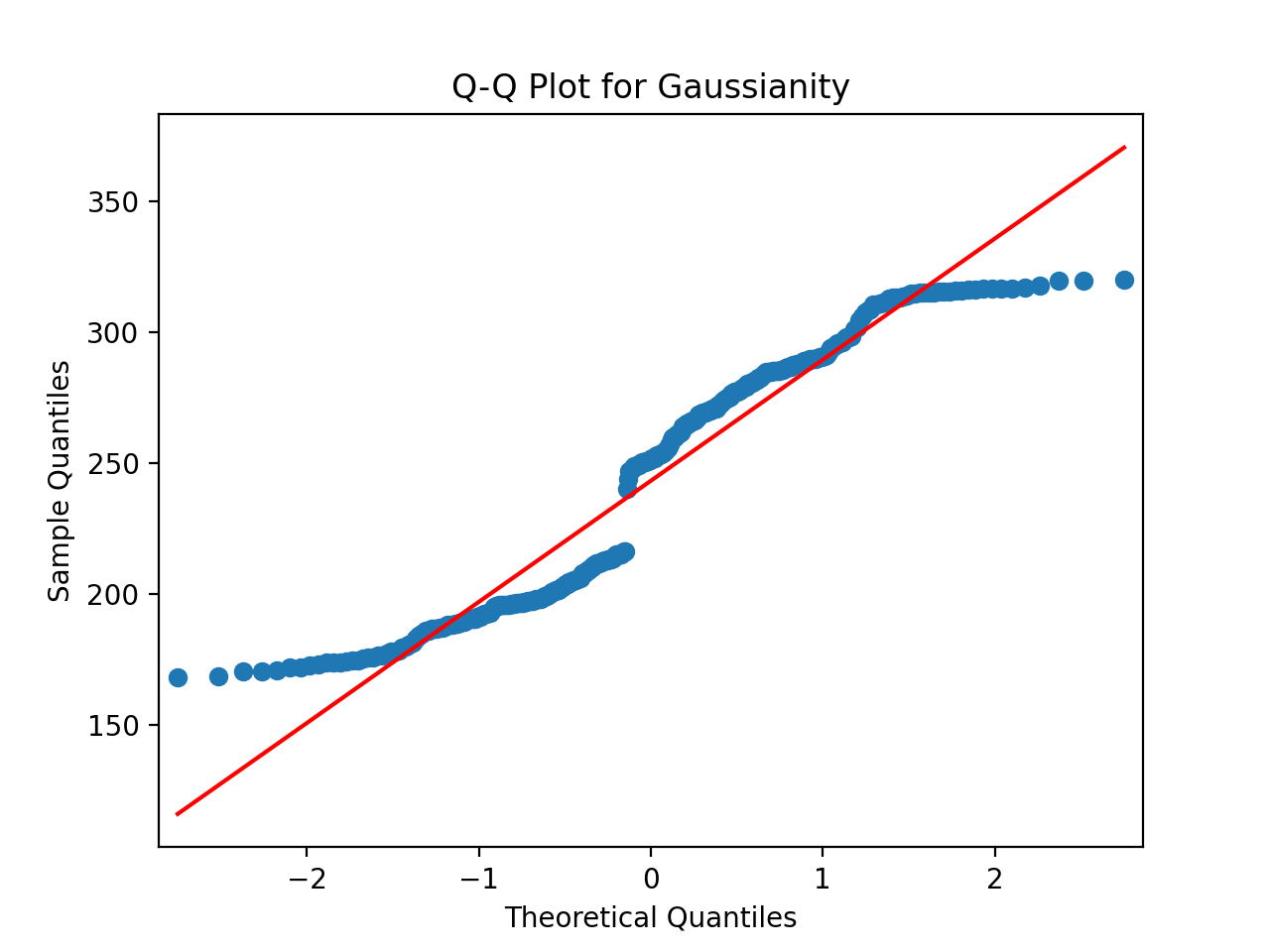}
			\caption{\textbf{Window 12}}
		\end{subfigure}
		
		\caption{\textbf{Local Q-Q plots of the driving function increments $\Delta U_t$ (network).}
			Representative outer spatial windows (2, 3, 4, 7, 8, 12) computed on the largest connected component of the extracted network. 
			As for the pseudopod component, the empirical quantiles closely follow the theoretical Gaussian quantiles, with only mild deviations at extreme tails.}
		\label{QQPlotsLocauxSLENetworkBiggestComponent}
	\end{figure*}

	
	\begin{figure*}[h!]
		\centering
		
		\begin{subfigure}[t]{0.32\textwidth}
			\centering
			\includegraphics[width=\linewidth]{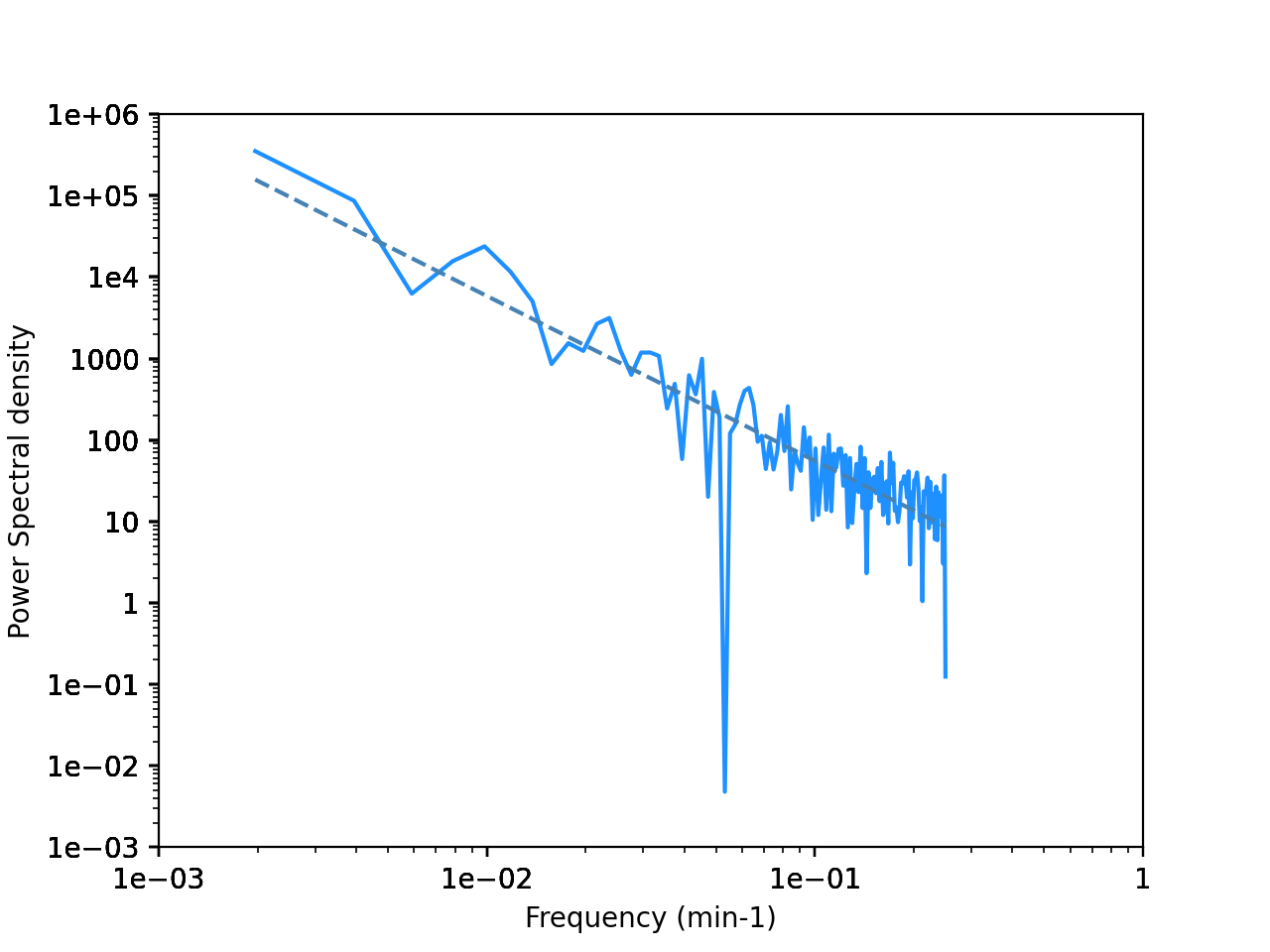}
			\caption{\textbf{Window 2}}
		\end{subfigure}
		\hfill
		\begin{subfigure}[t]{0.32\textwidth}
			\centering
			\includegraphics[width=\linewidth]{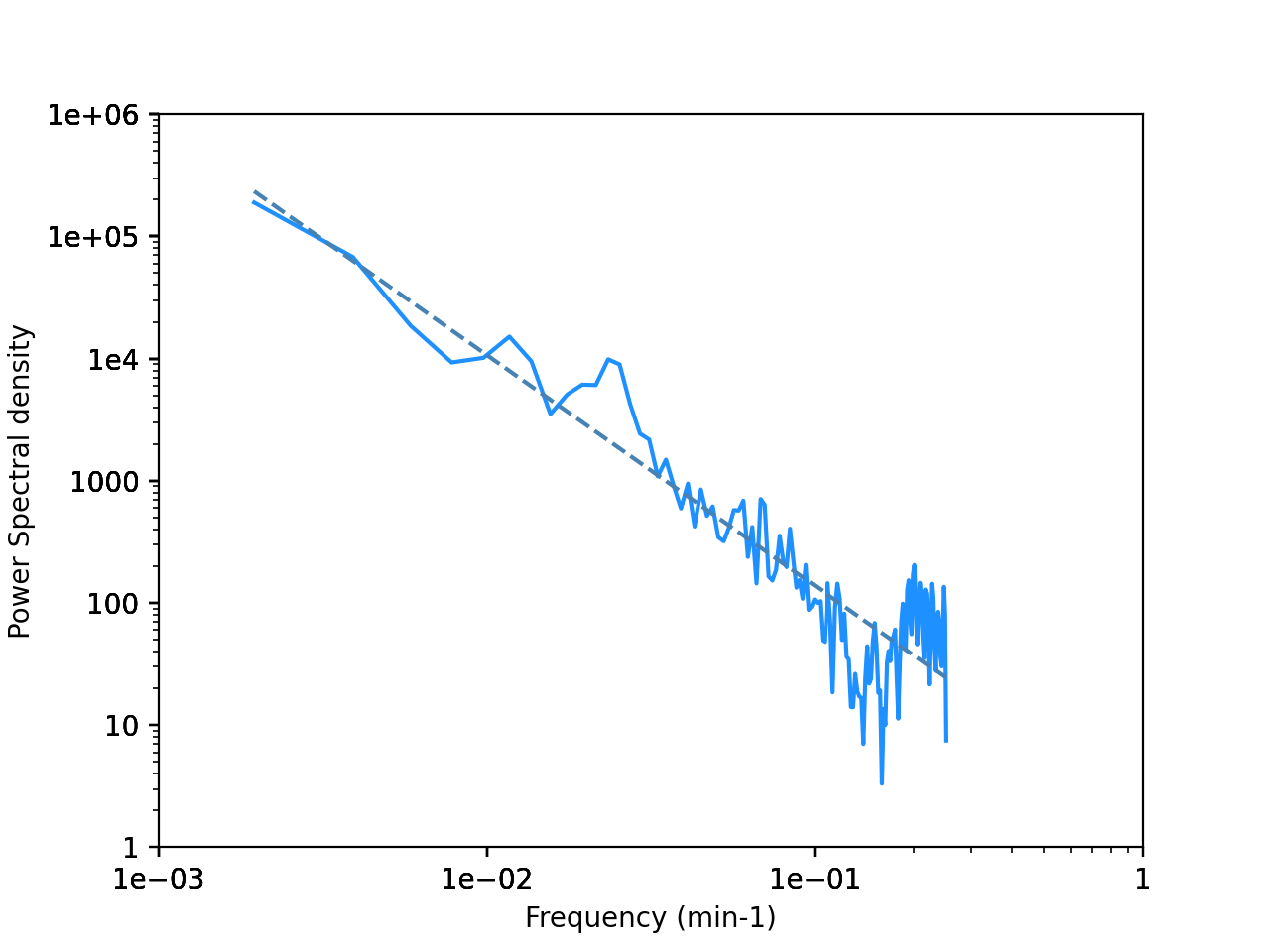}
			\caption{\textbf{Window 3}}
		\end{subfigure}
		\hfill
		\begin{subfigure}[t]{0.32\textwidth}
			\centering
			\includegraphics[width=\linewidth]{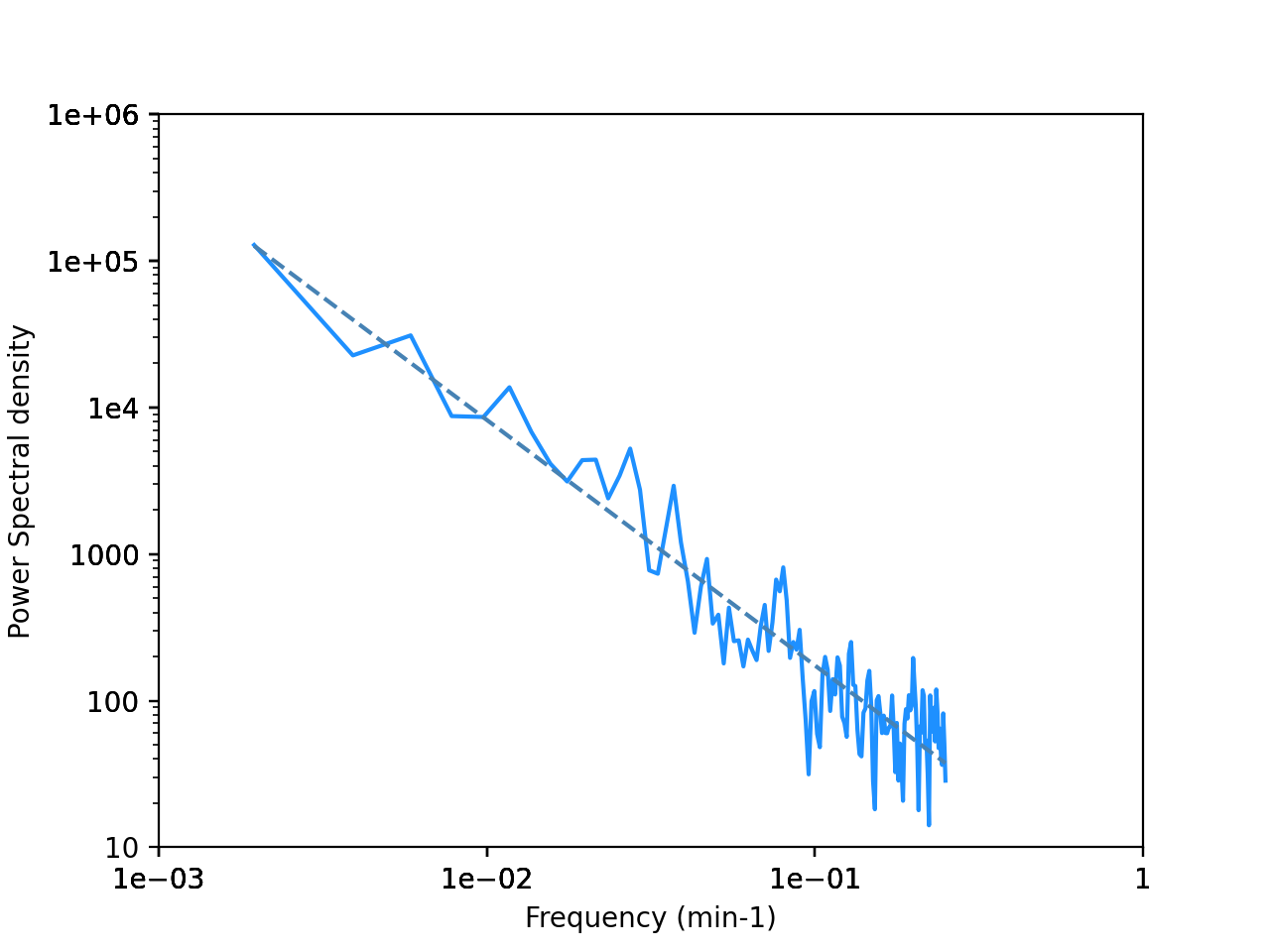}
			\caption{\textbf{Window 4}}
		\end{subfigure}
		
		\vspace{0.35cm}
		
		\begin{subfigure}[t]{0.32\textwidth}
			\centering
			\includegraphics[width=\linewidth]{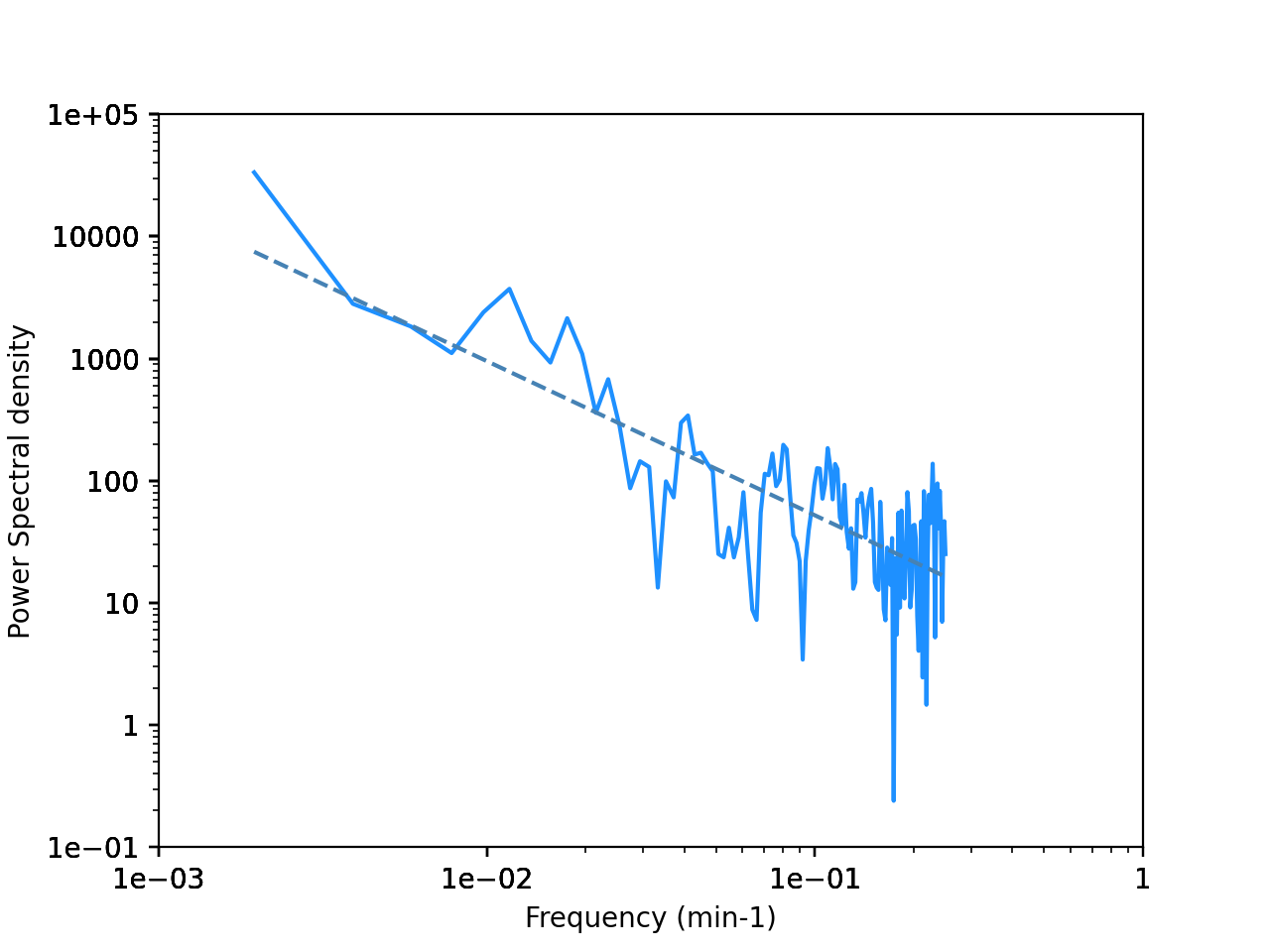}
			\caption{\textbf{Window 7}}
		\end{subfigure}
		\hfill
		\begin{subfigure}[t]{0.32\textwidth}
			\centering
			\includegraphics[width=\linewidth]{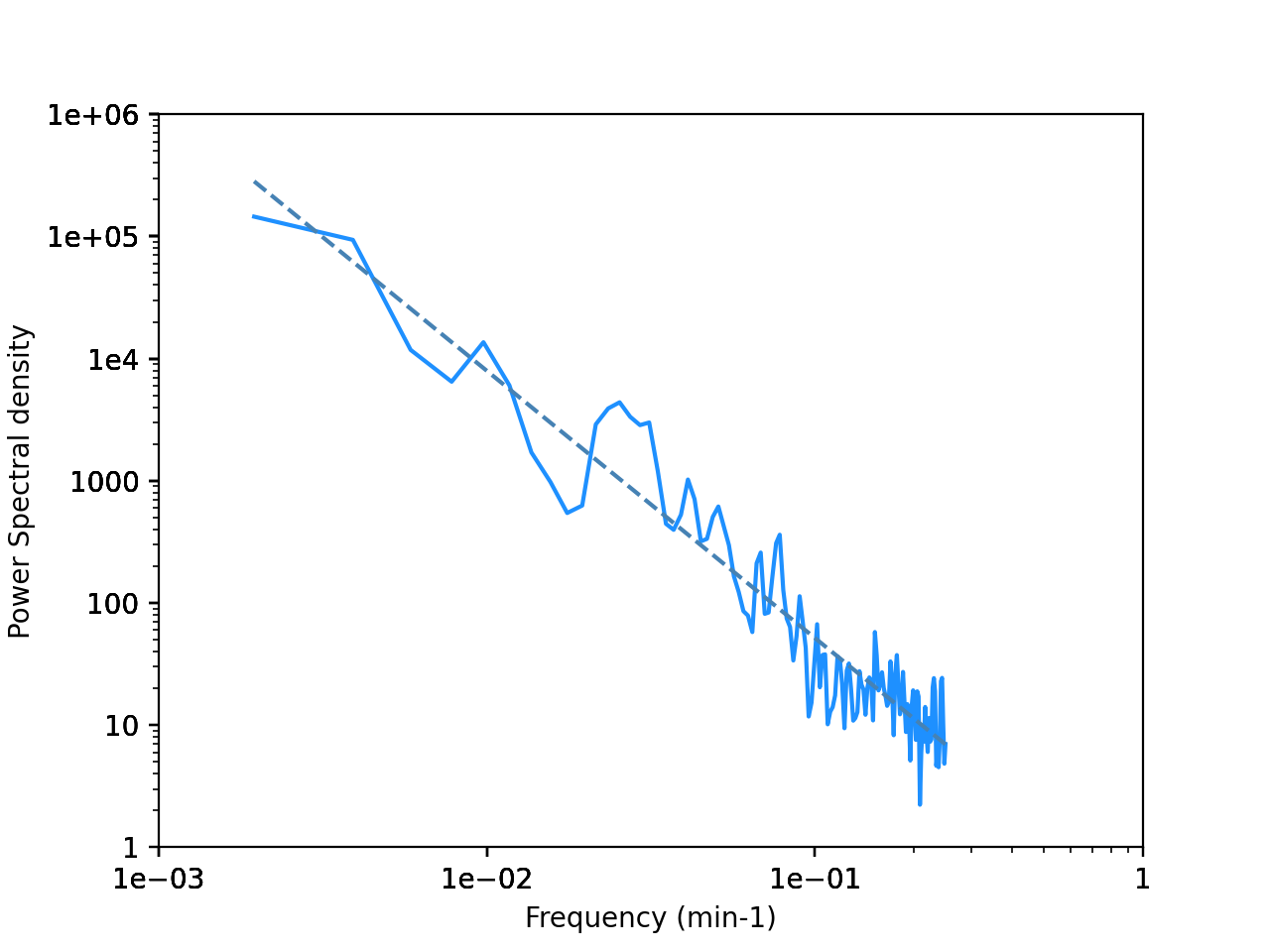}
			\caption{\textbf{Window 8}}
		\end{subfigure}
		\hfill
		\begin{subfigure}[t]{0.32\textwidth}
			\centering
			\includegraphics[width=\linewidth]{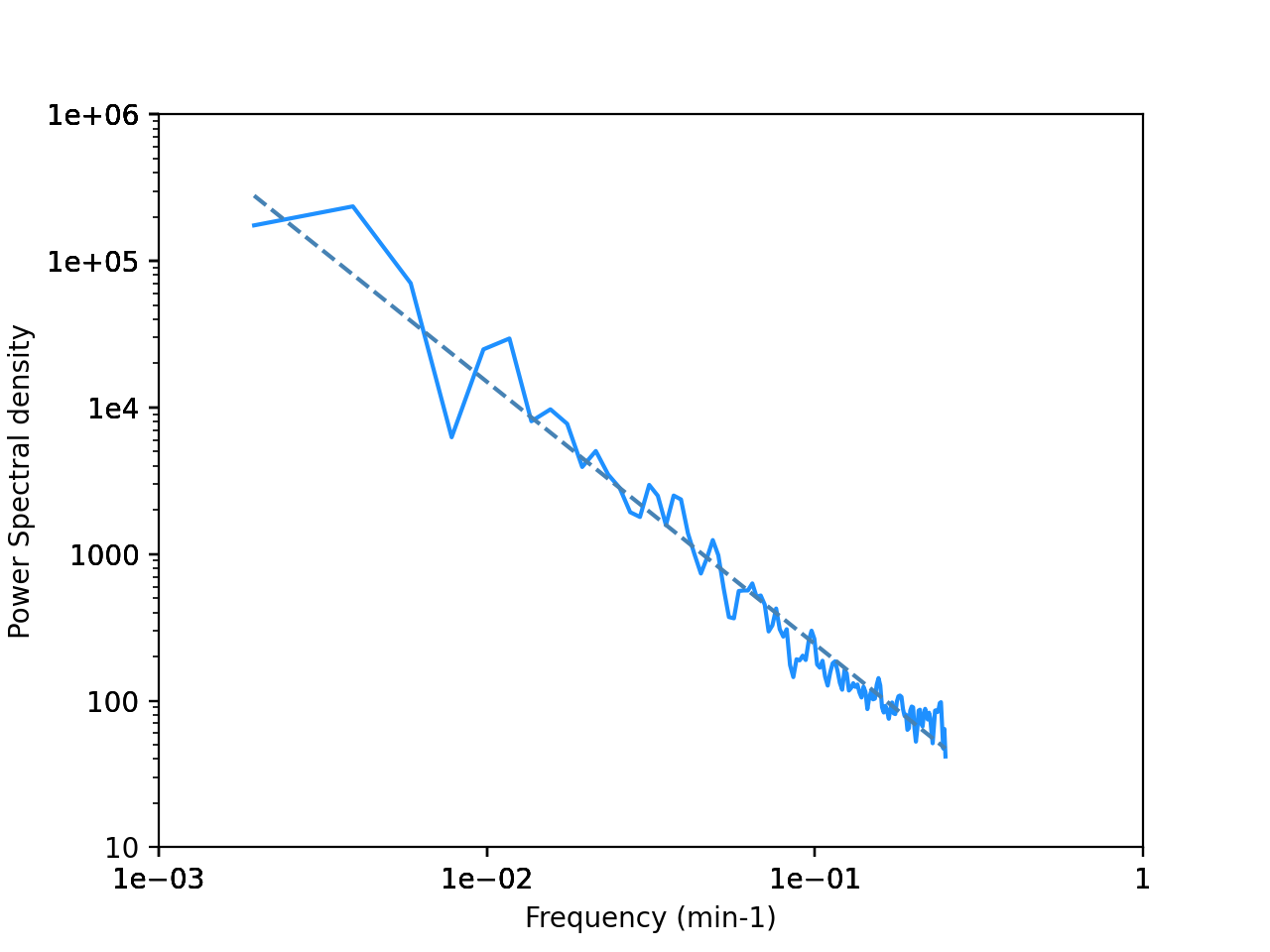}
			\caption{\textbf{Window 12}}
		\end{subfigure}
		
		\caption{\textbf{Local power spectral density (PSD) diagnostics for the network component.}
			Log--log PSD plots with linear regression for representative outer windows (2, 3, 4, 7, 8, 12), computed on the largest connected component of the extracted network. 
			In each window, the fitted slope provides an estimate of the local scaling exponent $\beta$ in $S(\omega)\propto\omega^{-\beta}$, which remains broadly compatible with Brownian-type scaling ($\beta\approx 2$), while exhibiting inter-window variability due to structural heterogeneity and finite-size effects.}
		\label{PSDWithRegressionLocauxSLENetworkBiggestComponent}
	\end{figure*}
	
	

	\begin{figure*}[h!]
		\centering
		
		\begin{subfigure}[t]{0.32\textwidth}
			\centering
			\includegraphics[width=\linewidth]{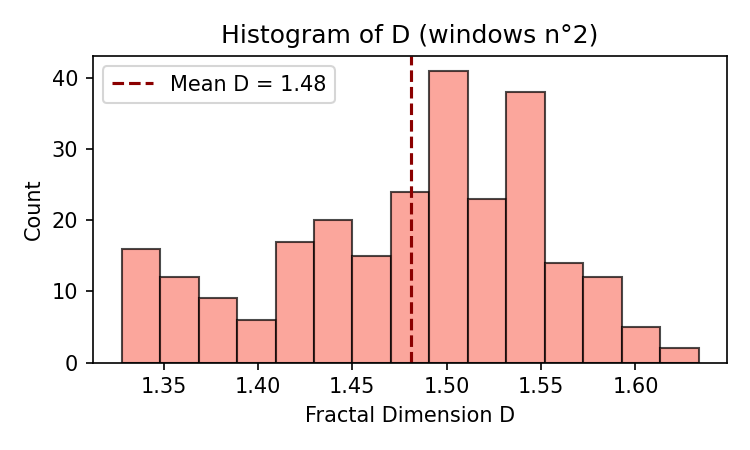}
			\caption{\textbf{Window 2}}
		\end{subfigure}
		\hfill
		\begin{subfigure}[t]{0.32\textwidth}
			\centering
			\includegraphics[width=\linewidth]{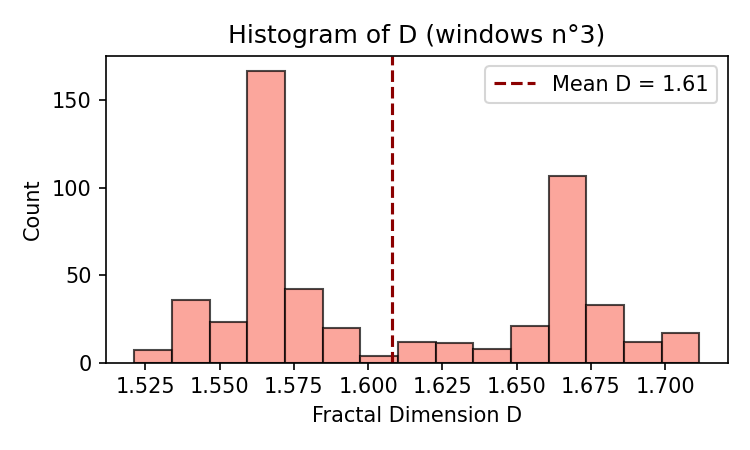}
			\caption{\textbf{Window 3}}
		\end{subfigure}
		\hfill
		\begin{subfigure}[t]{0.32\textwidth}
			\centering
			\includegraphics[width=\linewidth]{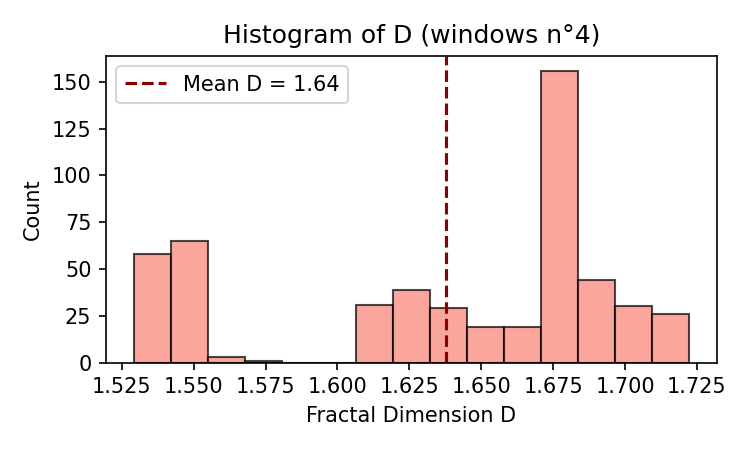}
			\caption{\textbf{Window 4}}
		\end{subfigure}
		
		\vspace{0.35cm}
		
		\begin{subfigure}[t]{0.32\textwidth}
			\centering
			\includegraphics[width=\linewidth]{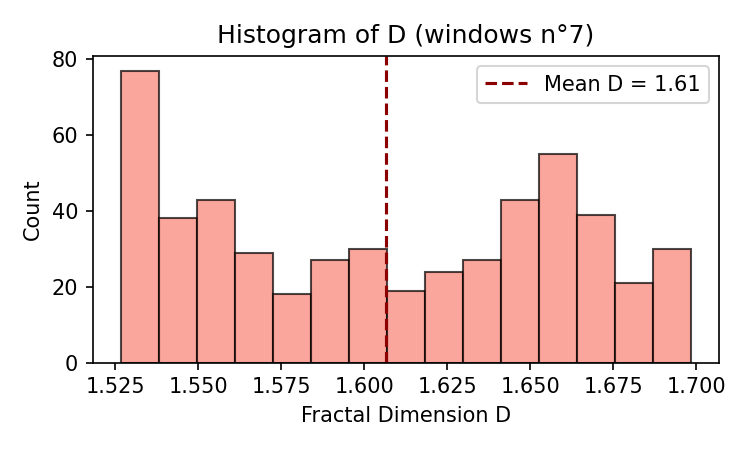}
			\caption{\textbf{Window 7}}
		\end{subfigure}
		\hfill
		\begin{subfigure}[t]{0.32\textwidth}
			\centering
			\includegraphics[width=\linewidth]{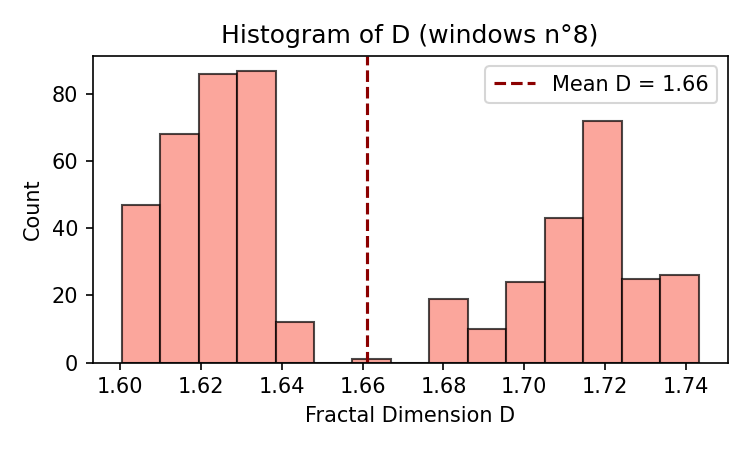}
			\caption{\textbf{Window 8}}
		\end{subfigure}
		\hfill
		\begin{subfigure}[t]{0.32\textwidth}
			\centering
			\includegraphics[width=\linewidth]{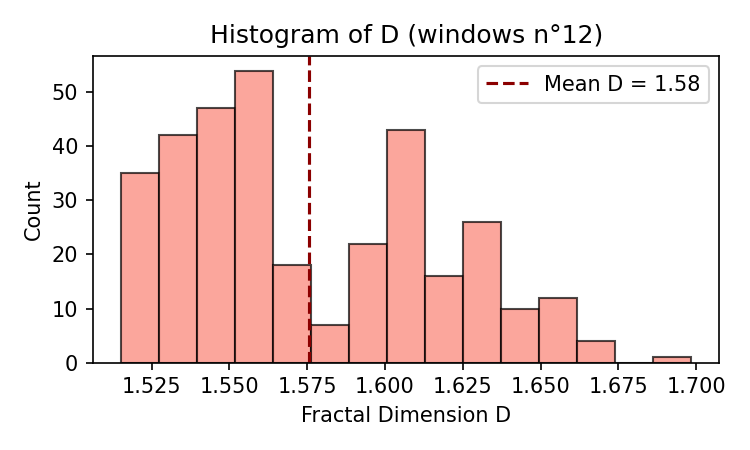}
			\caption{\textbf{Window 12}}
		\end{subfigure}
		
		\caption{\textbf{Local histograms of the fractal dimension (network component).}
			Distributions of the estimated local fractal dimension computed on representative outer windows (2, 3, 4, 7, 8, 12), for the largest connected component of the extracted network. 
			Compared with the pseudopod level, the distributions are shifted toward higher values and remain relatively stable across windows, suggesting a denser and more homogeneous multiscale organization of the transport architecture.}
		\label{HistogrammesLocauxDimensionFractaleSLENetworkBiggestComponent}
	\end{figure*}
	
	

	\begin{figure*}[h!]
		\centering
		
		\begin{subfigure}[t]{0.32\textwidth}
			\centering
			\includegraphics[width=\linewidth]{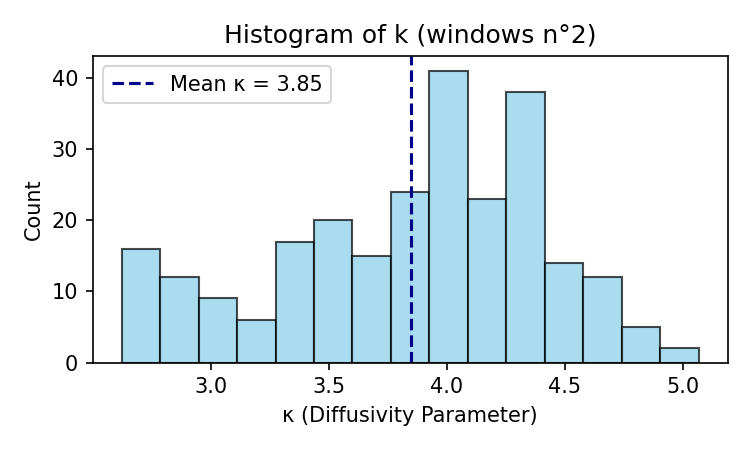}
			\caption{\textbf{Window 2}}
		\end{subfigure}
		\hfill
		\begin{subfigure}[t]{0.32\textwidth}
			\centering
			\includegraphics[width=\linewidth]{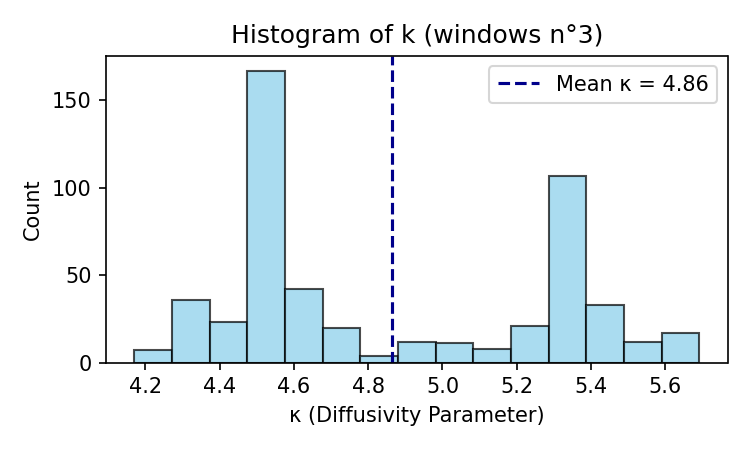}
			\caption{\textbf{Window 3}}
		\end{subfigure}
		\hfill
		\begin{subfigure}[t]{0.32\textwidth}
			\centering
			\includegraphics[width=\linewidth]{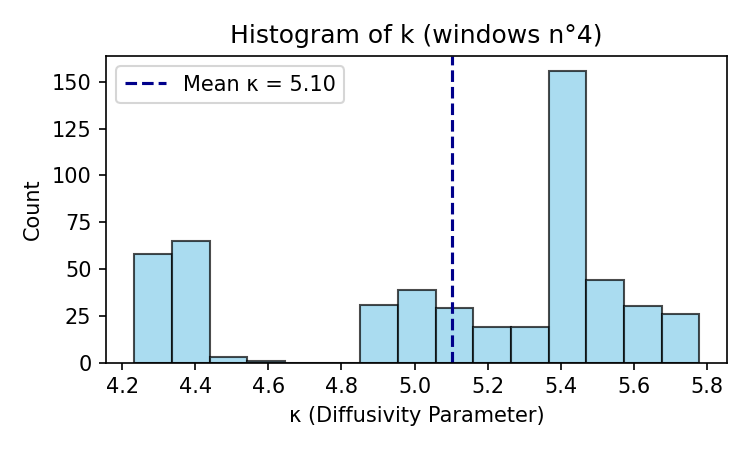}
			\caption{\textbf{Window 4}}
		\end{subfigure}
		
		\vspace{0.35cm}
		
		\begin{subfigure}[t]{0.32\textwidth}
			\centering
			\includegraphics[width=\linewidth]{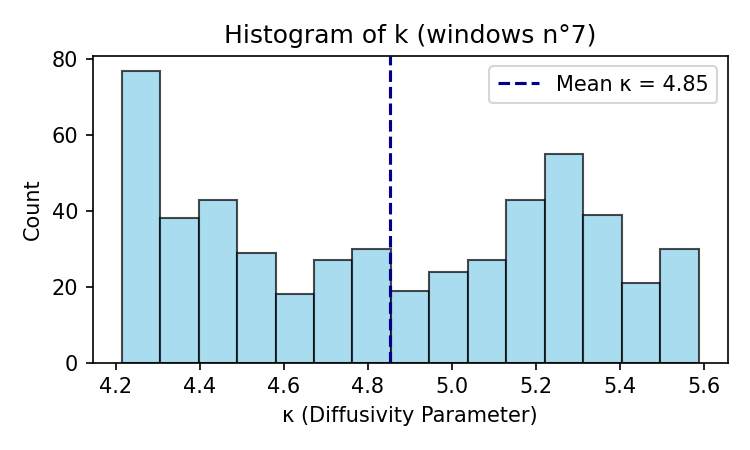}
			\caption{\textbf{Window 7}}
		\end{subfigure}
		\hfill
		\begin{subfigure}[t]{0.32\textwidth}
			\centering
			\includegraphics[width=\linewidth]{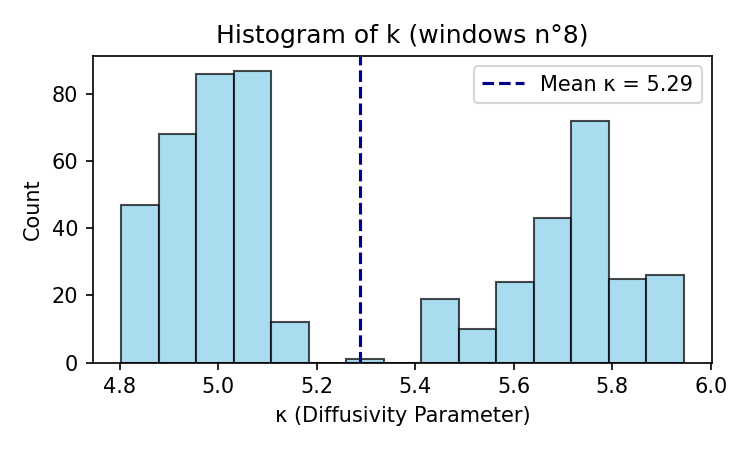}
			\caption{\textbf{Window 8}}
		\end{subfigure}
		\hfill
		\begin{subfigure}[t]{0.32\textwidth}
			\centering
			\includegraphics[width=\linewidth]{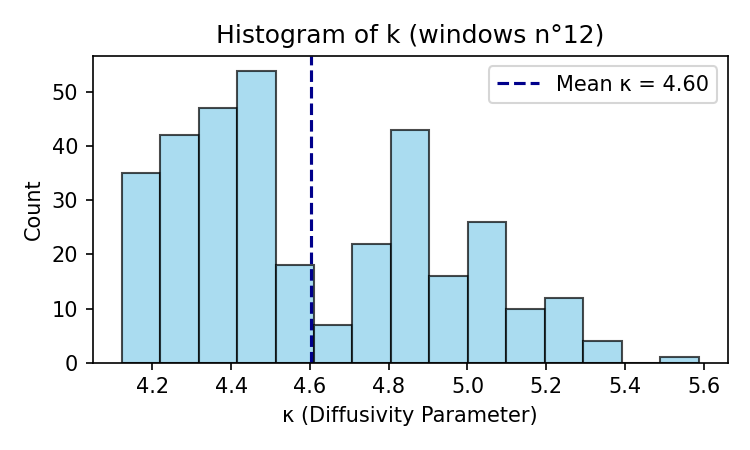}
			\caption{\textbf{Window 12}}
		\end{subfigure}
		
		\caption{\textbf{Local histograms of the diffusivity parameter $\kappa$ (network component).}
			Estimated values of $\kappa$ for representative outer windows (2, 3, 4, 7, 8, 12), computed on the largest connected component of the extracted network. 
			Compared with the pseudopod-level estimates, the distributions are shifted toward larger values of $\kappa$ and remain relatively stable across windows, suggesting a more coherent and homogeneous stochastic organization of the transport backbone across scales, despite residual dispersion and finite-size effects.}
		\label{HistogrammesLocauxKappaSLENetworkBiggestComponent}
	\end{figure*}
	
	

	\begin{figure*}[h!]
		\centering
		
		\begin{subfigure}[t]{0.48\textwidth}
			\centering
			\includegraphics[width=\linewidth]{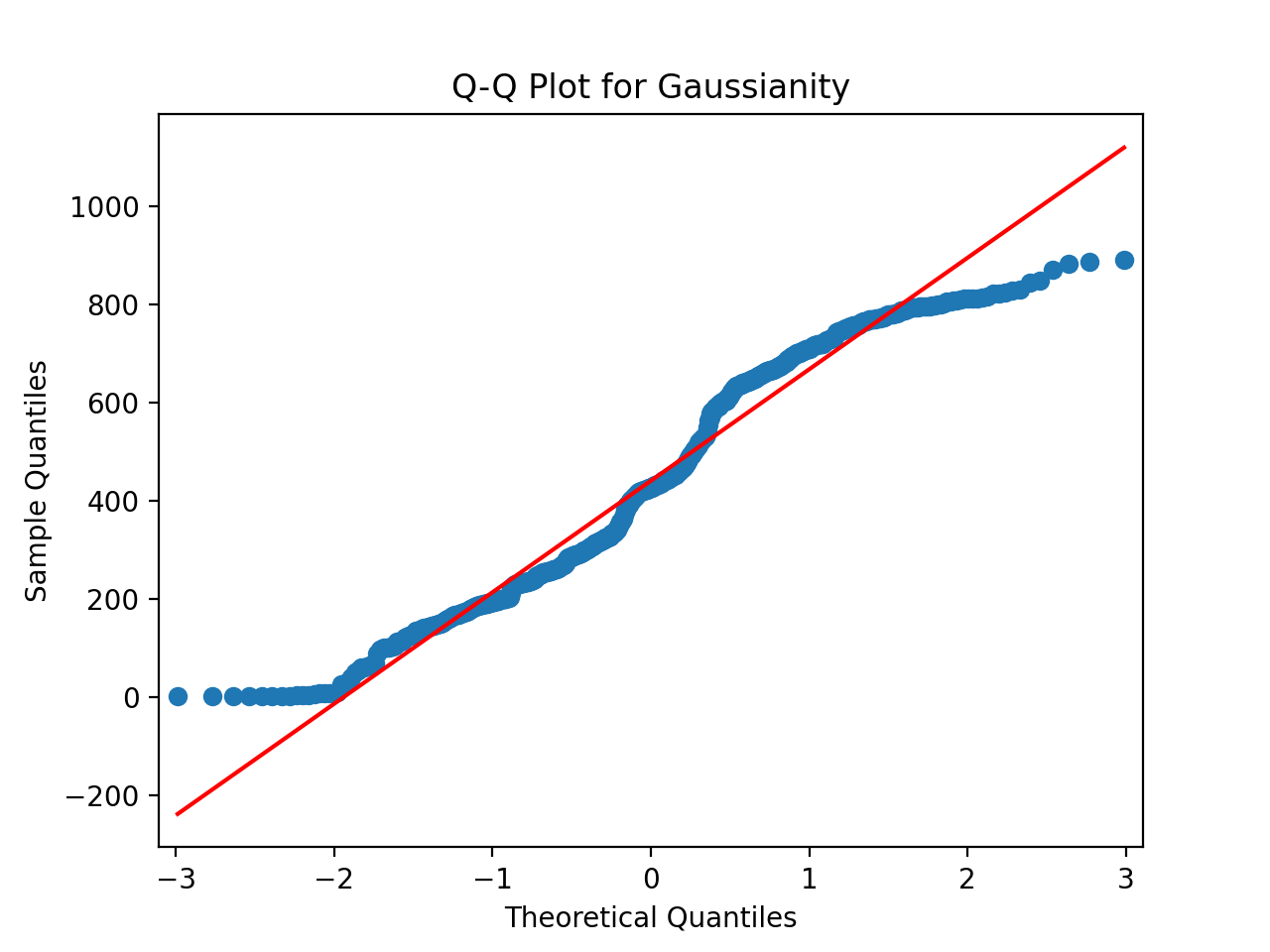}
			\caption{\textbf{Global Q-Q plot}}
		\end{subfigure}
		\hfill
		\begin{subfigure}[t]{0.48\textwidth}
			\centering
			\includegraphics[width=\linewidth]{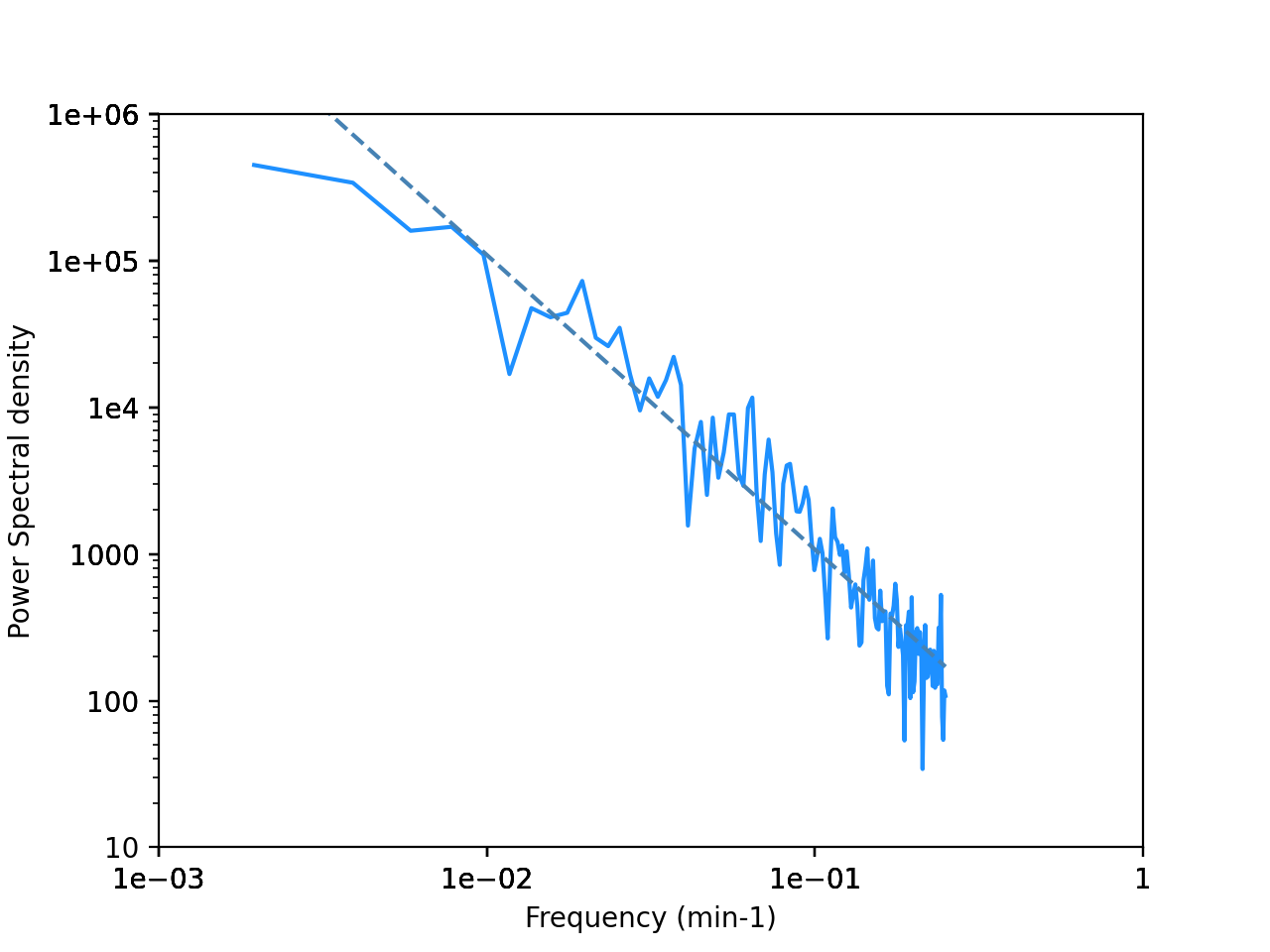}
			\caption{\textbf{Global PSD with regression}}
		\end{subfigure}
		
		\vspace{0.4cm}
		
		\begin{subfigure}[t]{0.60\textwidth}
			\centering
			\includegraphics[width=\linewidth]{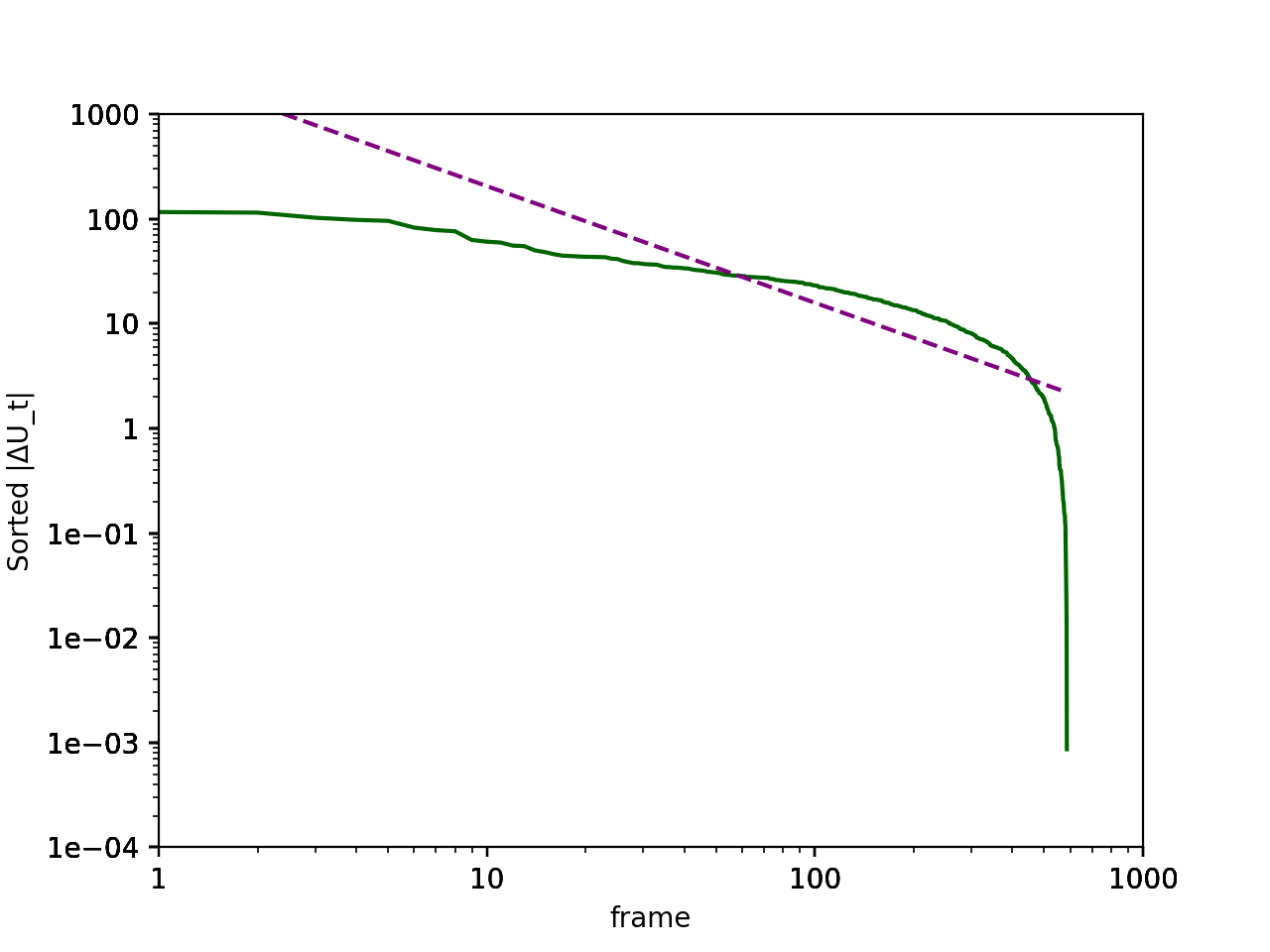}
			\caption{\textbf{Hill-type tail diagnostic.}
				Log--log plot of ordered absolute increments $\Delta U_t$. 
				The absence of a stable scaling plateau and the rapid decay of the tail 
				are incompatible with heavy-tailed behavior and support Gaussian statistics.}
		\end{subfigure}
		
		\vspace{0.4cm}
		
		\begin{subfigure}[t]{0.48\textwidth}
			\centering
			\includegraphics[width=\linewidth]{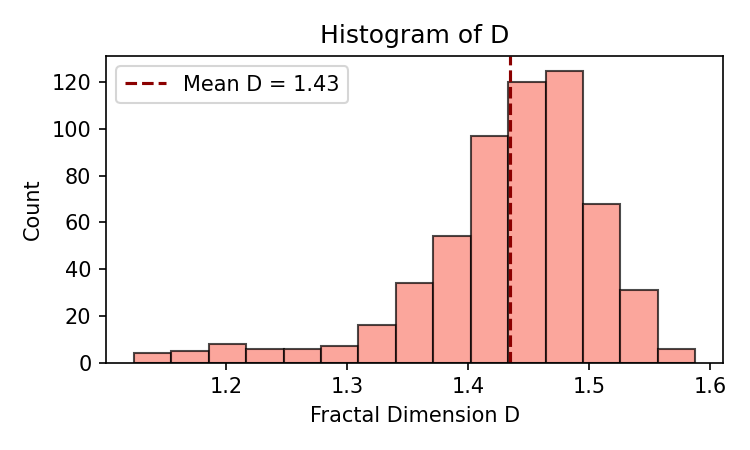}
			\caption{\textbf{Global fractal dimension}}
		\end{subfigure}
		\hfill
		\begin{subfigure}[t]{0.48\textwidth}
			\centering
			\includegraphics[width=\linewidth]{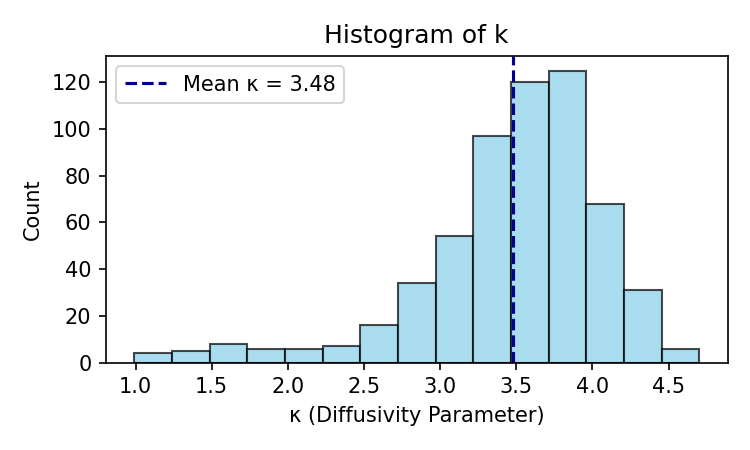}
			\caption{\textbf{Global diffusivity parameter $\kappa$}}
		\end{subfigure}
		
		\caption{\textbf{Global statistical diagnostics for the largest connected component of the network configuration during the retraction phase.}
			Q-Q plot, power spectral density with regression, Hill-type tail analysis, and histograms of the global fractal dimension and diffusivity parameter $\kappa$. 
			Here, \emph{global} refers to the reconstruction performed on the largest connected component of the network configuration during the retraction phase.}
		\label{QQPlotsPSDEtHistogrammesGlobauxSLESlimmingBiggestComponent}
	\end{figure*}
	


	\begin{figure*}[h!]
		\centering
		
		\begin{subfigure}[t]{0.32\textwidth}
			\centering
			\includegraphics[width=\linewidth]{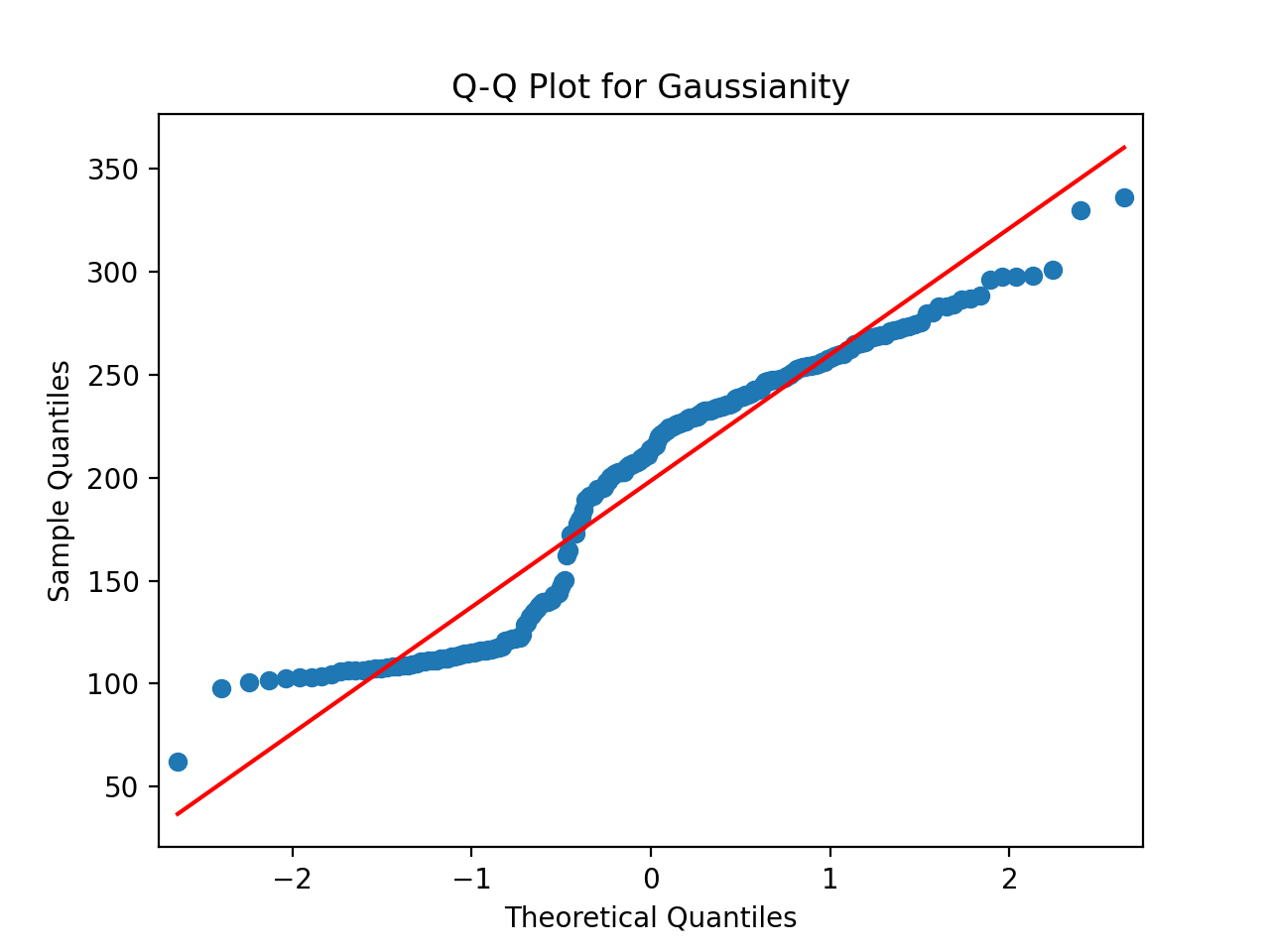}
			\caption{\textbf{Window 2}}
		\end{subfigure}
		\hfill
		\begin{subfigure}[t]{0.32\textwidth}
			\centering
			\includegraphics[width=\linewidth]{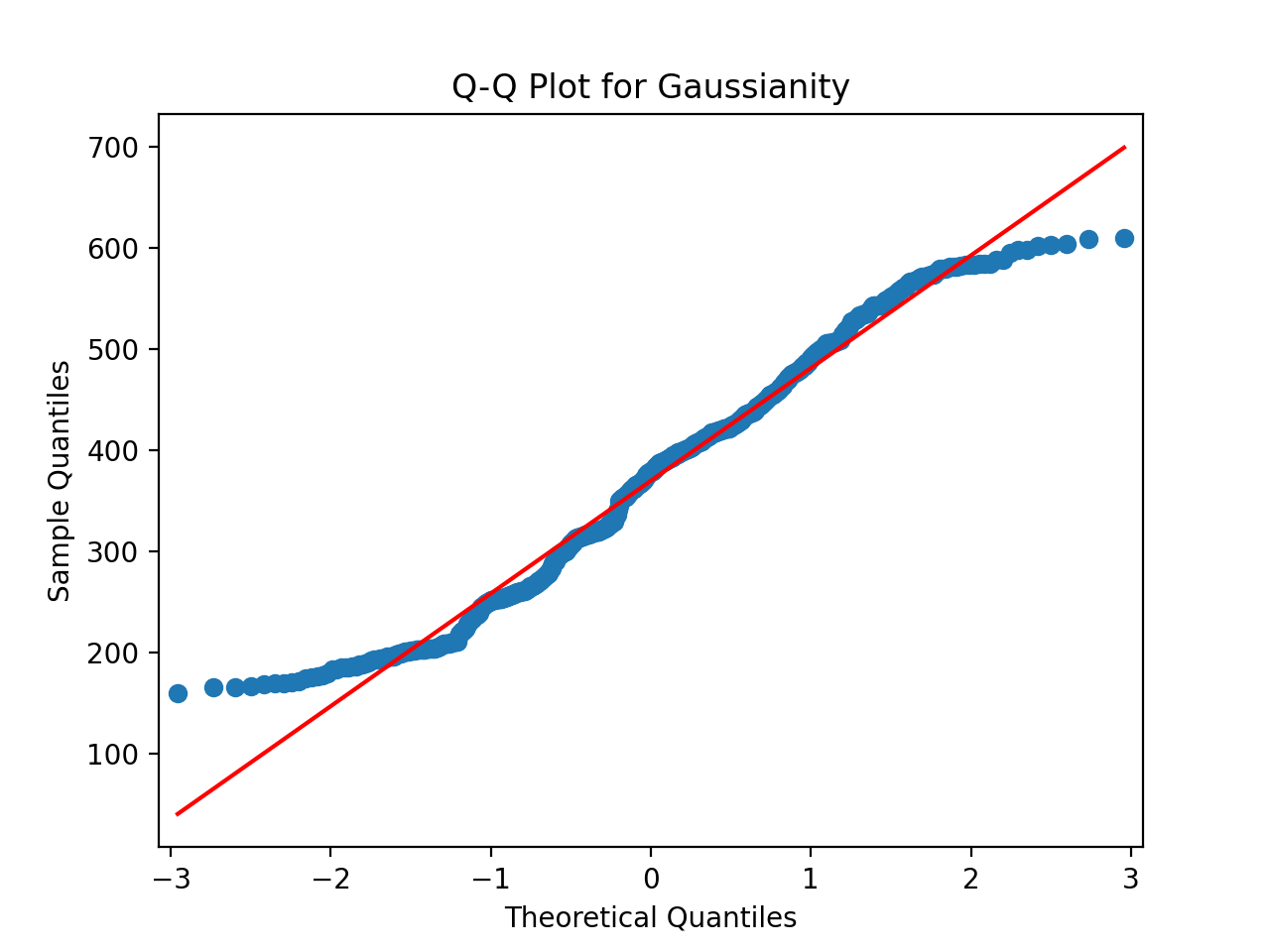}
			\caption{\textbf{Window 3}}
		\end{subfigure}
		\hfill
		\begin{subfigure}[t]{0.32\textwidth}
			\centering
			\includegraphics[width=\linewidth]{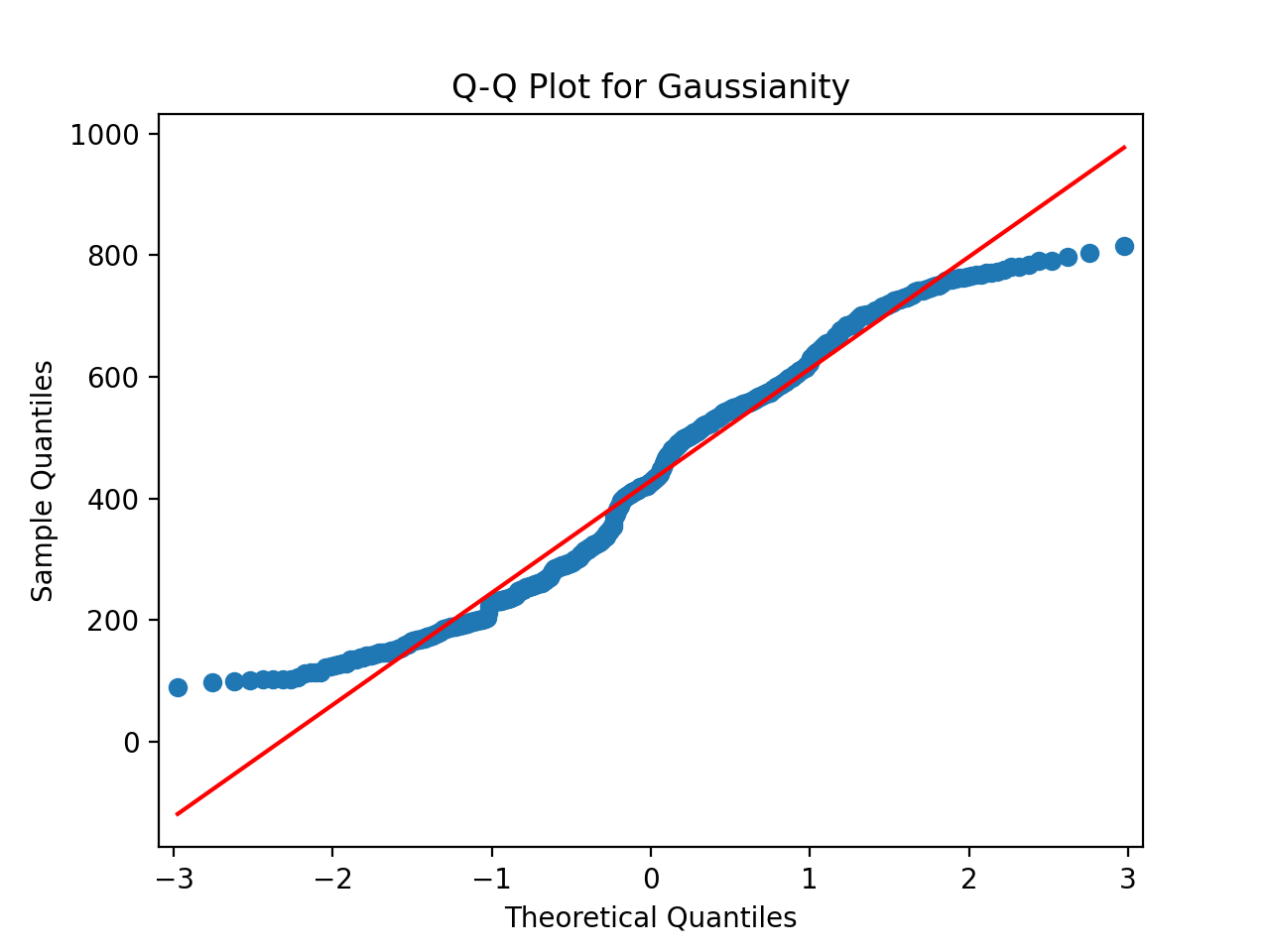}
			\caption{\textbf{Window 4}}
		\end{subfigure}
		
		\vspace{0.35cm}
		
		\begin{subfigure}[t]{0.32\textwidth}
			\centering
			\includegraphics[width=\linewidth]{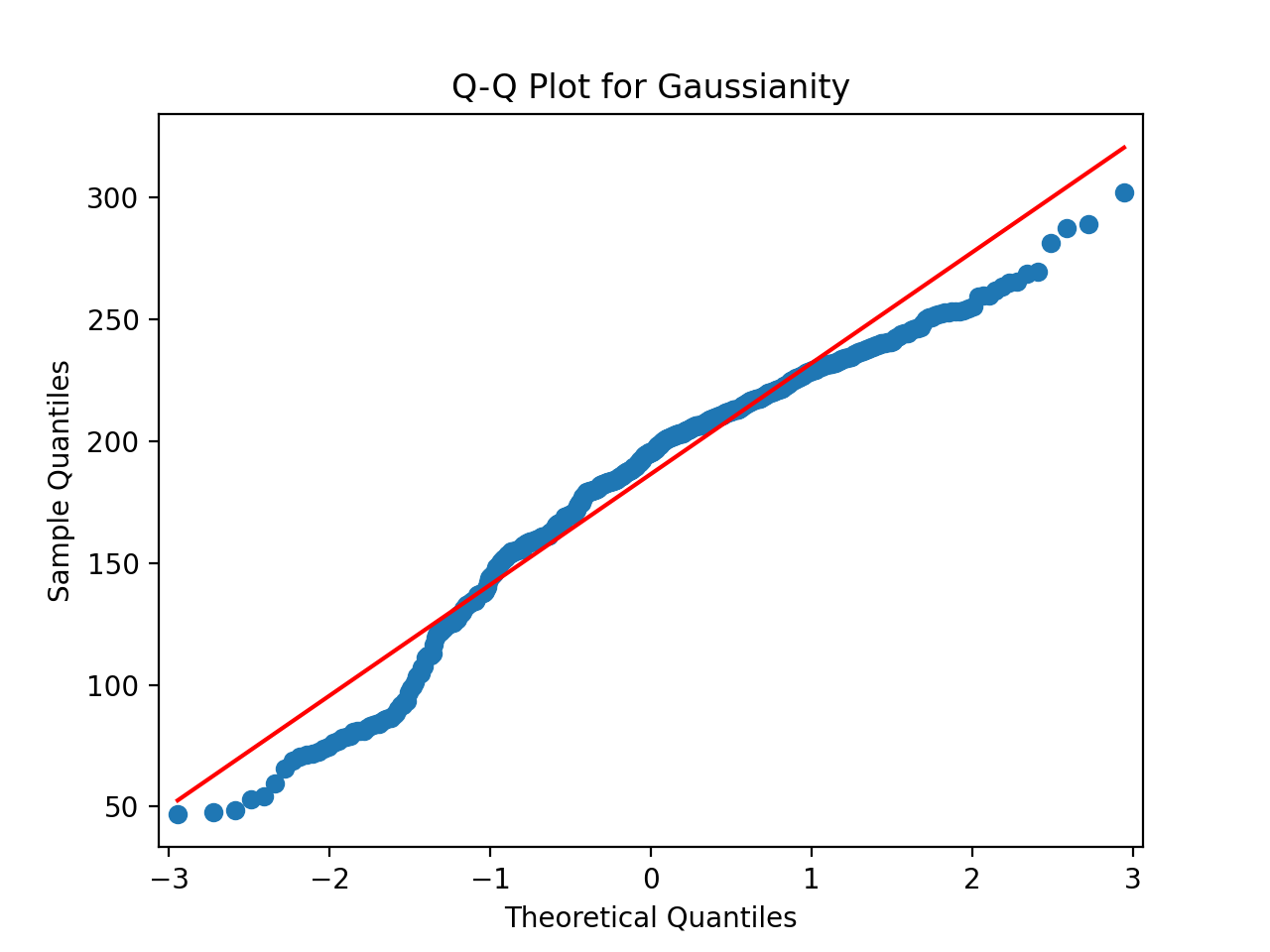}
			\caption{\textbf{Window 7}}
		\end{subfigure}
		\hfill
		\begin{subfigure}[t]{0.32\textwidth}
			\centering
			\includegraphics[width=\linewidth]{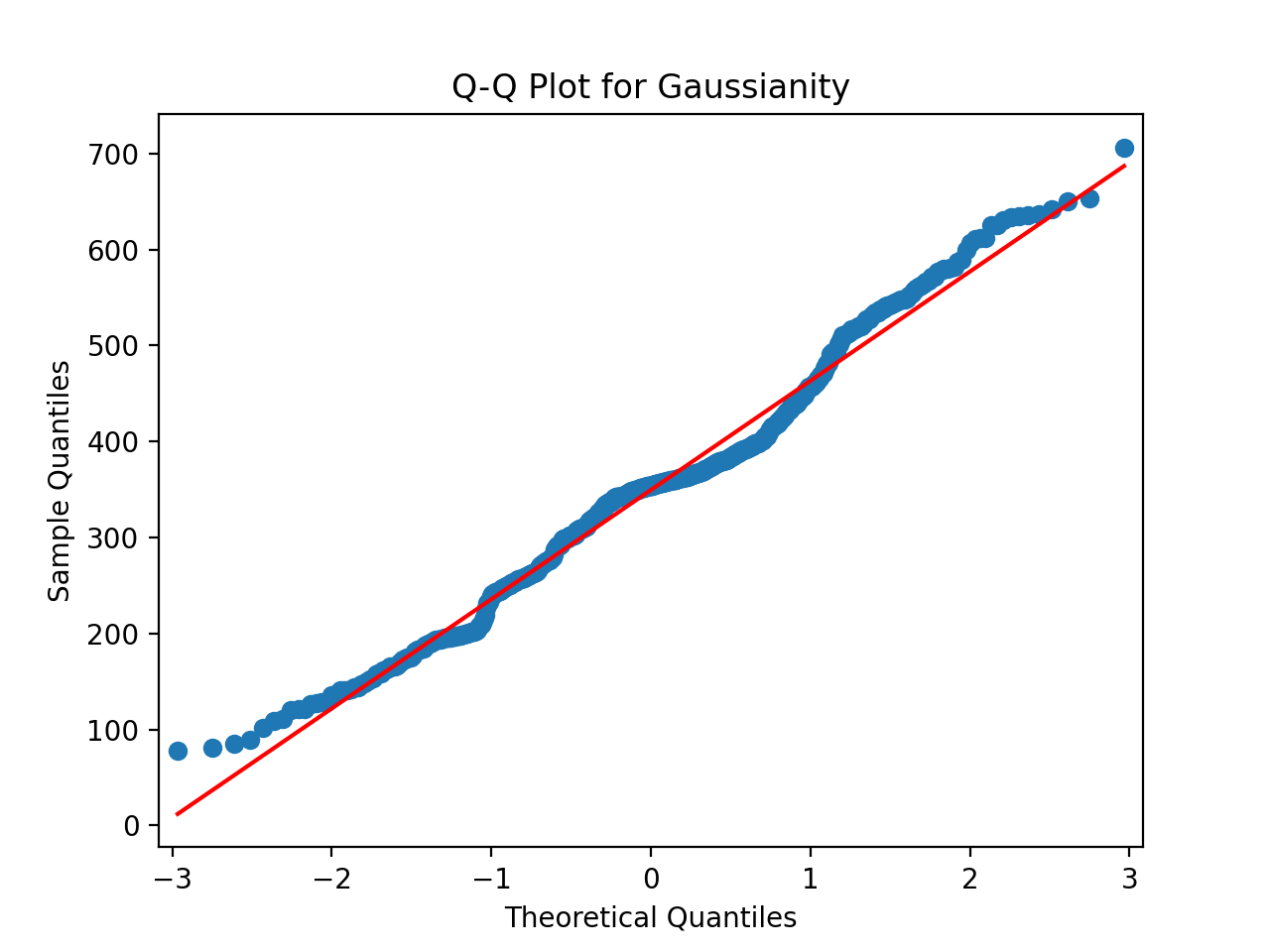}
			\caption{\textbf{Window 8}}
		\end{subfigure}
		\hfill
		\begin{subfigure}[t]{0.32\textwidth}
			\centering
			\includegraphics[width=\linewidth]{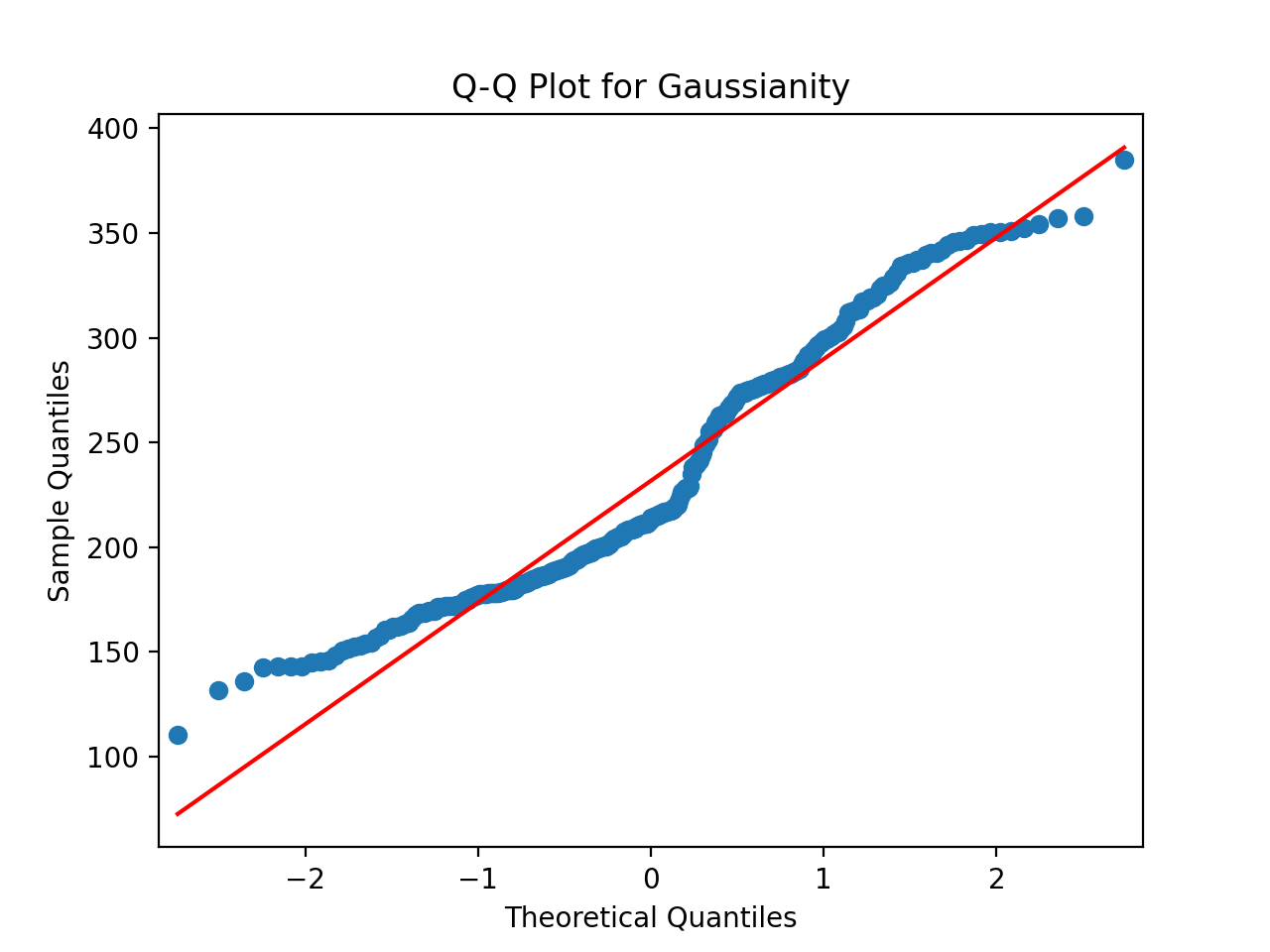}
			\caption{\textbf{Window 12}}
		\end{subfigure}
		
		\caption{\textbf{Local Q-Q plots of the driving function increments $\Delta U_t$.}
			Representative outer spatial windows (2, 3, 4, 7, 8, 12) computed on the largest connected component of the network configuration during the retraction phase. 
			As for the pseudopod and network components, the empirical quantiles closely follow the theoretical Gaussian quantiles, with only mild deviations at extreme tails.}
		\label{QQPlotsLocauxSLESlimmingBiggestComponent}
	\end{figure*}
	

	\begin{figure*}[h!]
		\centering
		
		\begin{subfigure}[t]{0.32\textwidth}
			\centering
			\includegraphics[width=\linewidth]{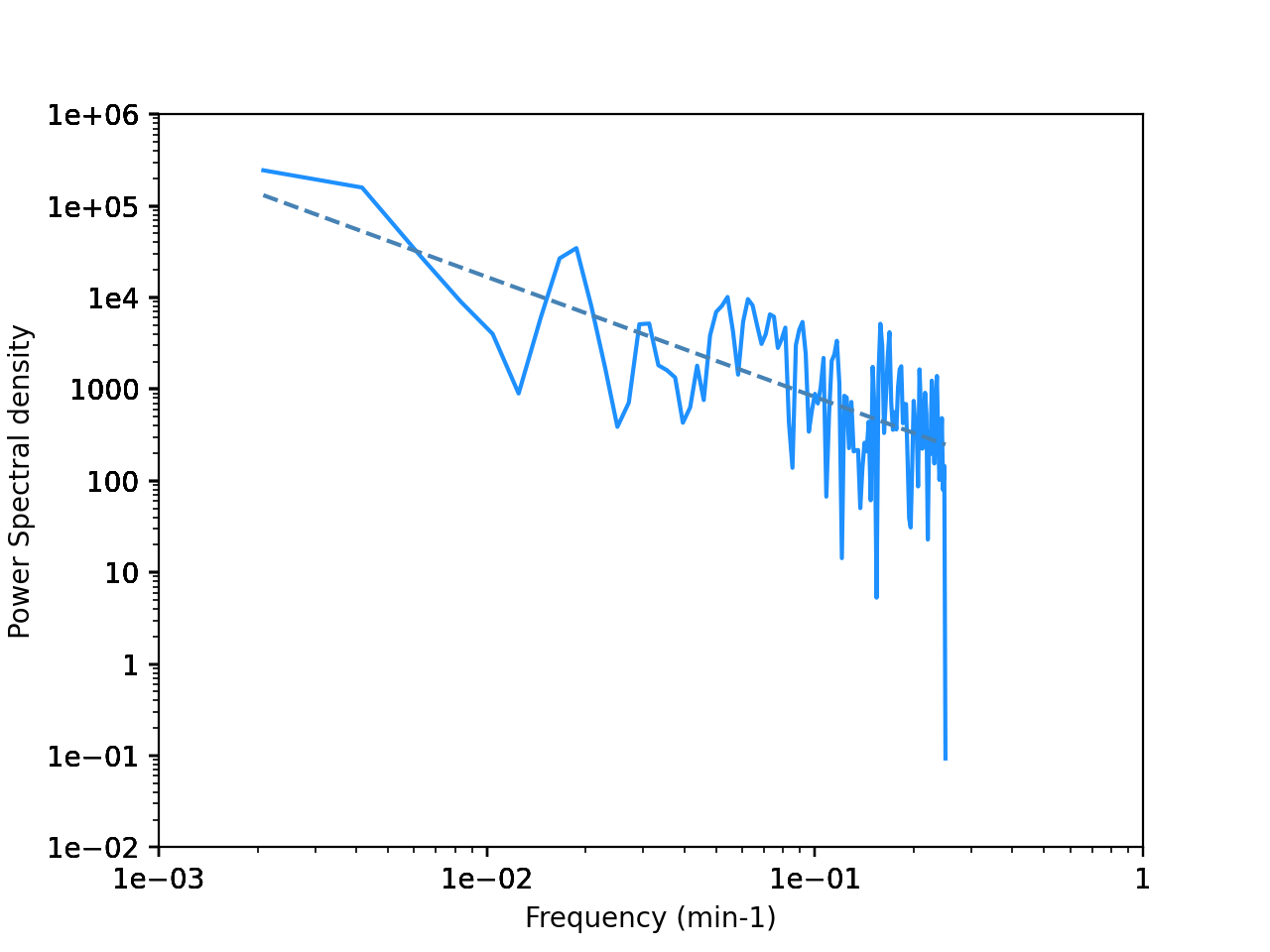}
			\caption{\textbf{Window 2}}
		\end{subfigure}
		\hfill
		\begin{subfigure}[t]{0.32\textwidth}
			\centering
			\includegraphics[width=\linewidth]{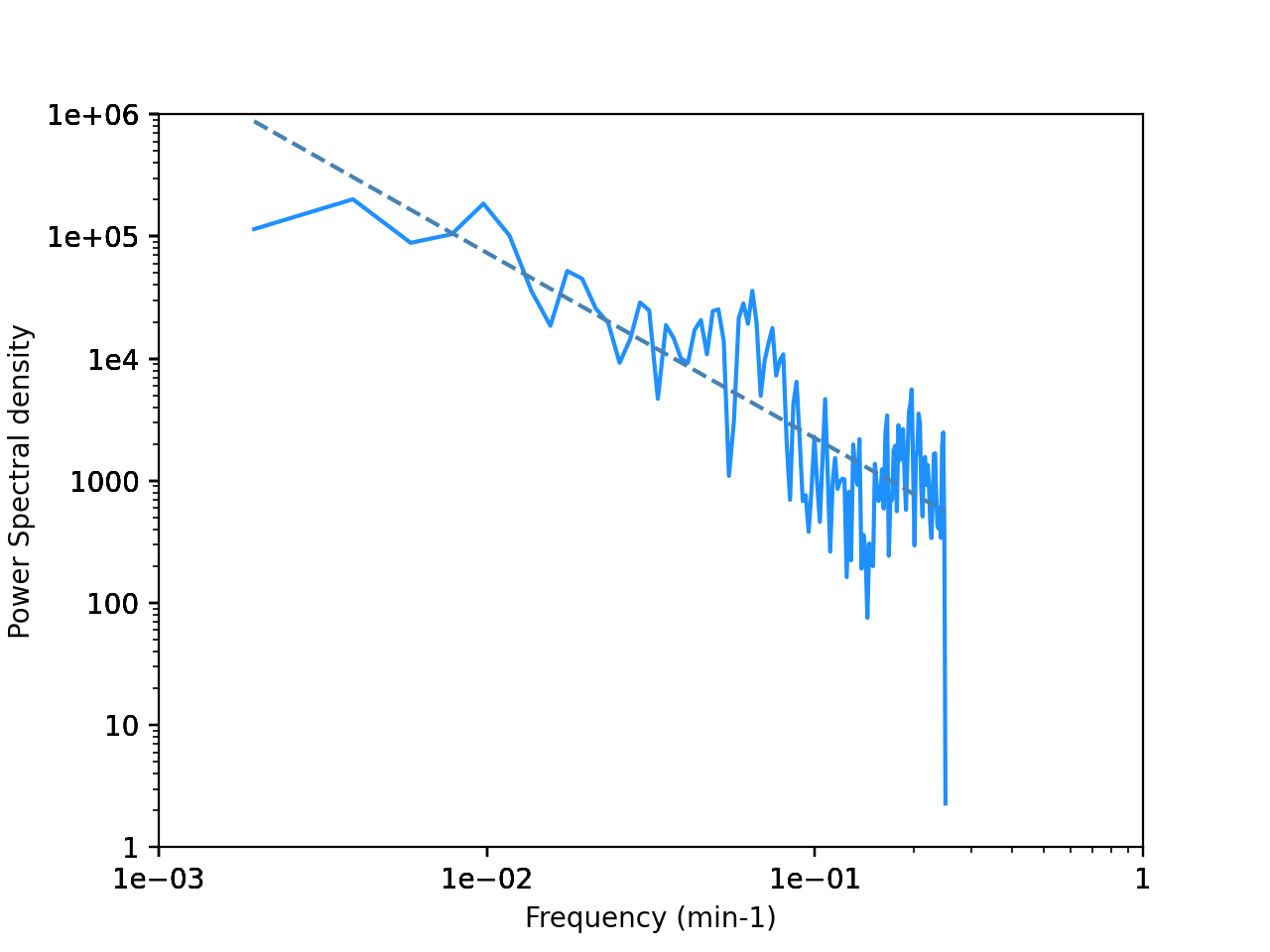}
			\caption{\textbf{Window 3}}
		\end{subfigure}
		\hfill
		\begin{subfigure}[t]{0.32\textwidth}
			\centering
			\includegraphics[width=\linewidth]{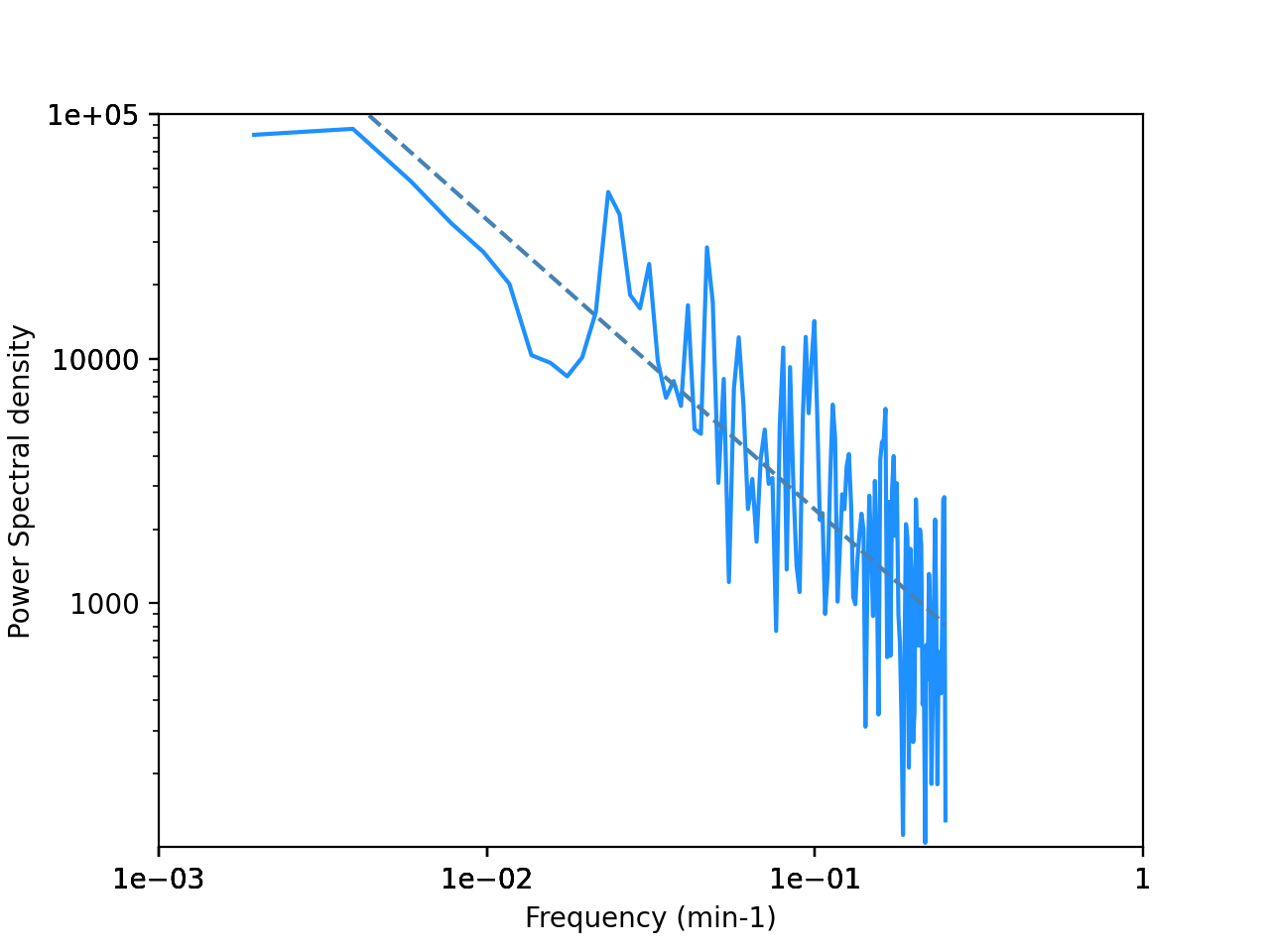}
			\caption{\textbf{Window 4}}
		\end{subfigure}
		
		\vspace{0.35cm}
		
		\begin{subfigure}[t]{0.32\textwidth}
			\centering
			\includegraphics[width=\linewidth]{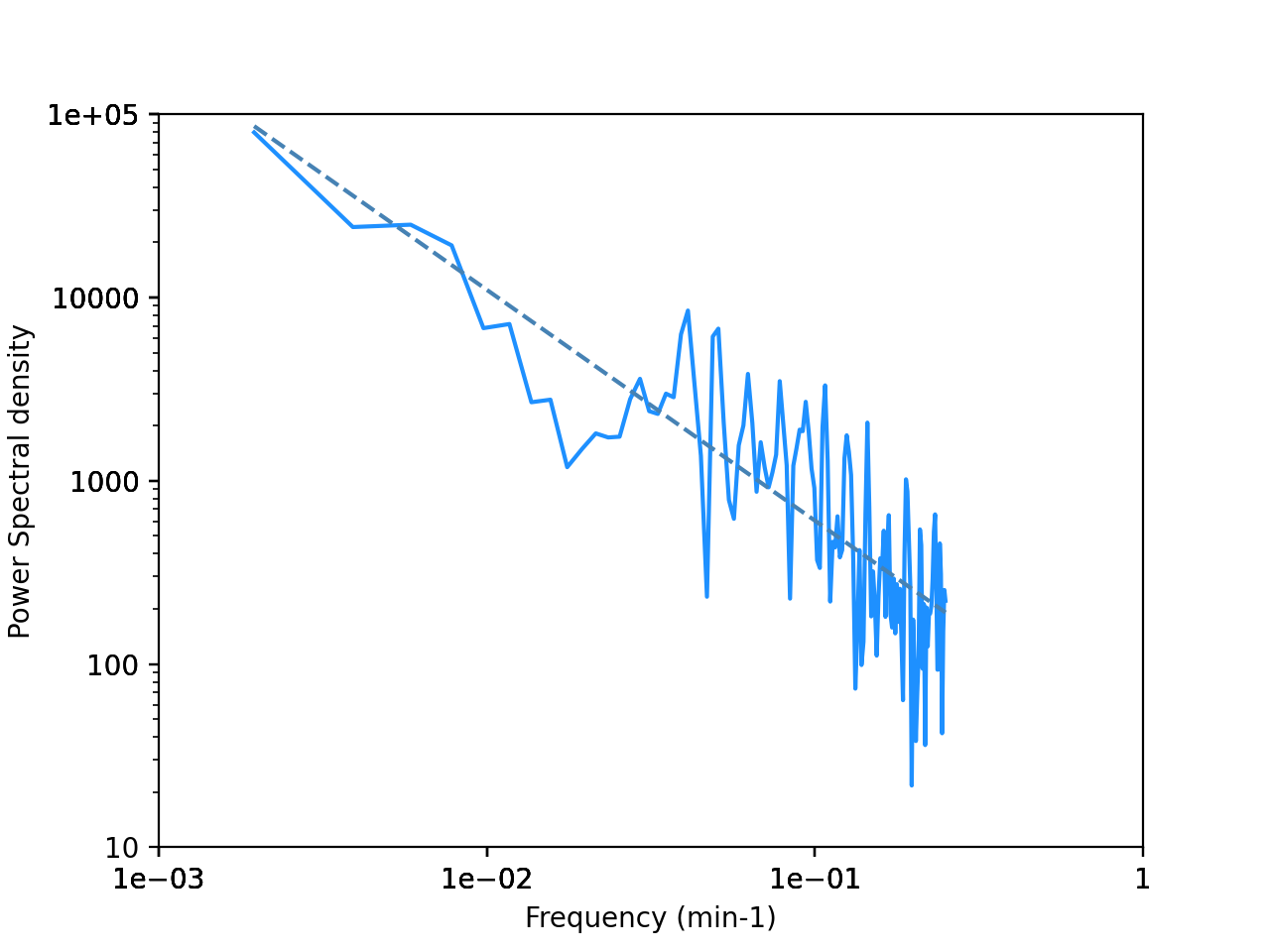}
			\caption{\textbf{Window 7}}
		\end{subfigure}
		\hfill
		\begin{subfigure}[t]{0.32\textwidth}
			\centering
			\includegraphics[width=\linewidth]{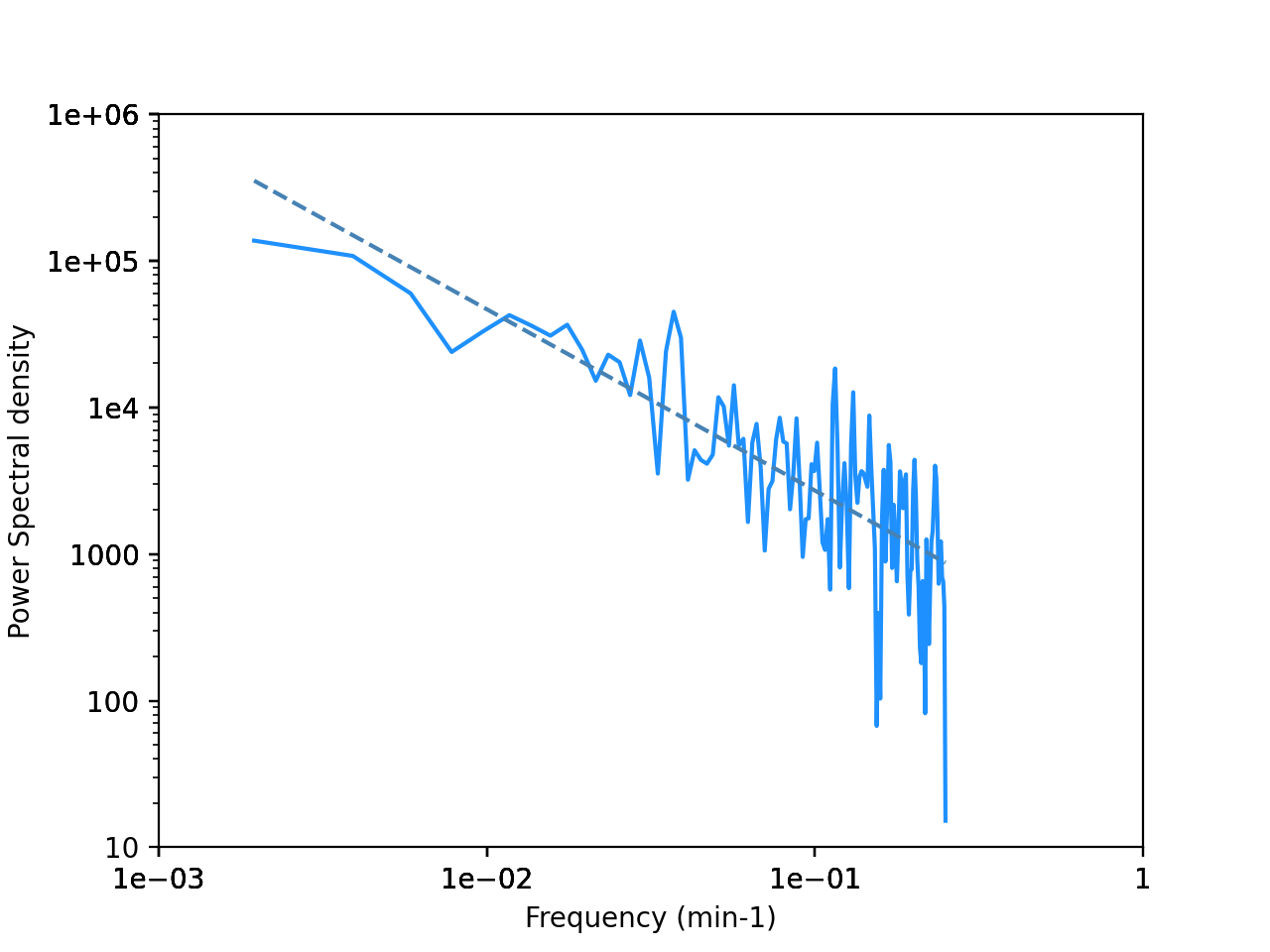}
			\caption{\textbf{Window 8}}
		\end{subfigure}
		\hfill
		\begin{subfigure}[t]{0.32\textwidth}
			\centering
			\includegraphics[width=\linewidth]{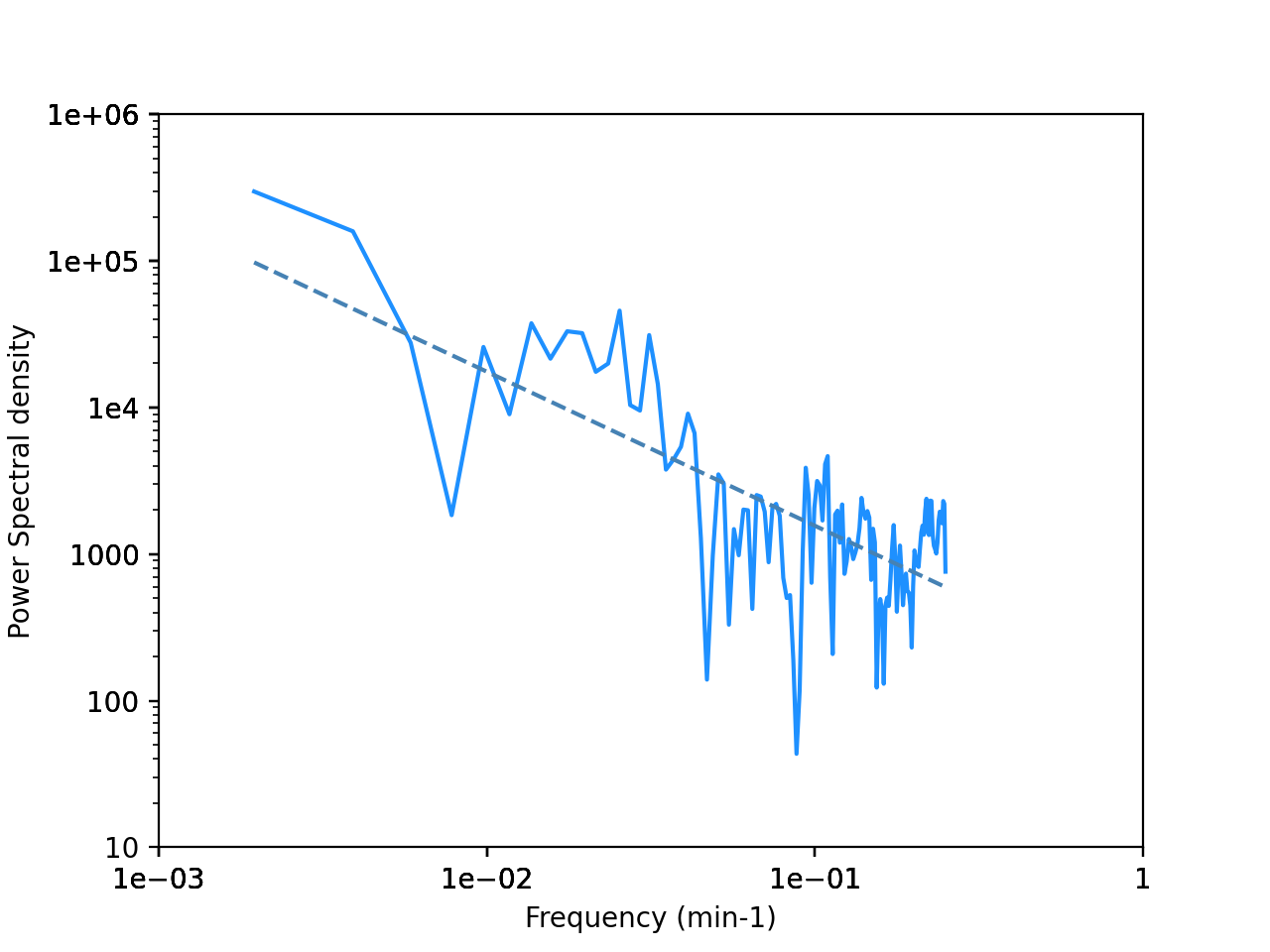}
			\caption{\textbf{Window 12}}
		\end{subfigure}
		
		\caption{\textbf{Local power spectral density (PSD) diagnostics for the network configuration during the retraction phase component.}
			Log--log PSD plots with linear regression for representative outer windows (2, 3, 4, 7, 8, 12), computed on the largest connected component of the network configuration during the retraction phase. 
			In each window, the fitted slope provides an estimate of the local scaling exponent $\beta$ in $S(\omega)\propto\omega^{-\beta}$. 
			The slopes remain broadly compatible with a power-law regime and exhibit reduced inter-window variability compared with the pseudopod and network levels, suggesting that the backbone-dominated structure retains coherent Brownian-type dynamics.}
		\label{PSDWithRegressionLocauxSLESlimmingBiggestComponent}
	\end{figure*}
	
	

	\begin{figure*}[h!]
		\centering
		
		\begin{subfigure}[t]{0.32\textwidth}
			\centering
			\includegraphics[width=\linewidth]{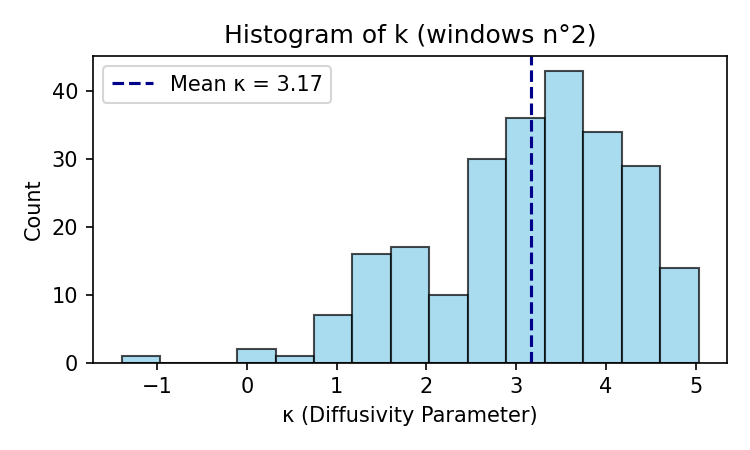}
			\caption{\textbf{Window 2}}
		\end{subfigure}
		\hfill
		\begin{subfigure}[t]{0.32\textwidth}
			\centering
			\includegraphics[width=\linewidth]{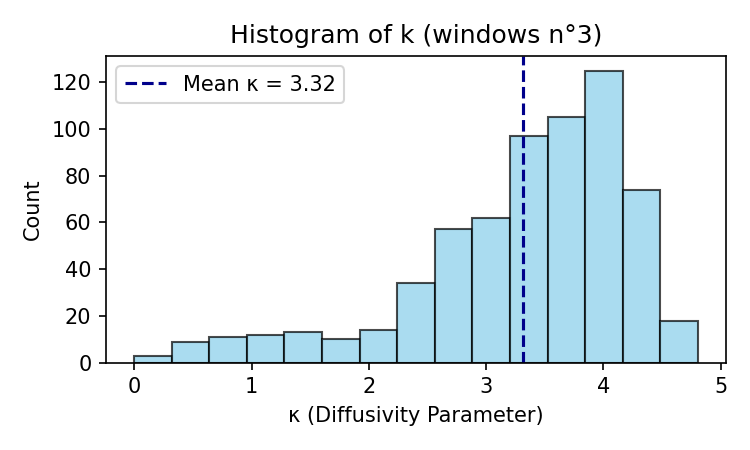}
			\caption{\textbf{Window 3}}
		\end{subfigure}
		\hfill
		\begin{subfigure}[t]{0.32\textwidth}
			\centering
			\includegraphics[width=\linewidth]{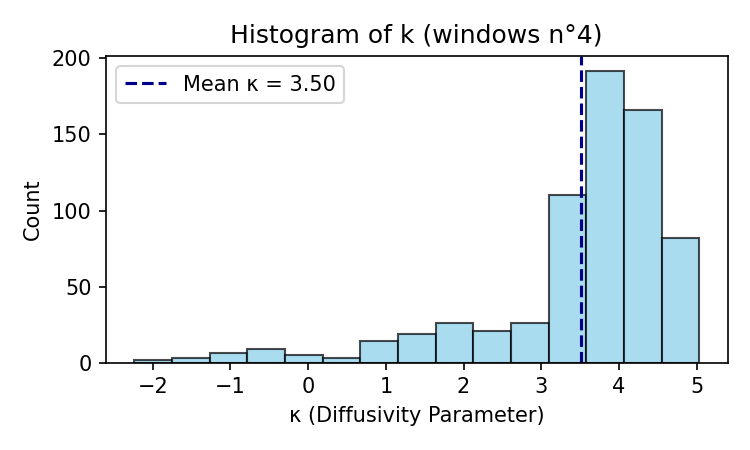}
			\caption{\textbf{Window 4}}
		\end{subfigure}
		
		\vspace{0.35cm}
		
		\begin{subfigure}[t]{0.32\textwidth}
			\centering
			\includegraphics[width=\linewidth]{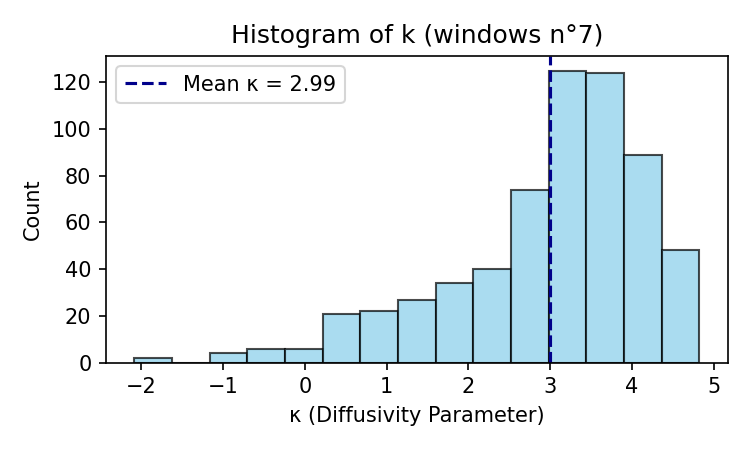}
			\caption{\textbf{Window 7}}
		\end{subfigure}
		\hfill
		\begin{subfigure}[t]{0.32\textwidth}
			\centering
			\includegraphics[width=\linewidth]{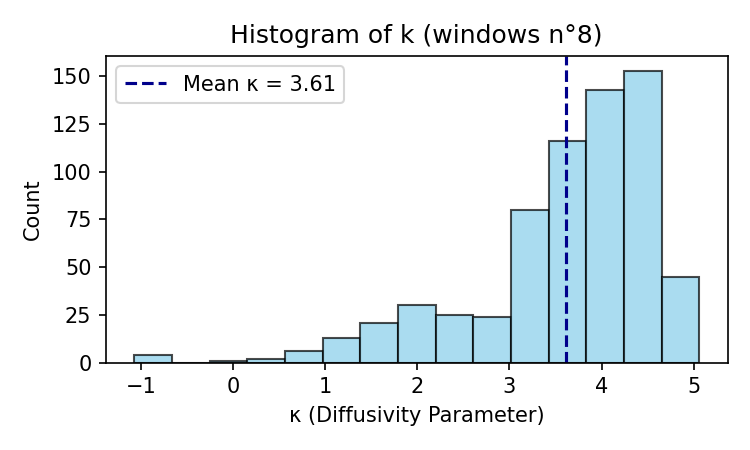}
			\caption{\textbf{Window 8}}
		\end{subfigure}
		\hfill
		\begin{subfigure}[t]{0.32\textwidth}
			\centering
			\includegraphics[width=\linewidth]{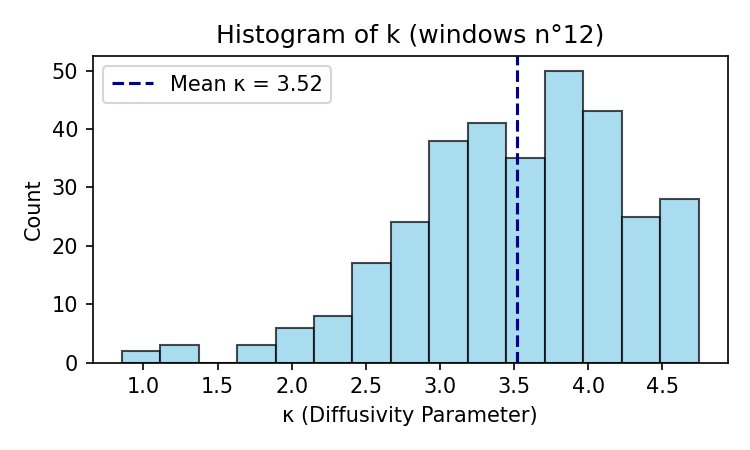}
			\caption{\textbf{Window 12}}
		\end{subfigure}
		
		\caption{\textbf{Local histograms of the diffusivity parameter $\kappa$ (network configuration during the expansion phase).}
			Estimated values of $\kappa$ for representative outer windows (2, 3, 4, 7, 8, 12), computed on the largest connected component of the network configuration during the retraction phase. 
			Because of the lower fractal dimension and the strongly filamentary structure of the network configuration during the retraction phase, local estimates of $\kappa$ may exhibit increased variability, including small or occasionally negative values. 
			These fluctuations reflect finite-size effects and limitations of the Brownian-type scaling approximation at these scales, rather than physically meaningful negative diffusivities.}
		\label{HistogrammesLocauxKappaSLESlimmingBiggestComponent}
	\end{figure*}
	


	\begin{figure*}[h!]
		\centering
		
		\begin{subfigure}[t]{0.48\textwidth}
			\centering
			\includegraphics[width=\linewidth]{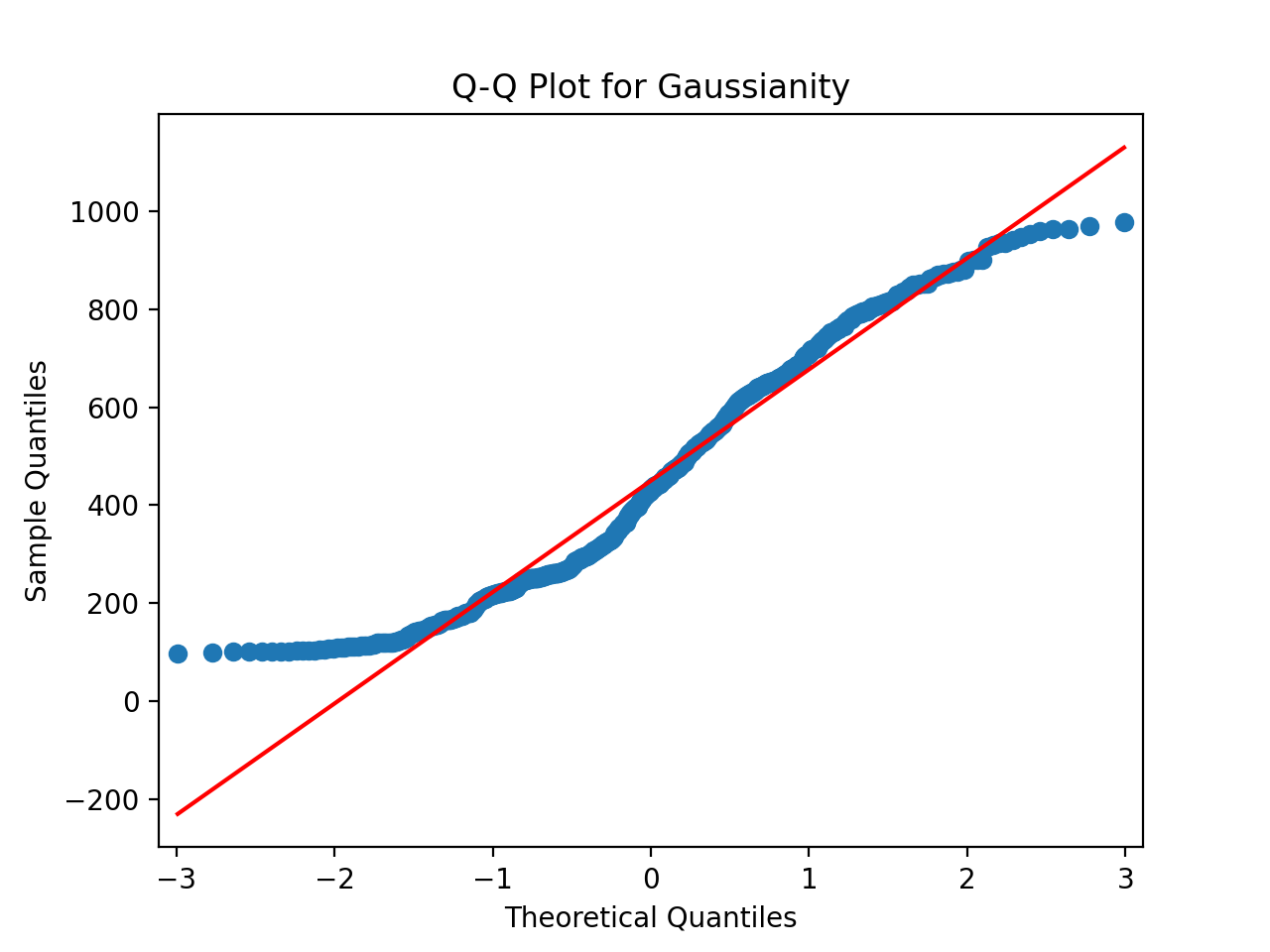}
			\caption{\textbf{Global Q-Q plot}}
		\end{subfigure}
		\hfill
		\begin{subfigure}[t]{0.48\textwidth}
			\centering
			\includegraphics[width=\linewidth]{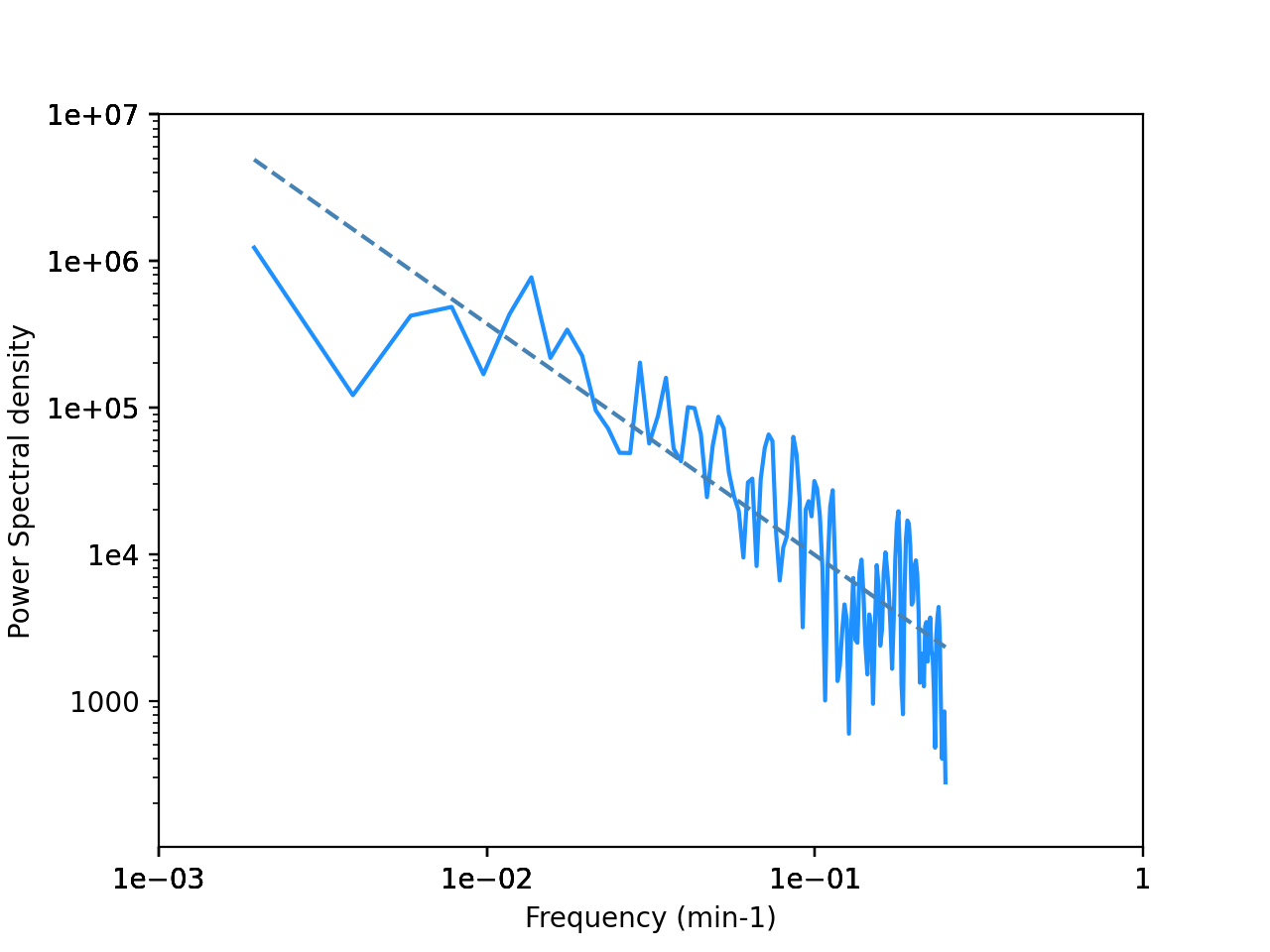}
			\caption{\textbf{Global PSD with regression}}
		\end{subfigure}
		
		\vspace{0.4cm}
		
		\begin{subfigure}[t]{0.60\textwidth}
			\centering
			\includegraphics[width=\linewidth]{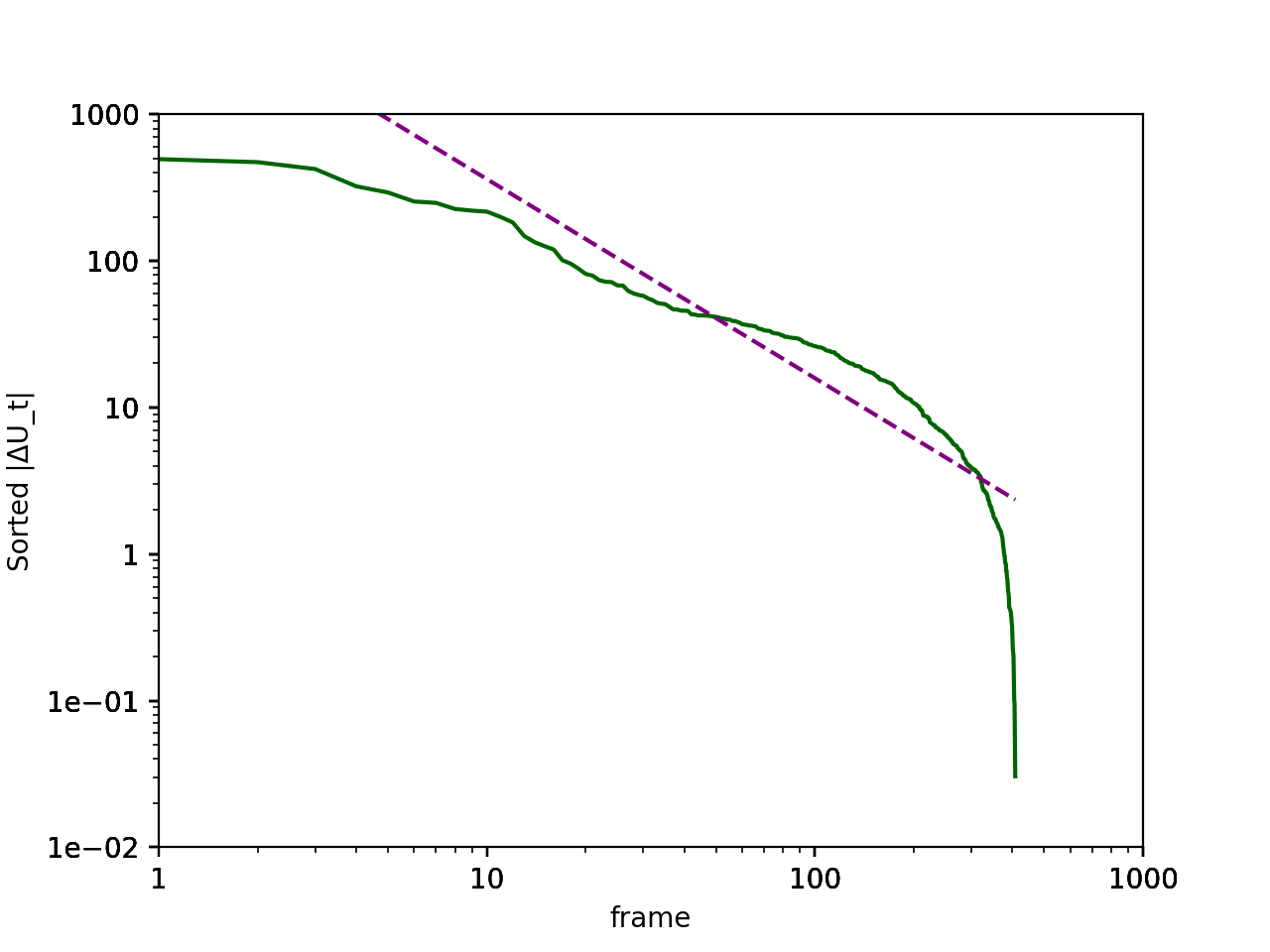}
			\caption{\textbf{Hill-type tail diagnostic.}
				Log--log plot of ordered absolute increments $\Delta U_t$. 
				The absence of a stable scaling plateau and the rapid decay of the tail 
				are incompatible with heavy-tailed behavior and support Gaussian statistics.}
		\end{subfigure}
		
		\vspace{0.4cm}
		
		\begin{subfigure}[t]{0.48\textwidth}
			\centering
			\includegraphics[width=\linewidth]{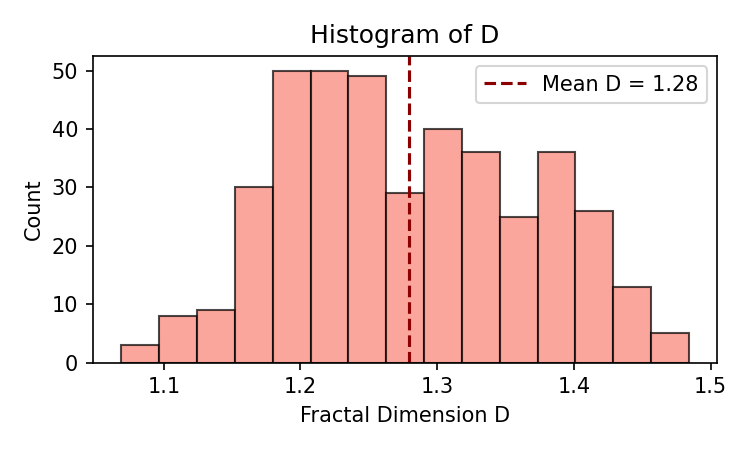}
			\caption{\textbf{Global fractal dimension}}
		\end{subfigure}
		\hfill
		\begin{subfigure}[t]{0.48\textwidth}
			\centering
			\includegraphics[width=\linewidth]{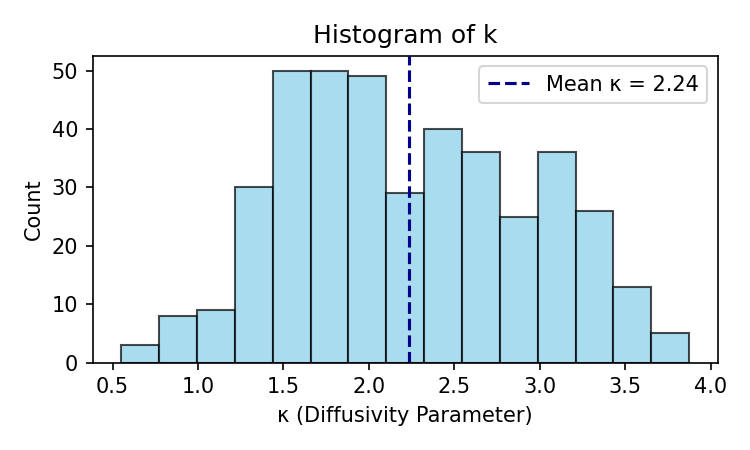}
			\caption{\textbf{Global diffusivity parameter $\kappa$}}
		\end{subfigure}
		
		\caption{\textbf{Global statistical diagnostics for the largest connected component of the network configuration during the expansion phase.}
			Q-Q plot, power spectral density with regression, Hill-type tail analysis, and histograms of the global fractal dimension and diffusivity parameter $\kappa$. 
			Here, \emph{global} refers to the reconstruction performed on the largest connected component of the network configuration during the expansion phase.}
		\label{QQPlotsPSDEtHistogrammesGlobauxSLEThickeningBiggestComponent}
	\end{figure*}


	\begin{figure*}[h!]
		\centering
		
		\begin{subfigure}[t]{0.32\textwidth}
			\centering
			\includegraphics[width=\linewidth]{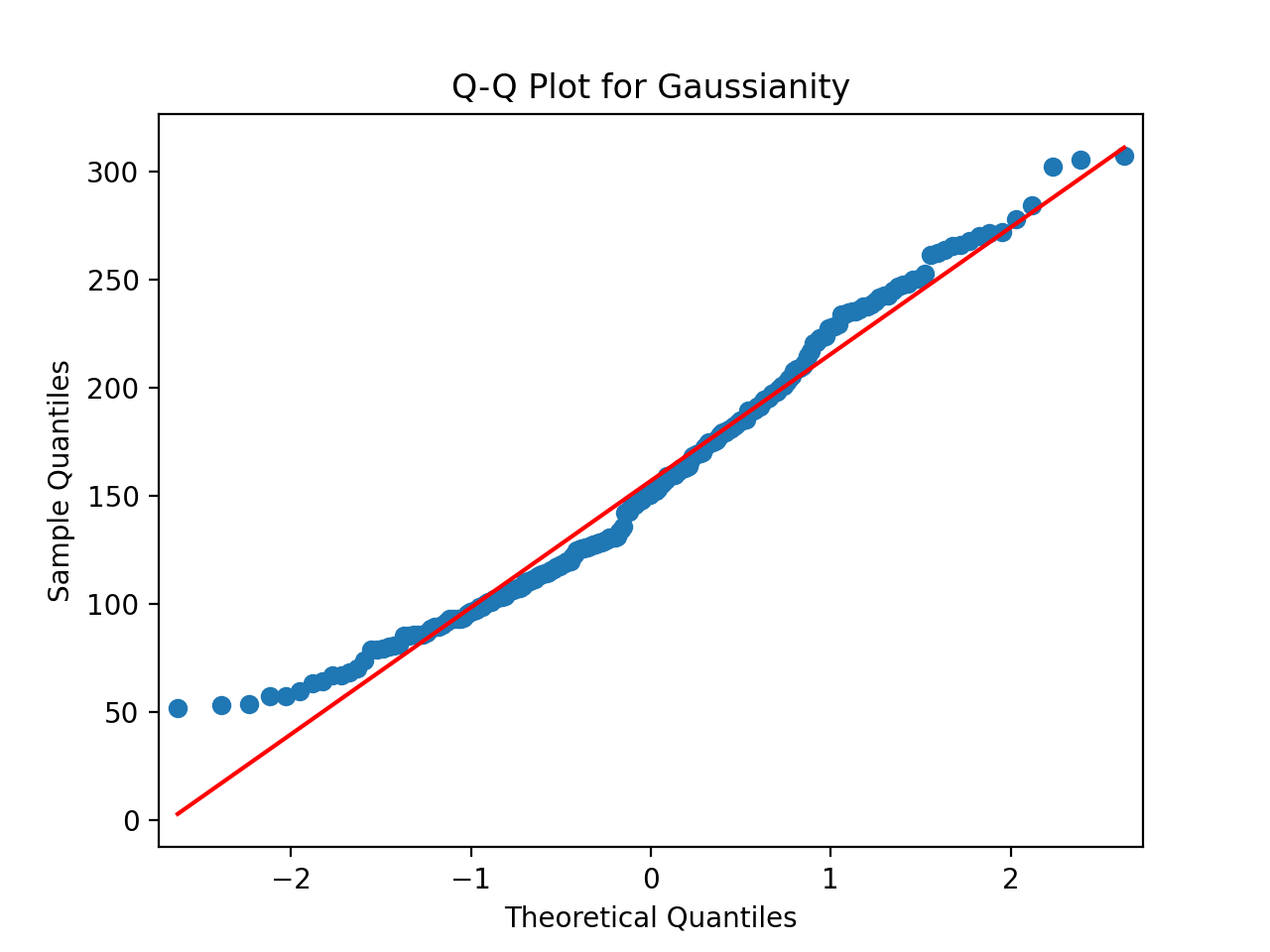}
			\caption{\textbf{Window 2}}
		\end{subfigure}
		\hfill
		\begin{subfigure}[t]{0.32\textwidth}
			\centering
			\includegraphics[width=\linewidth]{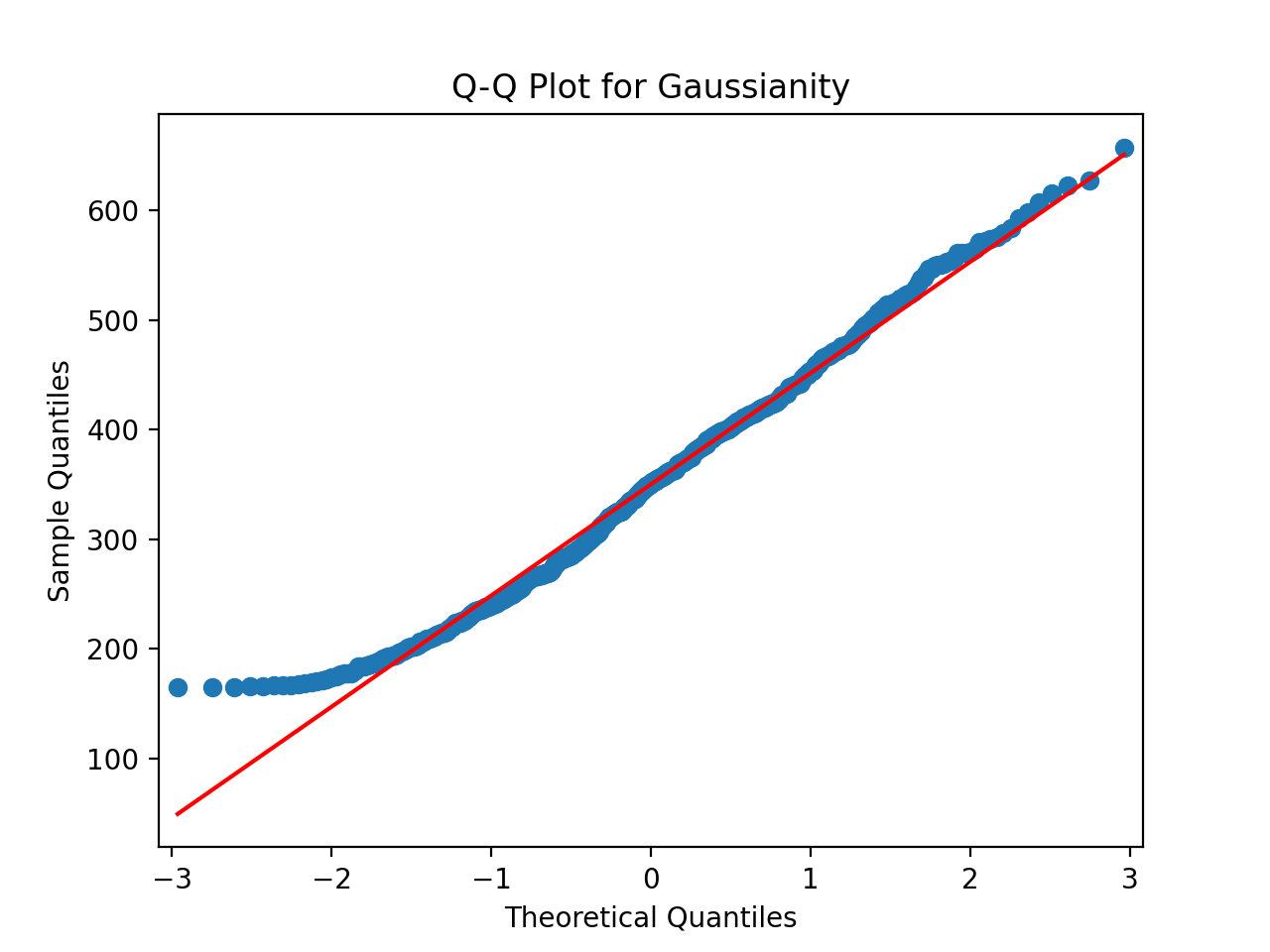}
			\caption{\textbf{Window 3}}
		\end{subfigure}
		\hfill
		\begin{subfigure}[t]{0.32\textwidth}
			\centering
			\includegraphics[width=\linewidth]{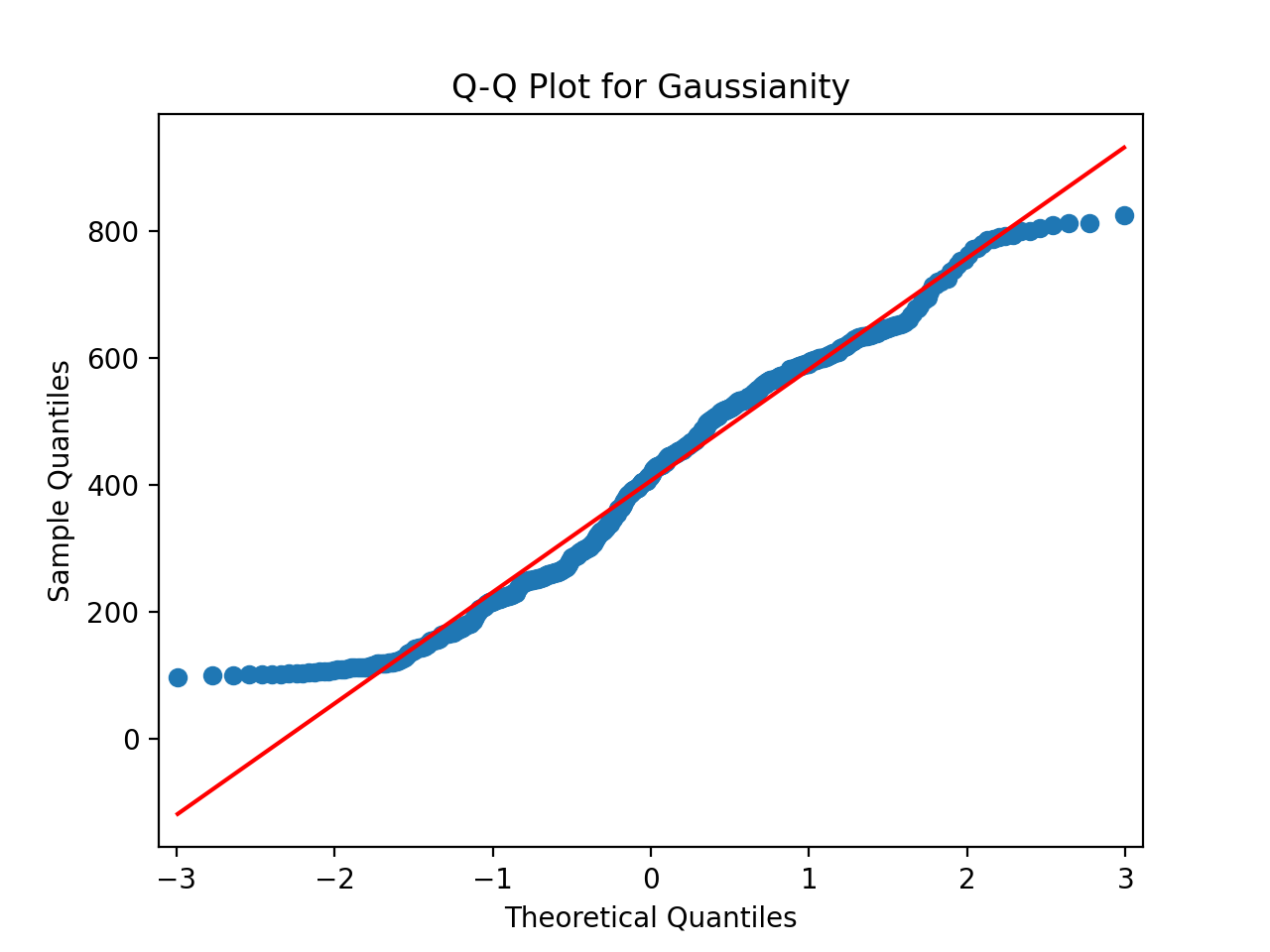}
			\caption{\textbf{Window 4}}
		\end{subfigure}
		
		\vspace{0.35cm}
		
		\begin{subfigure}[t]{0.32\textwidth}
			\centering
			\includegraphics[width=\linewidth]{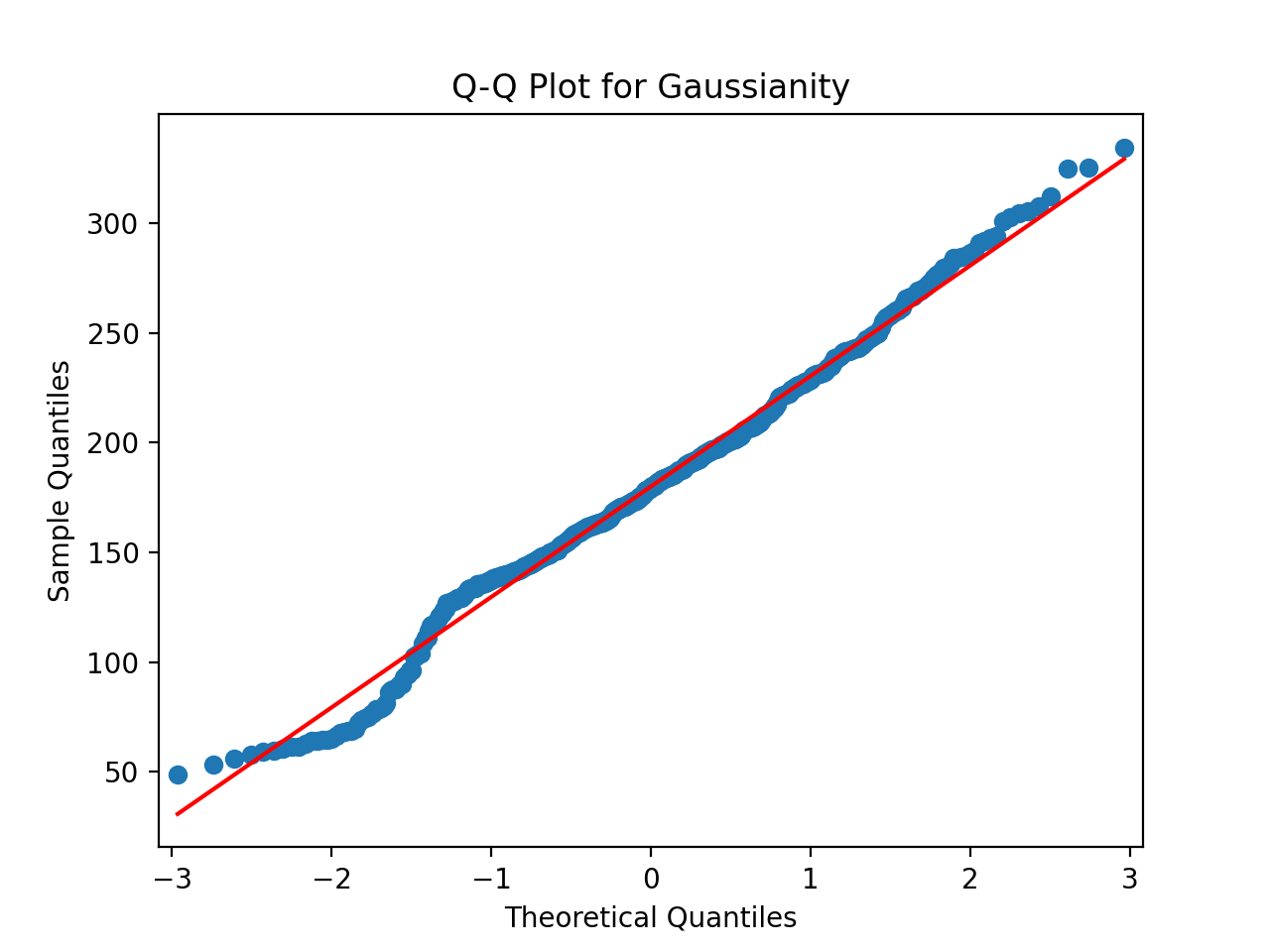}
			\caption{\textbf{Window 7}}
		\end{subfigure}
		\hfill
		\begin{subfigure}[t]{0.32\textwidth}
			\centering
			\includegraphics[width=\linewidth]{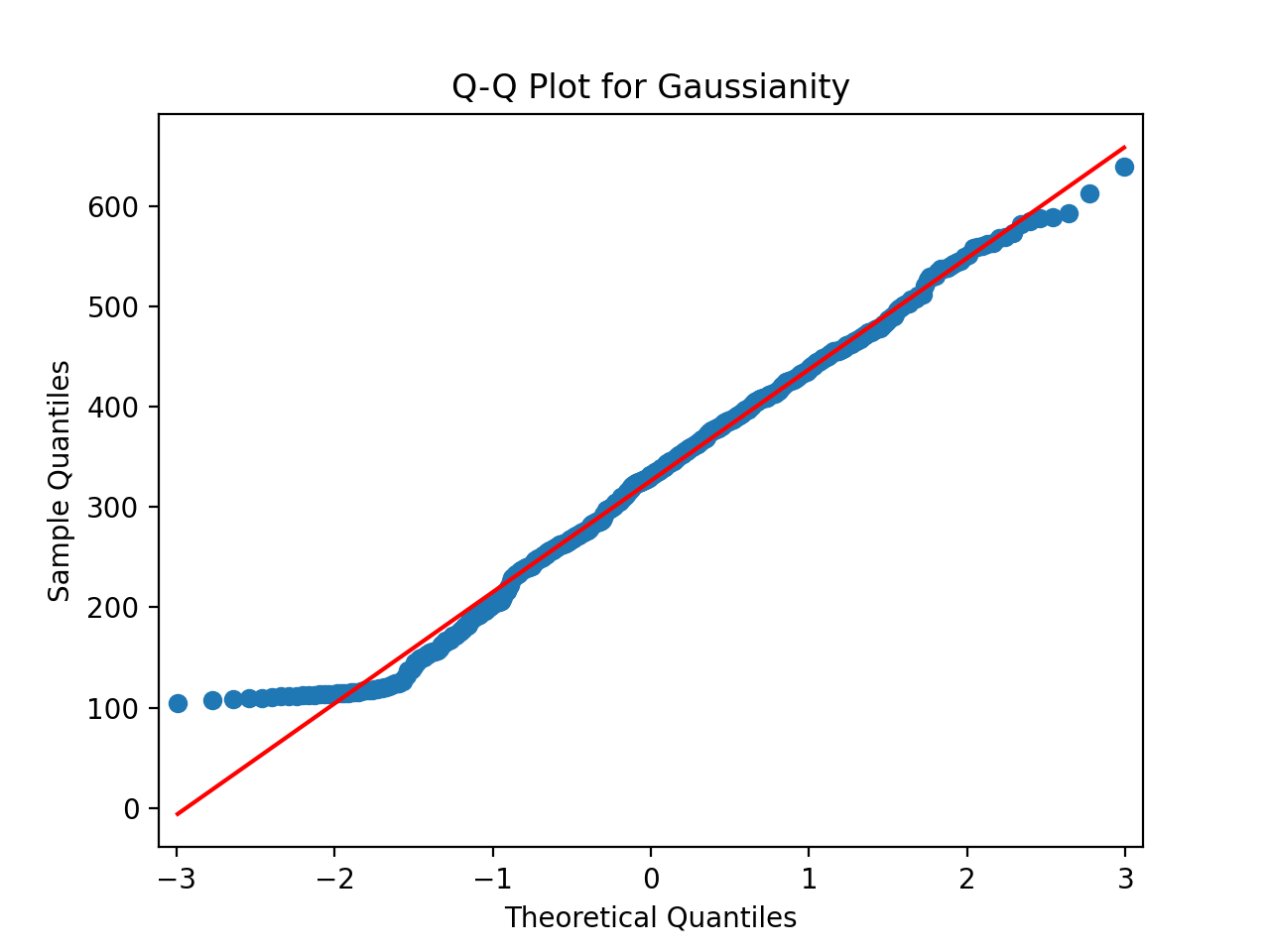}
			\caption{\textbf{Window 8}}
		\end{subfigure}
		\hfill
		\begin{subfigure}[t]{0.32\textwidth}
			\centering
			\includegraphics[width=\linewidth]{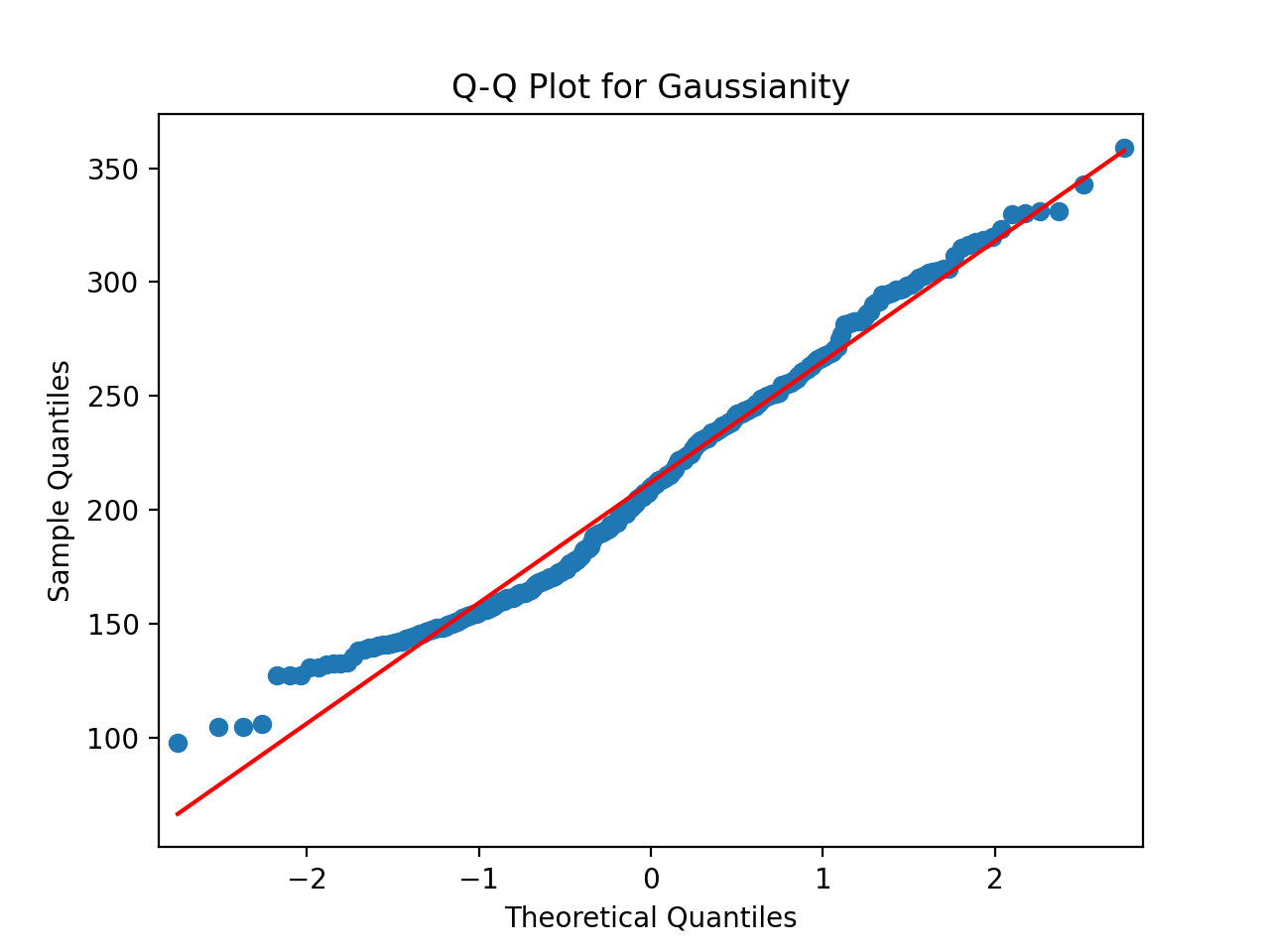}
			\caption{\textbf{Window 12}}
		\end{subfigure}
		
		\caption{\textbf{Local Q-Q plots of the driving function increments $\Delta U_t$ (reconstructed  (thickened)  network).}
			Representative outer spatial windows (2, 3, 4, 7, 8, 12), computed on the largest connected component of the network configuration during the expansion phase. 
			The empirical quantiles closely follow the theoretical Gaussian quantiles, with mild deviations confined to extreme tails.}
		\label{QQPlotsLocauxSLEThickeningBiggestComponent}
	\end{figure*}
	

	
	\begin{figure*}[h!]
		\centering
		
		\begin{subfigure}[t]{0.32\textwidth}
			\centering
			\includegraphics[width=\linewidth]{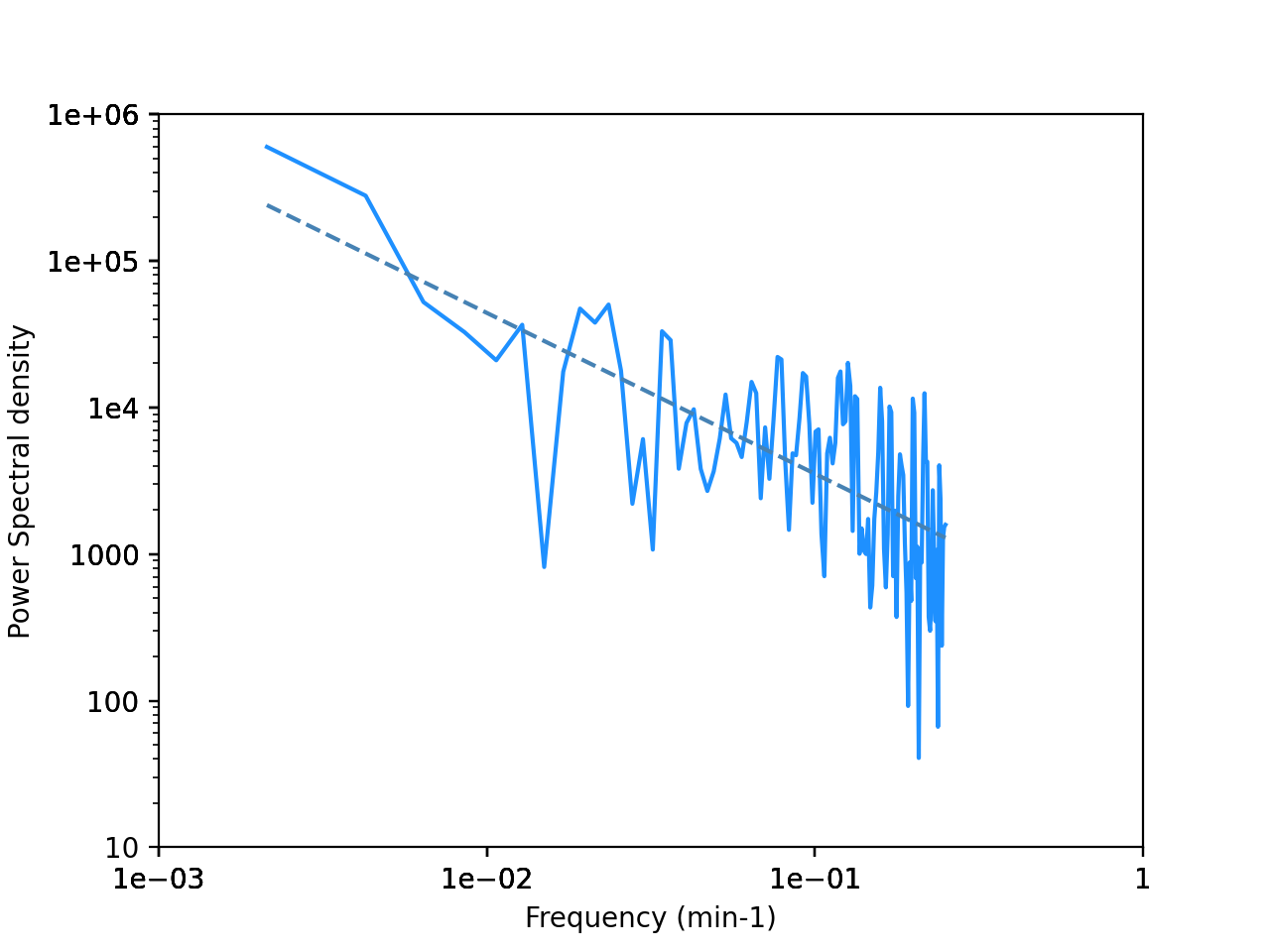}
			\caption{\textbf{Window 2}}
		\end{subfigure}
		\hfill
		\begin{subfigure}[t]{0.32\textwidth}
			\centering
			\includegraphics[width=\linewidth]{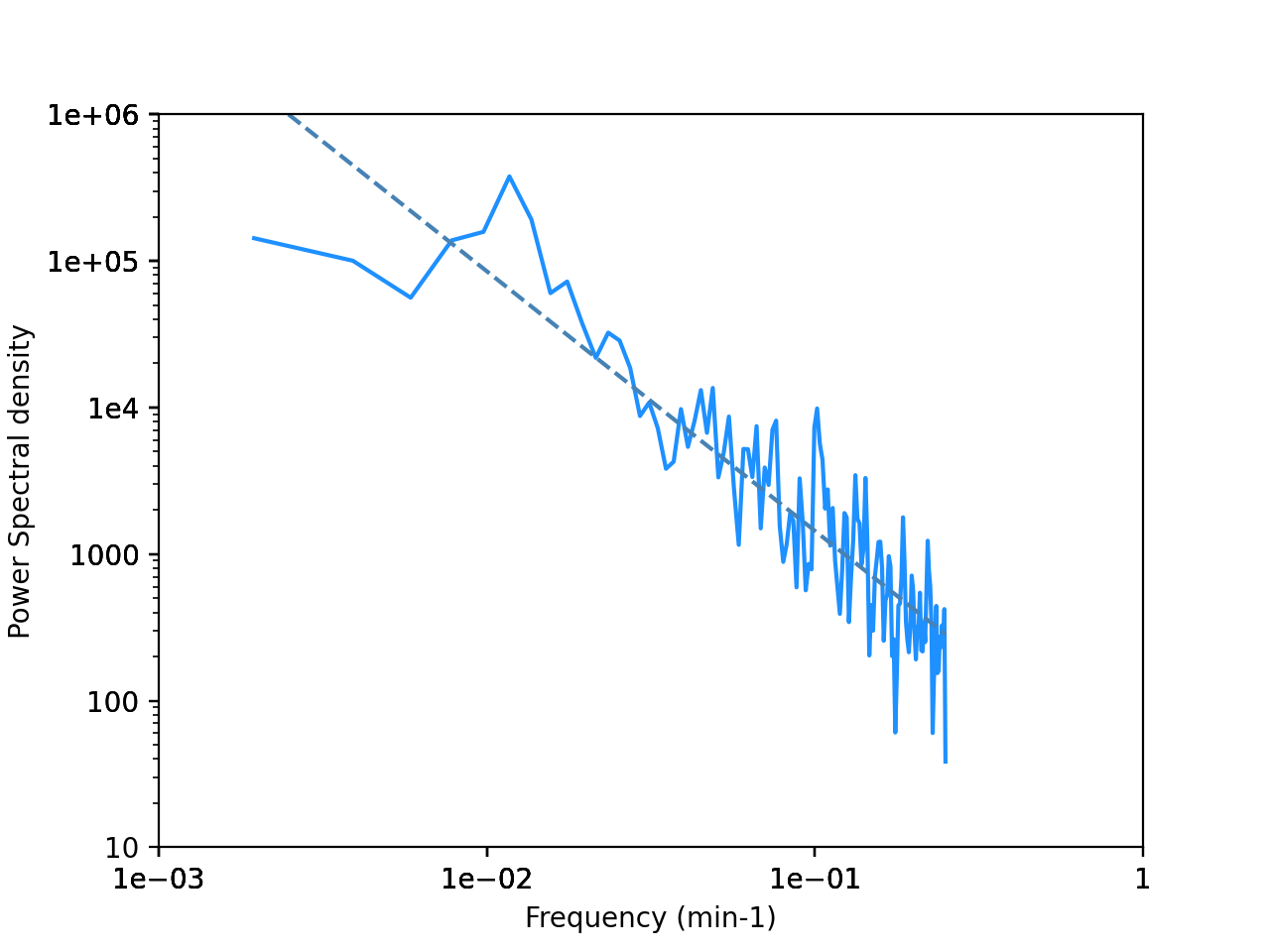}
			\caption{\textbf{Window 3}}
		\end{subfigure}
		\hfill
		\begin{subfigure}[t]{0.32\textwidth}
			\centering
			\includegraphics[width=\linewidth]{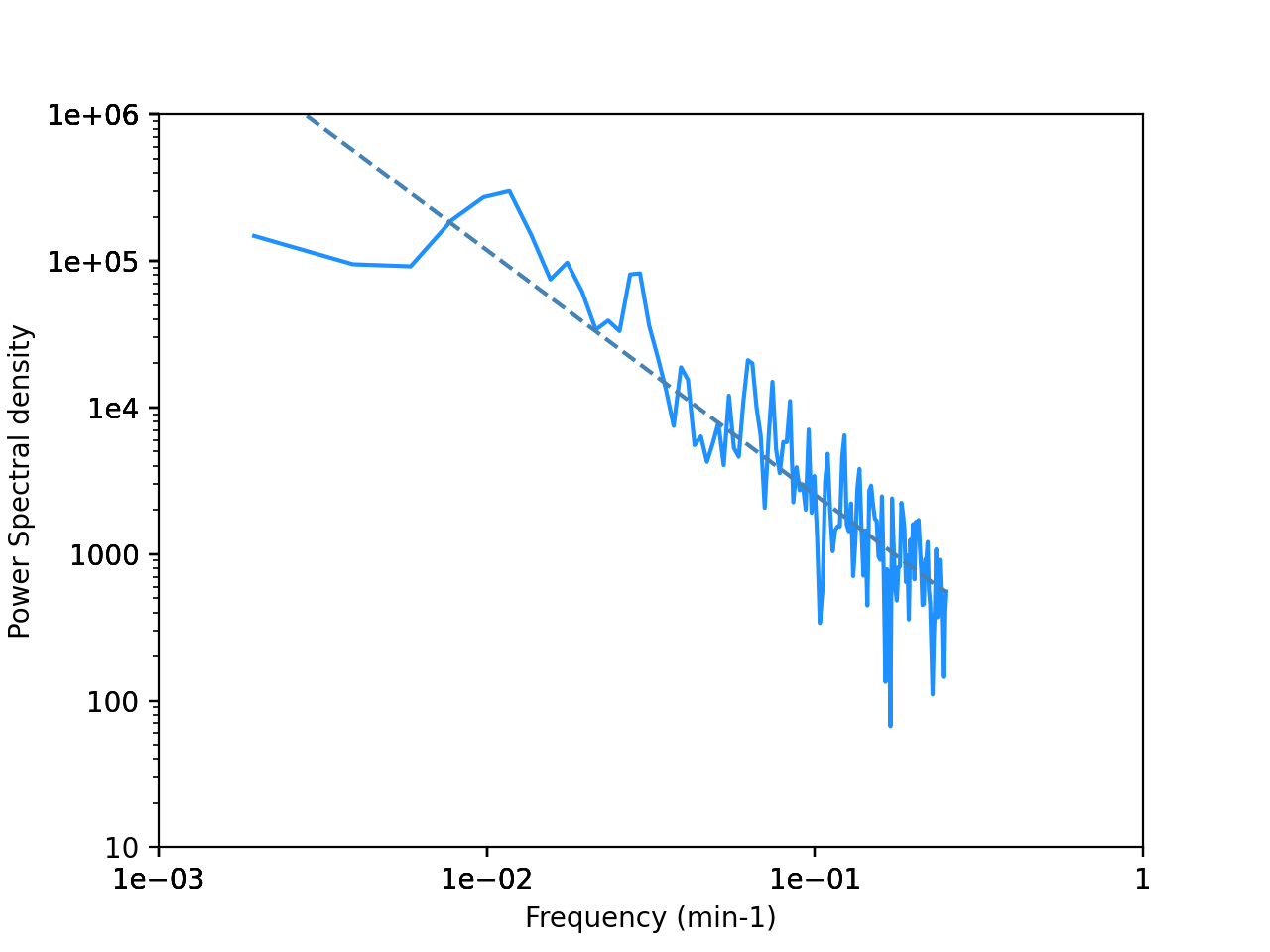}
			\caption{\textbf{Window 4}}
		\end{subfigure}
		
		\vspace{0.35cm}
		
		\begin{subfigure}[t]{0.32\textwidth}
			\centering
			\includegraphics[width=\linewidth]{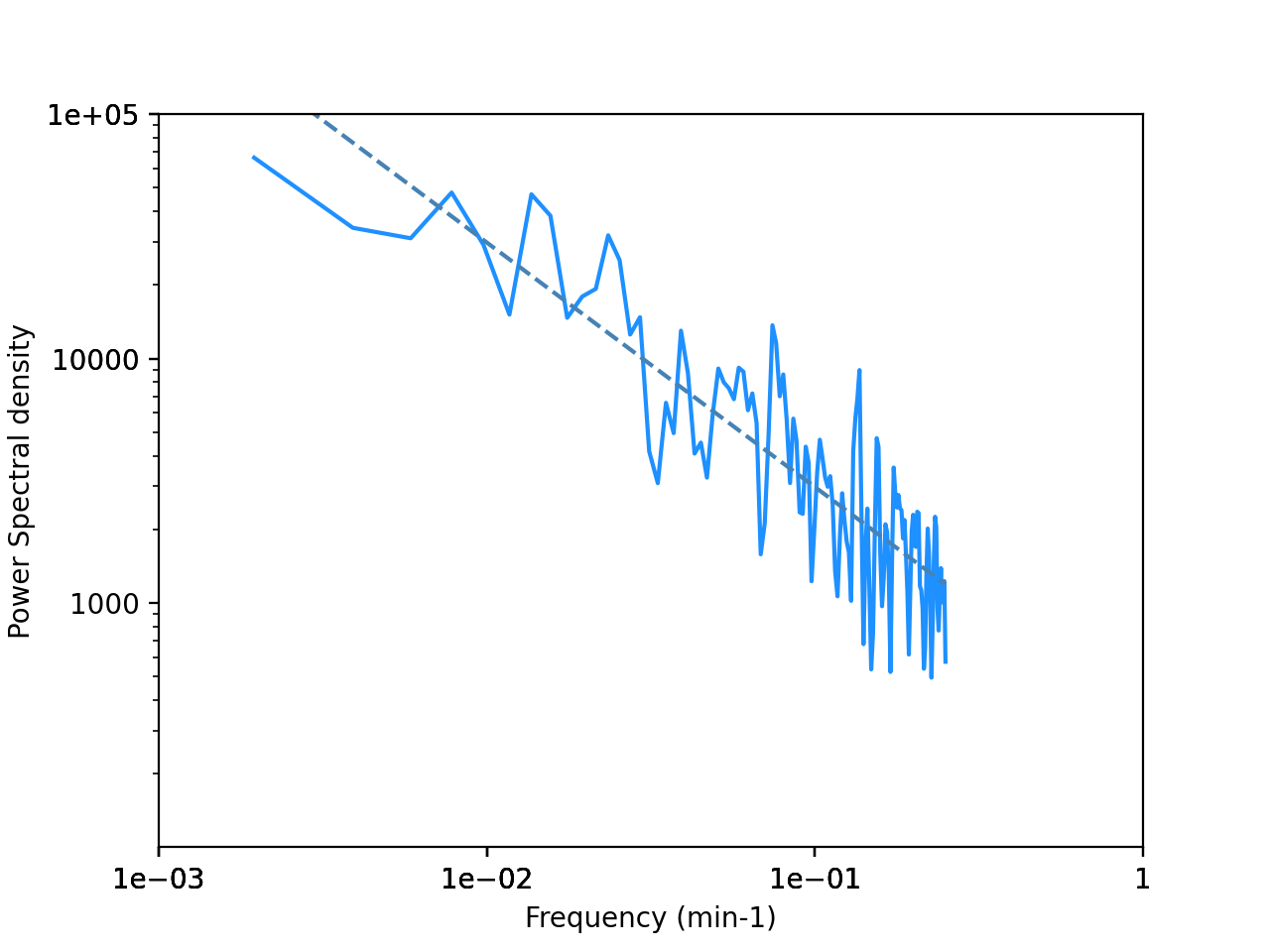}
			\caption{\textbf{Window 7}}
		\end{subfigure}
		\hfill
		\begin{subfigure}[t]{0.32\textwidth}
			\centering
			\includegraphics[width=\linewidth]{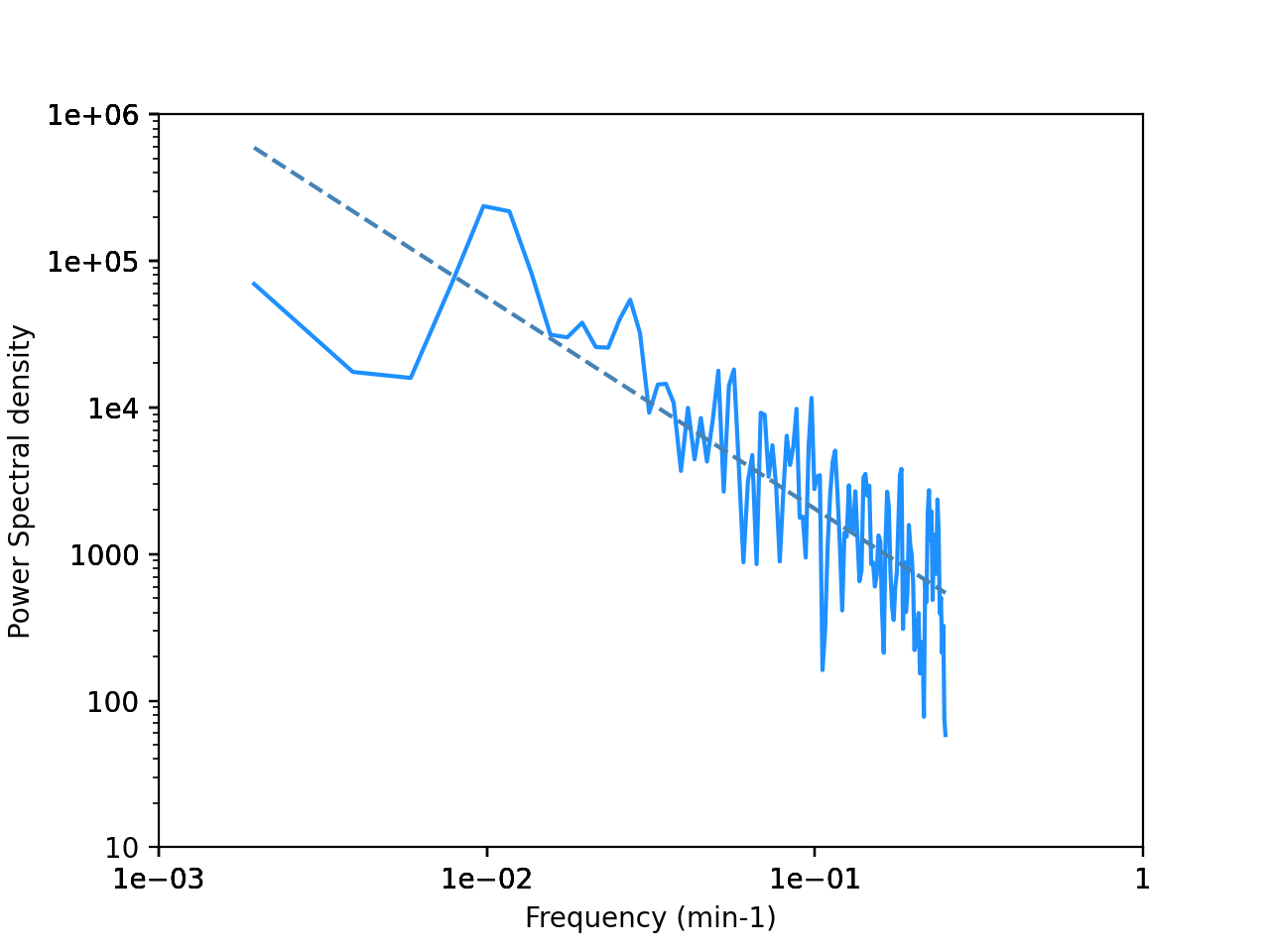}
			\caption{\textbf{Window 8}}
		\end{subfigure}
		\hfill
		\begin{subfigure}[t]{0.32\textwidth}
			\centering
			\includegraphics[width=\linewidth]{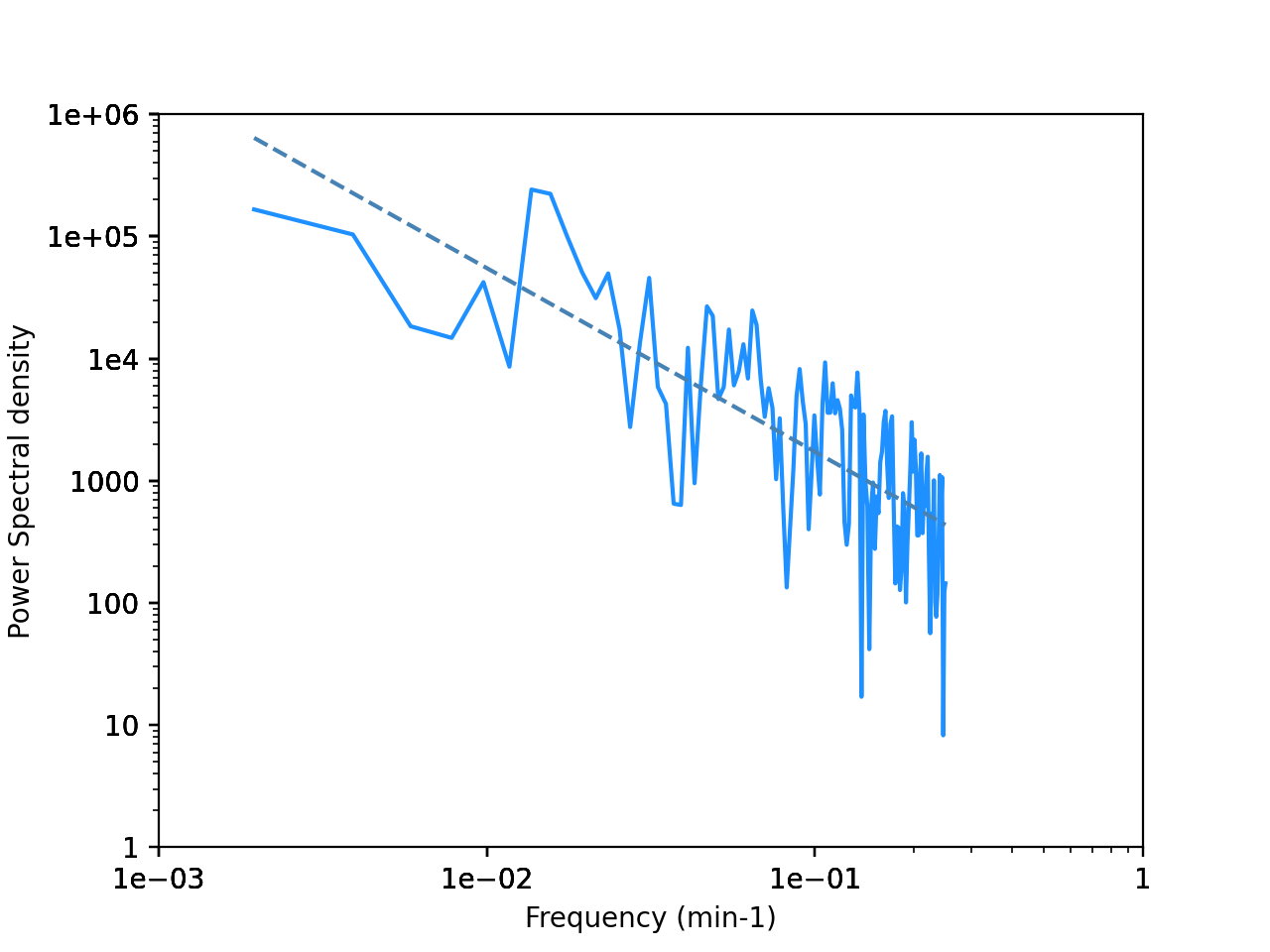}
			\caption{\textbf{Window 12}}
		\end{subfigure}
		
		\caption{\textbf{Local power spectral density (PSD) diagnostics for network configuration during the expansion phase.}
			Log--log PSD plots with linear regression for representative outer windows (2, 3, 4, 7, 8, 12), computed on the largest connected component of the network configuration during the expansion phase. 
			The fitted slopes provide local estimates of the scaling exponent $\beta$ in $S(\omega)\propto\omega^{-\beta}$, remaining broadly compatible with Brownian-type scaling ($\beta\approx 2$), while exhibiting inter-window variability due to structural heterogeneity and finite-size effects.}
		\label{PSDWithRegressionLocauxSLEThickeningBiggestComponent}
	\end{figure*}
	
	

	\begin{figure*}[h!]
		\centering
		
		\begin{subfigure}[t]{0.32\textwidth}
			\centering
			\includegraphics[width=\linewidth]{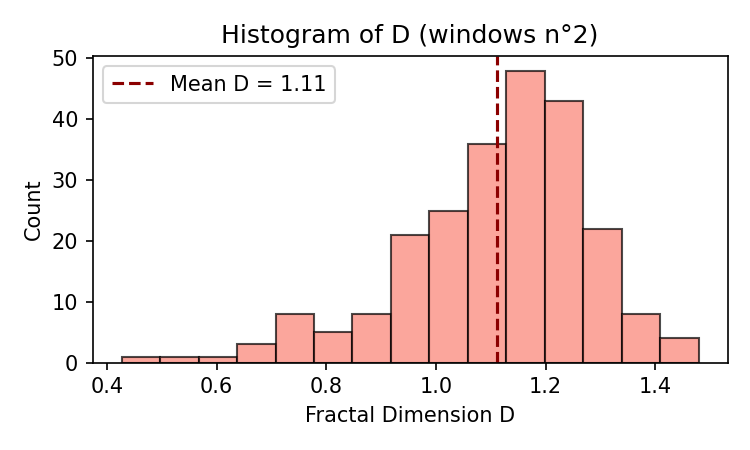}
			\caption{\textbf{Window 2}}
		\end{subfigure}
		\hfill
		\begin{subfigure}[t]{0.32\textwidth}
			\centering
			\includegraphics[width=\linewidth]{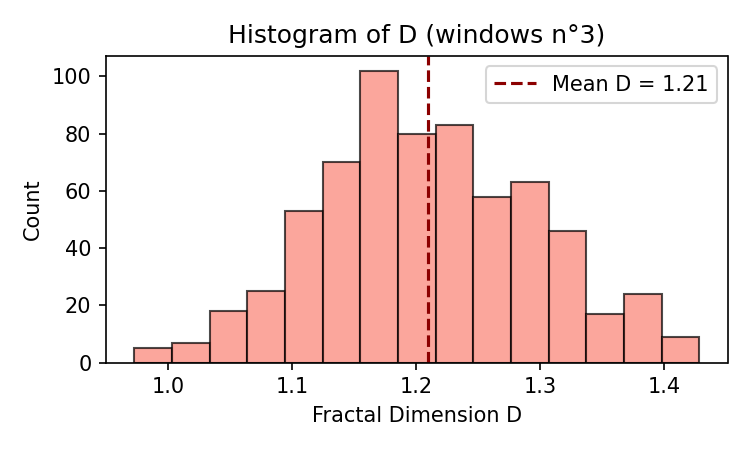}
			\caption{\textbf{Window 3}}
		\end{subfigure}
		\hfill
		\begin{subfigure}[t]{0.32\textwidth}
			\centering
			\includegraphics[width=\linewidth]{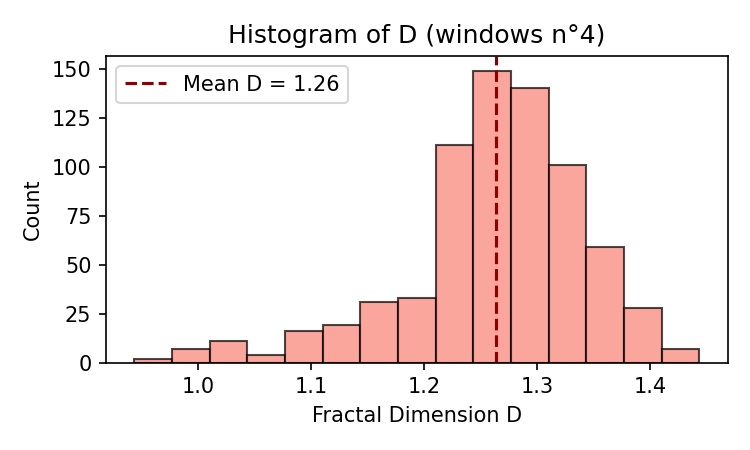}
			\caption{\textbf{Window 4}}
		\end{subfigure}
		
		\vspace{0.35cm}
		
		\begin{subfigure}[t]{0.32\textwidth}
			\centering
			\includegraphics[width=\linewidth]{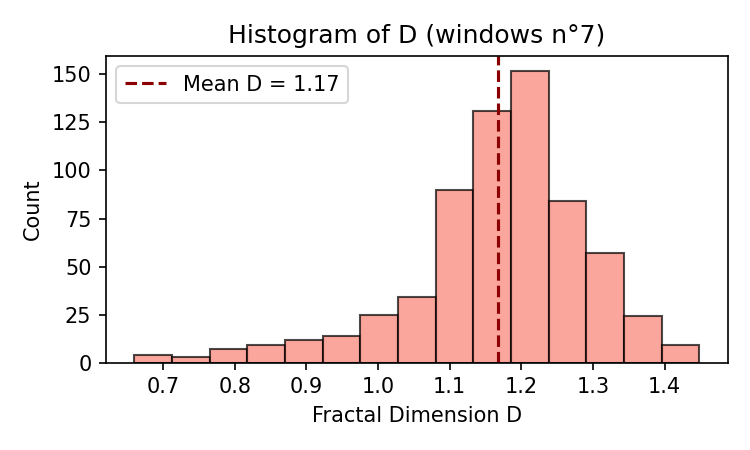}
			\caption{\textbf{Window 7}}
		\end{subfigure}
		\hfill
		\begin{subfigure}[t]{0.32\textwidth}
			\centering
			\includegraphics[width=\linewidth]{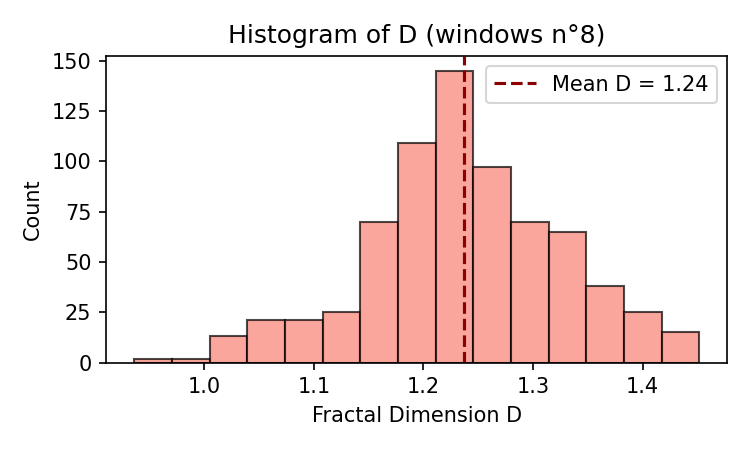}
			\caption{\textbf{Window 8}}
		\end{subfigure}
		\hfill
		\begin{subfigure}[t]{0.32\textwidth}
			\centering
			\includegraphics[width=\linewidth]{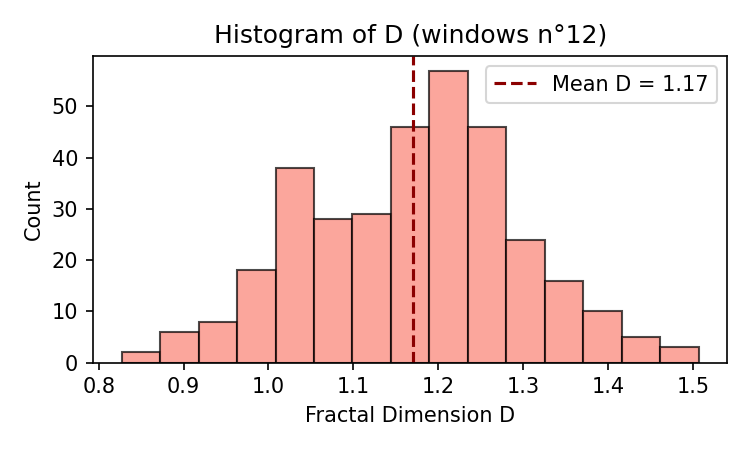}
			\caption{\textbf{Window 12}}
		\end{subfigure}
		
		\caption{\textbf{Local fractal dimension of the network configuration during the expansion phase.}
			Histograms computed on representative outer windows (2, 3, 4, 7, 8, 12) for the largest connected component of the network configuration during the expansion phase. 
			Compared with the network configuration during the retraction phase-level fractal dimensions, the present distributions are shifted toward higher values and exhibit reduced dispersion, reflecting the thickened and spatially reinforced character of the reconstructed  (thickened) structure.}
		\label{HistogrammesLocauxDimensionFractaleSLEThickeningBiggestComponent}
	\end{figure*}
	
	

	\begin{figure*}[h!]
		\centering
		
		\begin{subfigure}[t]{0.32\textwidth}
			\centering
			\includegraphics[width=\linewidth]{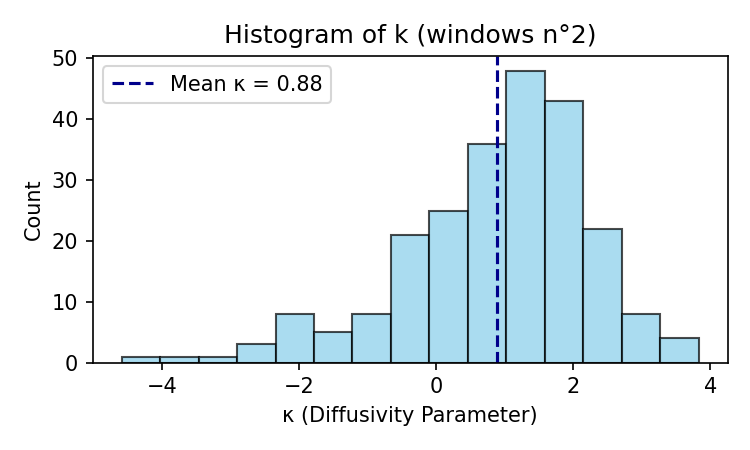}
			\caption{\textbf{Window 2}}
		\end{subfigure}
		\hfill
		\begin{subfigure}[t]{0.32\textwidth}
			\centering
			\includegraphics[width=\linewidth]{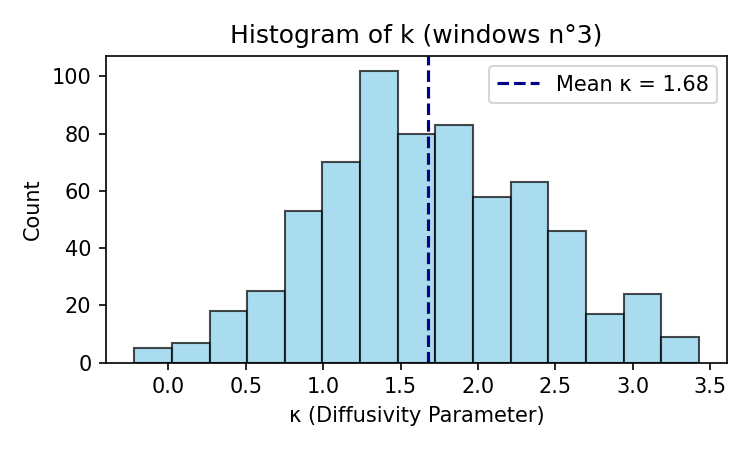}
			\caption{\textbf{Window 3}}
		\end{subfigure}
		\hfill
		\begin{subfigure}[t]{0.32\textwidth}
			\centering
			\includegraphics[width=\linewidth]{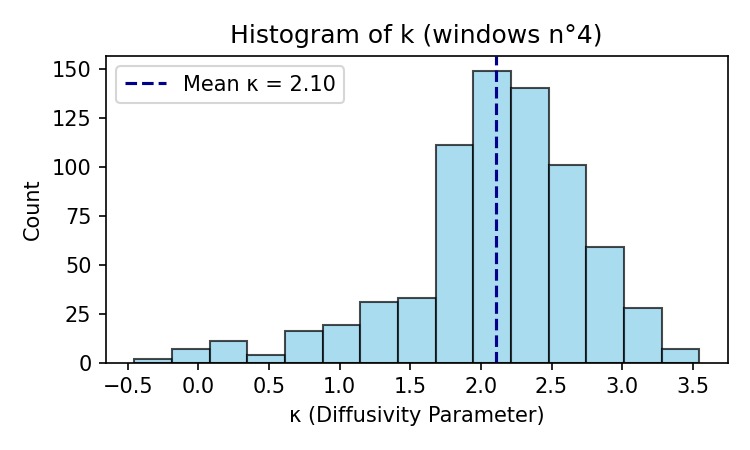}
			\caption{\textbf{Window 4}}
		\end{subfigure}
		
		\vspace{0.35cm}
		
		\begin{subfigure}[t]{0.32\textwidth}
			\centering
			\includegraphics[width=\linewidth]{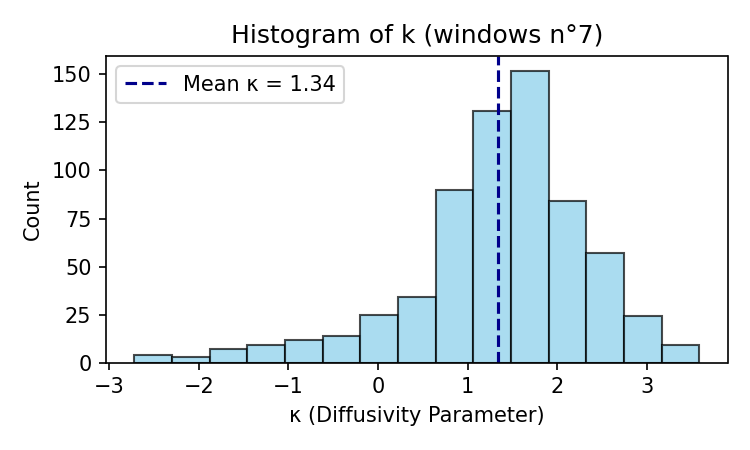}
			\caption{\textbf{Window 7}}
		\end{subfigure}
		\hfill
		\begin{subfigure}[t]{0.32\textwidth}
			\centering
			\includegraphics[width=\linewidth]{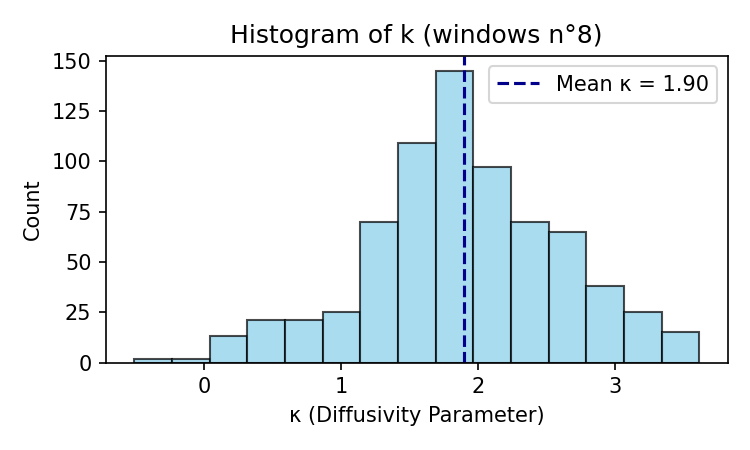}
			\caption{\textbf{Window 8}}
		\end{subfigure}
		\hfill
		\begin{subfigure}[t]{0.32\textwidth}
			\centering
			\includegraphics[width=\linewidth]{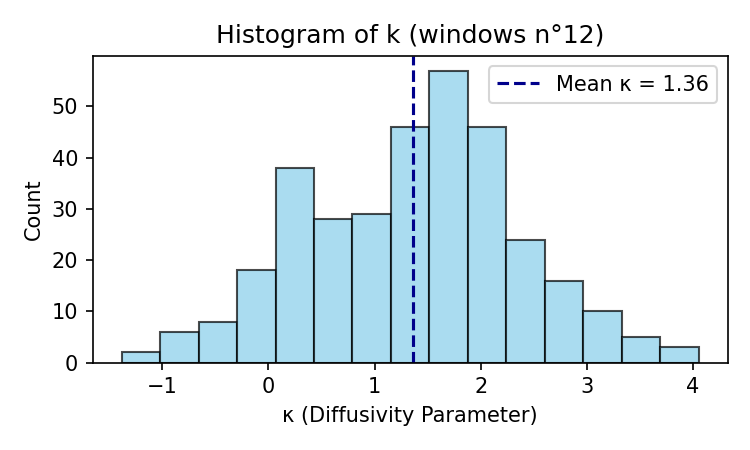}
			\caption{\textbf{Window 12}}
		\end{subfigure}
		
		\caption{\textbf{Local diffusivity parameter $\kappa$ for network configuration during the expansion phase.}
			Histograms estimated on representative outer windows (2, 3, 4, 7, 8, 12) for the largest connected component of the network configuration during the expansion phase. 
			Compared with the pseudopod- and retraction-level estimates, the distributions are shifted toward larger values and exhibit moderate dispersion, reflecting the structurally reinforced and thickened character of the  (expanded) network, while remaining compatible with Brownian-type scaling within the explored resolution range.}
		\label{HistogrammesLocauxKappaSLEThickeningBiggestComponent}
	\end{figure*}
	
	\clearpage
	
	\newpage

	\bibliographystyle{alpha}
	\bibliography{BibliographieClaire}

	\end{document}